# LARGEST SET OF UNITALS IN PROJECTIVE PLANES OF ORDER 16

Stoicho D. Stoichev

Department of Computer Systems, Technical University of Sofia, stoi@tu-sofia.bg

May 3, 2020

## Abstract

In this article we present the largest set of unitals (totally 553) in projective planes of order 16. An open question is what'is the number of the known unitals that are non-isomorphic to the reported ones. The results are obtained with a program that implements author's algorithm.

### *1. Introduction*

*A $t - (v; k; \lambda)$ design D [1] is a set X of points together with a family B of k-subsets of X called blocks with the property that every t points are contained in exactly $\lambda$ blocks. The design with $t = 2$ is called a block–design. The block-design is symmetric if the role of the points and blocks can be changed and the resulting configuration is still a block-design. A projective plane of order n is a symmetric 2-design with $v = n^2 + n + 1$, $k = n + 1$, $\lambda = 1$. The blocks of such a design are called lines. A unital in a projective plane of order $n = q^2$ is a set of $q^3 + 1$ points that meet every line in one or $q + 1$ points. We also call it point unital. The unital is called line unital if we consider the lines as points. In the case of projective planes of order $n = 16$ we have: $q = 4$, the projective plane is $2 - (273; 17; 1)$ design, the unital is a subset of $q^3 + 1 = 4^3 + 1 = 65$ points and every line meets 1 or 5 points from the subset.* In the case of projective planes of order $n = 25$ we have: $q = 5$, the projective plane is $2 - (651; 26; 1)$ design, the unital is a subset of $q^3 + 1 = 5^3 + 1 = 126$ points and every line meets 1 or 6 points from the subset.
Some experimental results of the algorithm in [5] for the known planes of order 16 [4] are given in [6]. The basic tool in [5] is the author's algorithm in [7] that determines the generators, orbits and order of the graph automorphism group, and graph isomorphism.
Recently, new unitals are found in [3] where is stated that the number of all known unitals in the projective planes of order 16 is 108. In this paper we present 553 found unitals:
i) 149 unitals in Gordon Royle's representation [4] of the planes (Table 1, Appendix 1);
ii) 148 unitals in Eric Moorhouse's representation [9] of the planes (Table 1, Appendix 2);
iii) 256 unitals in dreadnaut representation [10] of the planes (Table 3, Appendix 3).
All found unitals in a plane are non-isomorphic pairwise. The number of the new-found unitals is greater than the number of unitals in [6] and [3]. The reported results may be considered as an extension of the results in [5, 6]. Probably, some of the unitals (maybe all of them) found in [3] are isomorphic to a subset of the new unitals presented here. The isomorphism check for each pair of all found unitals is a next task.

## 2. Experimental results

In Appendix 1 we present experimental results from a program that finds the unitals in all 13 known projective planes of order 16 from [4]. The vertex labels in [4] start from 0, but in our program and results the starting label is 1. The program (named M65ALL16pp16DAT2) we use in present paper is based on the algorithm described in [5]. The algorithm is heuristic – it does not guarantees finding all possible unitals of a given plane. All found 149 unitals in Appenix 1 are non-isomorphic each other. This has been checked by our program UNITALSISOM9pp16DRE2. The current results are point unitals for all planes. In the previous two versions of this article the results are point unitals for some planes and line unitals for others.
The results for each plane contain the following information:
Plane names and order of plane automorphism group (two names:our and Gordon Royle's [4] naming)
Unital number (Unital No)
Order of the unital automorphism group (|Aut(G,M65)|)
Labels of the vertices of the unital (Notice: the label of the first vertex is 1)
***Example*** (for the plane G546A (our naming), MATH plane (Gordon Royle naming), see the first part of the Appendix 1 - from the line with 'G546A' to the line with 'G546B'): The order of the plane automorphism group is 12288. Then, the results for 16 unitals follow. The orders of their automorphism groups are 16, 16, 8, 8, 8, 8, 8, 32, 8, 32, 16, 16, 4, 8, 4, 4. For the first unital in this plane the results are:
PROGRAM_NAME=M65ALL16pp16DAT2 (the first line in the Appendix 1)
          g546A, MATH graph in Gordon Royle format, group order 12288
Unitals No    1
 |Aut(G,M65)|=   128
  1  42  43  44  45  46  47  48  49  58  59  60  61  62  63  64
 65  98 100 102 104 107 109 111 113 114 115 116 117 118 119 120
121 194 196 198 200 203 205 207 209 226 228 230 232 235 237 239
241 250 251 252 253 254 255 256 257 258 260 262 264 267 269 271
273

Table 1. Comparison of the numbers of the found unitals in projective planes of order 16 (with the author's algorithm on Gordon Royle's and Moorhouse's presentation and those in [3]) *(\* indicates difference with column two)*

| Plane | Number of the found unitals on Gordon Royle's Presentation of planes | Number of the found unitals on Moorhouse's presentation of planes | Number of the found unitals in [3] |
|---|---|---|---|
| MATH | 16 | 16 | 16 |
| JOHN | 21 | 25 | 8* |
| BBH1 | 16 | 13 | 6* |
| PG(2,16) | 2 | 2 | 2 |
| BBS4 | 13 | 13 | 3* |
| JOWK | 7 | 7 | 5* |
| DSFP | 2 | 2 | 2 |
| HALL | 7 | 6 | 6 |
| DEMP | 4 | 4 | 4 |
| SEMI2 | 21 | 21 | 21 |
| SEMI4 | 12 | 12 | 8* |
| LMRH | 2 | 2 | 2 |
| BBH2 | 26 | 25 | 12* |
| **TOTAL** | **149** | **148** | **95** |

*Table 2. Unitals found with the author's algorithm on G.Royle's presentation of planes*

| Plane | U No | |Aut(U)| | Plane | U No | |Aut(U)| | Plane | U No | |Aut(U)| | Plane | U No | |Aut(U)| |
|---|---|---|---|---|---|---|---|---|---|---|---|
| BBH2 | 1 | 16 | JOHN | 1 | 48 | SEMI2 | 1 | 64 | BBH1 | 1 | 16 |
|  | 2 | 16 |  | 2 | 32 |  | 2 | 64 |  | 2 | 16 |
|  | 3 | 16 |  | 3 | 32 |  | 3 | 64 |  | 3 | 8 |
|  | 4 | 32 |  | 4 | 16 |  | 4 | 192 |  | 4 | 8 |
|  | 5 | 32 |  | 5 | 8 |  | 5 | 48 |  | 5 | 8 |
|  | 6 | 4 |  | 6 | 8 |  | 6 | 192 |  | 6 | 8 |
|  | 7 | 4 |  | 7 | 8 |  | 7 | 48 |  | 7 | 8 |
|  | 8 | 80 |  | 8 | 8 |  | 8 | 192 |  | 8 | 32 |
|  | 9 | 20 |  | 9 | 24 |  | 9 | 48 |  | 9 | 8 |
|  | 10 | 10 |  | 10 | 16 |  | 10 | 64 |  | 10 | 32 |
|  | 11 | 8 |  | 11 | 16 |  | 11 | 16 |  | 11 | 16 |
|  | 12 | 8 |  | 12 | 8 |  | 12 | 64 |  | 12 | 16 |
|  | 13 | 16 |  | 13 | 8 |  | 13 | 64 |  | 13 | 4 |
|  | 14 | 16 |  | 14 | 4 |  | 14 | 64 |  | 14 | 8 |
|  | 15 | 16 |  | 15 | 4 |  | 15 | 64 |  | 15 | 4 |
|  | 16 | 16 |  | 16 | 4 |  | 16 | 64 |  | 16 | 4 |
|  | 17 | 16 |  | 17 | 8 |  | 17 | 64 | HALL | 1 | 1200 |
|  | 18 | 8 |  | 18 | 16 |  | 18 | 32 |  | 2 | 48 |
|  | 19 | 8 |  | 19 | 8 |  | 19 | 32 |  | 3 | 32 |
|  | 20 | 8 |  | 20 | 8 |  | 20 | 64 |  | 4 | 100 |
|  | 21 | 8 |  | 21 | 8 |  | 21 | 64 |  | 5 | 16 |
|  | 22 | 8 | DSFP | 1 | 12 | DEMP | 1 | 48 |  | 6 | 80 |
|  | 23 | 8 |  | 2 | 24 |  | 2 | 16 |  |  |  |
|  | 24 | 8 | LMRH | 1 | 32 |  | 3 | 12 | PG(2,16) | 1 | 249600 |
|  | 25 | 4 |  | 2 | 16 |  | 4 | 24 |  | 2 | 768 |
|  | 26 | 4 |  |  |  |  |  |  |  |  |  |

*Table 2. Unitals found with the author's algorithm on G.Royle's presentation of planes (continued 1)*

| Plane | U No | |Aut(U)| | Plane | U No | |Aut(U)| | Plane | U No | |Aut(U)| | Plane | U No | |Aut(U)| |
|---|---|---|---|---|---|---|---|---|---|---|---|
| MATH | 1 | 16 | BBS4 | 1 | 8 | SEMI4 | 1 | 8 | JOWK | 1 | 8 |
|  | 2 | 128 |  | 2 | 8 |  | 2 | 8 |  | 2 | 16 |
|  | 3 | 64 |  | 3 | 8 |  | 3 | 128 |  | 3 | 4 |
|  | 4 | 128 |  | 4 | 8 |  | 4 | 64 |  | 4 | 48 |
|  | 5 | 16 |  | 5 | 4 |  | 5 | 48 |  | 5 | 24 |
|  | 6 | 16 |  | 6 | 4 |  | 6 | 192 |  | 6 | 32 |
|  | 7 | 64 |  | 7 | 8 |  | 7 | 8 |  | 7 | 12 |
|  | 8 | 128 |  | 8 | 24 |  | 8 | 12 |  |  |  |
|  | 9 | 128 |  | 9 | 12 |  | 9 | 128 |  |  |  |
|  | 10 | 128 |  | 10 | 4 |  | 10 | 128 |  |  |  |
|  | 11 | 128 |  | 11 | 4 |  | 11 | 12 |  |  |  |
|  | 12 | 16 |  | 12 | 4 |  | 12 | 8 |  |  |  |
|  | 13 | 32 |  | 13 | 12 |  |  |  |  |  |  |
|  | 14 | 64 |  |  |  |  |  |  |  |  |  |
|  | 15 | 16 |  |  |  |  |  |  |  |  |  |
|  | 16 | 16 |  |  |  |  |  |  |  |  |  |

In Table 1 we present the numbers of the found unitals in projective planes of order 16 for three cases: i) found with the author's algorithm on Gordon Royle's planes presentation; ii) found with the author's algorithm on Moorhouse's planes presentation; iii) found in [3]). The total numbers of the found unitals are 149, 148 and 95, respectively. For seven planes (MATH, PG(2,16), DSFP, HALL, DEMP, SEMI2, LMRH) the numbers in cases (ii) and (iii) are equal. For other planes (BBH1, BBS4, SEMI4, BBH2) the difference in the numbers is quite large. We suppose that all unitals in [3] are isomorphic to a subset of the new unitals (presented in Appendix 1) but this has not been checked until now. This is based on the fact that we have done this check for a few unitals from [3].

In Table 2 we give the numbers of unitals found with the author's algorithm on G.Royle's presentation of planes. For each plane we give each found unital, its number (U No) and group order (|Aut(U)|). The unitals themselves and their orders for this case are in Appendix 2.

*Table 3 . Unitals found with the author's algorithm on dreadnaut presentation of projective planes of order 16 [10]*

| Projective plane G | |Aut(G)| | Number of the found unitals |
|---|---|---|
| pp-16-1 | 34217164800 | 2 |
| pp-16-2 | 147456 | 21 |
| pp-16-3 | 884736 | 11 |
| pp-16-4 | 921600 | 6 |
| pp-16-5 | 921600 | 6 |
| pp-16-6 | 258048 | 2 |
| pp-16-7 | 258048 | 2 |
| pp-16-8 | 258048 | 7 |
| pp-16-9 | 258048 | 7 |
| pp-16-10 | 55296 | 2 |
| pp-16-11 | 55296 | 2 |
| pp-16-12 | 92160 | 3 |
| pp-16-13 | 92160 | 3 |
| pp-16-14 | 2304 | 27 |
| pp-16-15 | 2304 | 29 |
| pp-16-16 | 3456 | 13 |
| pp-16-17 | 3456 | 13 |
| pp-16-18 | 3840 | 26 |
| pp-16-19 | 3840 | 26 |
| pp-16-20 | 12288 | 16 |
| pp-16-21 | 12288 | 16 |
| pp-16-22 | 18432 | 16 |
| Total | | 256 |

In Table 3 we present the numbers of the found unitals by the author's algorithm on dreadnaut presentation of projective planes of order 16 [10], where for each plane presented is also its dual plane if any. The unitals themselves and their orders for this case are in Appendix 3.

### 3. Concluding remarks

By our algorithm we have found new unitals in projective plane of order 16. These results are obtained by the author's heuristic algorithm, i.e. more results can be found by some improved algorithm. An open question is what'is the number of the known unitals that are isomorphic to the reported ones.The following approaches can be used to find more or all unitals: (a) development of improved algorithms by finding new conditions for pruning the search tree; (b) transformation of the solution for one plane to solution for another plane (R. Mathon's approach – in private communication); (c) development of parallel algorithms; (d) use of different representation of the planes; (e) use of relabeling of unital's points and lines.

## Acknowledgements

The author would like to thank Vladimir Tonchev for suggesting the problem of developing an algorithm for finding unitals in projective planes and for giving the general idea of such an algorithm – use of unions of combinations of orbits of different subgroup of the plane automorphism group, and for extensive discussions and exchanges for many years. Thanks to Mustafa Gezek for the intensive communication and for checking the correctness of the current results.

**Appendix 1: Unitals of projective planes of order 16, found with author's algorithm on G. Royle's representation of the planes [4]**

```
  G546A, MATH graph in G. Royle's format, |AutG|=12288
Unitals No       1
 |Aut(G,M65)|=   16
  4   7   9  12  14  19  28  29  31  32  33  38  39  41  42  62
 67  69  76  77  79  81  83  88  90  93 103 113 115 118 119 121
132 145 148 150 152 158 162 168 170 172 175 179 182 184 186 189
193 196 198 205 206 217 228 231 233 238 239 248 250 251 252 255
268
Unitals No       2
 |Aut(G,M65)|=  128
  6   8   9  12  16  18  19  22  30  32  33  36  38  41  42  51
 52  58  59  64  67  68  70  72  77  87  89  90  94  96 103 125
132 154 161 163 167 168 170 180 185 187 189 190 193 198 199 205
206 215 216 219 221 224 225 227 235 238 239 241 245 248 249 251
```

```
 268
Unitals No        3
 |Aut(G,M65)|=  64
  5   7   9  12  14  18  19  24  29  31  33  36  40  41  42  51
 52  58  59  62  67  68  69  71  79  85  89  90  93  95 101 127
130 156 162 164 167 168 172 178 186 188 189 190 195 198 199 205
206 215 216 217 221 224 226 228 233 238 239 243 245 248 250 252
 268
Unitals No        4
 |Aut(G,M65)|=  128
 12  18  36  37  39  46  47  53  56  58  61  63  67  68  74  78
 80  84  86  88  89  90  98 105 106 110 111 115 116 117 120 124
133 137 140 141 142 146 147 151 152 159 167 168 169 171 173 177
179 183 189 190 194 196 199 201 204 210 211 218 220 221 239 245
 268
Unitals No        5
 |Aut(G,M65)|=  16
  1   2   5  12  16  18  22  27  28  31  36  38  40  47  48  53
 54  58  62  64  65  68  74  75  79  81  84  85  90  91 112 118
130 132 133 136 140 146 154 156 158 159 161 162 165 168 170 180
187 188 190 191 198 224 225 231 232 233 238 243 248 251 253 254
 268
Unitals No        6
 |Aut(G,M65)|=  16
  6   7   8  11  12  17  18  29  30  32  34  38  40  43  46  49
 56  60  62  64  65  72  73  76  79  82  83  85  91  94  97 101
102 103 109 119 123 125 127 128 143 149 162 164 165 167 176 182
186 188 189 191 194 220 229 230 232 235 237 241 247 254 255 256
 268
Unitals No        7
 |Aut(G,M65)|=  64
  1   2   3   9  13  19  23  25  27  28  33  35  37  40  48  54
 57  59  62  63  67  69  70  76  79  82  85  89  95  96 112 118
139 145 161 162 168 172 175 178 181 187 188 190 194 195 198 200
203 209 217 220 222 224 229 230 232 234 238 244 248 254 255 256
 268
Unitals No        8
 |Aut(G,M65)|=  128
  3  25  37  40  43  46  48  49  54  56  62  63  69  73  75  76
 79  81  82  83  85  95  97 101 105 106 112 115 116 118 123 127
130 134 139 141 142 145 151 152 156 160 162 172 174 175 176 178
181 182 184 188 194 195 201 203 208 209 211 214 217 220 232 254
 268
Unitals No        9
 |Aut(G,M65)|=  128
  3   8   9  11  15  17  19  21  25  30  33  59  65  66  69  72
 75  81  91  92  94  95 101 102 104 108 109 114 119 126 127 128
129 130 131 138 143 148 149 153 155 156 162 163 165 166 176 182
185 188 191 192 198 224 227 232 236 238 240 242 246 248 249 254
 268
Unitals No       10
 |Aut(G,M65)|=  128
  3   4   5   7  15  21  25  26  29  31  37  63  69  71  73  76
 78  82  83  88  93  95  97 100 104 105 106 115 116 122 123 126
131 134 135 141 142 151 152 153 157 160 162 164 169 174 175 179
181 184 186 188 194 220 226 228 231 232 236 242 250 252 253 254
 268
Unitals No       11
 |Aut(G,M65)|=  128
  6  18  22  25  30  31  37  38  41  43  48  52  58  59  60  62
 69  81  85  90  93  96 101 102 106 108 111 115 121 123 124 125
129 130 140 141 143 152 154 158 159 160 162 178 182 186 187 189
193 194 203 206 208 215 217 221 223 224 225 241 245 249 252 254
 268
Unitals No       12
 |Aut(G,M65)|=  16
 15  22  28  29  30  32  33  36  43  45  46  49  52  56  58  60
 68  73  75  77  78  85  87  91  92  94 103 113 118 122 123 127
132 145 150 156 157 160 167 169 170 174 176 178 180 185 191 192
198 199 201 202 208 211 214 215 221 223 236 241 249 250 251 255
 268
Unitals No       13
 |Aut(G,M65)|=  32
  3   4  11  13  15  19  22  24  30  31  33  51  52  54  55  57
 73  81  91  92  94  95  99 101 103 107 108 118 119 123 126 128
130 132 136 143 144 145 148 153 155 160 174 182 185 188 191 192
```

```
       198 209 212 215 216 222 231 232 234 236 240 241 243 248 249 252
       268
 Unitals No         14
  |Aut(G,M65)|=   64
  15  18  22  24  28  29  33  36  38  43  44  52  53  55  56  60
  80  81  85  87  91  94  98  99 101 107 108 115 118 119 120 123
 130 133 136 143 144 145 147 148 152 160 171 178 180 182 185 192
 193 198 199 207 208 210 211 212 215 223 236 241 243 245 250 255
 268
 Unitals No         15
  |Aut(G,M65)|=   16
   1   6   7  10  15  17  19  25  26  31  37  38  39  40  42  51
  55  56  58  59  79  81  87  88  91  93 100 107 108 111 112 123
 135 139 140 143 144 160 172 179 180 182 186 192 193 194 195 196
 205 211 212 216 221 224 225 228 230 236 237 246 248 252 253 254
 268
 Unitals No         16
  |Aut(G,M65)|=   16
   8  11  13  14  16  18  21  25  28  32  33  34  42  46  48  49
  71  76  77  78  79  81  86  90  91  95  97  98 105 109 111 114
 133 134 138 140 142 149 164 169 170 172 175 177 182 188 189 192
 197 198 201 203 205 214 227 233 234 235 240 242 245 251 254 255
 268
    G546B,JOHN graph in G. Royle's format, |AutG|=2304
   Unitals No          1
 |Aut(G,M65)|=   48
   2   5  12  15  17  18  21  22  27  28  31  32  60  81  84  86
  87  90  91  93  96 113 116 118 119 122 123 125 128 132 135 138
 141 146 149 156 159 162 165 172 175 180 183 186 189 210 212 213
 215 218 220 221 223 225 227 230 232 233 235 238 240 260 263 268
 270
   Unitals No          2
 |Aut(G,M65)|=   32
   2  20  35  49  63  67  70  78  84  86  87  89  90  91  93  98
 100 111 113 118 119 122 123 127 128 132 137 141 142 143 146 151
 155 159 160 161 163 172 175 176 179 186 187 189 190 202 205 208
 210 212 213 215 217 218 221 225 227 230 231 232 233 238 242 249
 254
   Unitals No          3
 |Aut(G,M65)|=   32
   3   4   7  12  14  22  26  27  28  31  33  37  38  41  44  52
  57  59  61  62  63  71  74  80  83  84  85  87  90  92  93  97
 106 112 113 115 117 118 123 126 128 133 157 169 177 198 201 206
 211 212 215 218 219 221 224 226 227 229 232 233 237 238 244 247
 252
   Unitals No          4
 |Aut(G,M65)|=   16
   8   9  10  12  13  20  22  27  28  31  37  40  43  44  48  50
  55  57  58  62  64  74  77  82  84  85  87  90  93  97 102 113
 118 121 123 126 128 134 135 142 152 155 156 161 165 167 180 181
 190 200 201 209 212 215 218 221 224 227 229 232 233 236 238 244
 253
   Unitals No          5
 |Aut(G,M65)|=    8
   2  26  46  50  54  65  66  69  70  75  80  83  88  89  90  93
  94  99 104 123 128 133 135 137 139 140 151 153 154 156 160 161
 162 163 170 174 177 183 188 190 192 195 197 200 204 209 210 213
 214 215 218 220 223 226 229 232 233 243 244 246 247 248 249 251
 254
   Unitals No          6
 |Aut(G,M65)|=    8
   3   4   5   9  13  21  22  28  29  32  34  38  41  43  44  49
  50  51  55  61  62  65  70  71  75  77  80  82  83  88  89  92
  94  97 107 120 126 133 134 137 140 141 146 151 152 154 160 162
 163 167 171 174 177 181 183 184 192 210 213 220 223 243 248 249
 254
   Unitals No          7
 |Aut(G,M65)|=    8
   2   3   9  17  26  27  34  44  46  50  52  54  58  70  71  74
  80  82  83  89  95  99 102 103 104 105 107 109 110 113 114 115
 118 123 124 126 128 135 137 139 153 156 160 161 162 170 183 188
 190 193 208 211 222 231 234 245 252 257 260 262 263 265 267 271
 272
   Unitals No          8
 |Aut(G,M65)|=    8
   7   8  10  14  15  17  18  23  27  31  33  35  37  47  48  50
```

```
  52  56  57  58  59  66  68  69  74  85  90  93  95  99 102 104
 112 115 121 123 128 129 132 137 147 159 160 162 164 166 178 183
 185 195 197 200 204 209 214 215 218 226 229 232 233 244 246 247
 251
  Unitals No       9
|Aut(G,M65)|=  24
   4   7  10  13  49  52  54  55  58  59  60  61  64  83  84  87
  88  89  90  93  94 113 114 117 118 123 124 127 128 162 164 165
 167 170 172 173 175 178 180 181 183 186 188 189 191 193 198 203
 208 210 213 220 223 225 230 235 240 242 245 252 255 260 263 268
 270
  Unitals No      10
|Aut(G,M65)|=  16
   3   8   9  14  20  23  26  29  49  52  54  55  58  59  61  63
  64  99 104 105 110 130 132 133 135 138 140 141 143 146 149 156
 159 163 164 167 168 169 170 173 174 194 195 197 200 201 204 206
 207 209 210 213 214 219 220 223 224 225 230 235 240 258 266 270
 273
  Unitals No      11
|Aut(G,M65)|=  16
  13  31  48  62  67  69  72  73  77  78  79  81  84  85  86  90
  91  96  97  99 105 118 121 128 131 132 136 137 138 139 142 149
 158 159 161 164 171 180 181 183 184 186 189 190 194 197 198 202
 205 217 218 220 221 223 225 227 228 230 232 243 248 251 255 256
 259
  Unitals No      12
|Aut(G,M65)|=   8
   3   5  13  22  28  29  38  41  44  49  50  51  53  60  61  63
  65  67  69  70  71  75  80  90  91  92  97  99 111 119 126 128
 140 142 144 151 155 158 165 173 174 177 180 191 200 203 207 210
 211 213 214 220 221 223 229 231 238 247 255 256 259 262 265 268
 273
  Unitals No      13
|Aut(G,M65)|=   8
   2   4   8  17  21  26  43  46  47  50  53  54  58  60  62  63
  65  67  69  70  74  75  80  85  90  91  99 111 112 115 119 128
 129 142 143 147 148 155 164 165 169 178 180 182 195 198 207 209
 210 213 220 221 222 223 226 234 238 242 244 256 259 262 265 268
 273
  Unitals No      14
|Aut(G,M65)|=   4
  17  22  27  32  35  40  41  46  50  53  60  63  67  70  74  79
  82  88  90  96 101 102 103 104 113 117 122 126 135 136 139 140
 150 151 158 159 164 165 171 174 177 179 181 183 199 200 207 208
 209 211 221 223 234 235 238 239 243 245 250 256 257 259 261 272
 273
  Unitals No      15
|Aut(G,M65)|=   4
  34  37  44  47  49  50  53  54  59  60  63  64  68  69  70  78
  81  84  88  92 103 108 110 112 117 119 120 123 129 135 136 140
 150 151 152 159 165 170 171 174 177 179 183 191 196 198 206 207
 217 219 220 221 228 229 238 240 244 248 255 256 257 259 261 272
 273
  Unitals No      16
|Aut(G,M65)|=   4
  17  22  27  32  35  40  41  46  50  53  60  63  68  69  70  78
  81  84  88  92 103 108 110 112 117 119 120 123 131 133 135 139
 145 154 158 159 162 163 164 171 179 180 181 192 193 195 196 197
 214 215 220 222 231 232 236 240 246 248 252 253 257 259 261 272
 273
  Unitals No      17
|Aut(G,M65)|=   8
   3   7   8   9  14  18  20  23  26  29  48  49  52  53  54  55
  57  58  59  61  64  66  72  73  79  86  87  90  91  98 104 105
 111 113 116 125 128 129 134 135 143 147 148 152 159 169 170 171
 173 179 182 188 191 194 203 208 211 212 216 227 234 237 241 242
 245
  Unitals No      18
|Aut(G,M65)|=  16
  10  12  15  26  28  29  35  40  43  49  54  57  69  74  79  82
  83  87  88  89  93  94 102 104 110 114 117 118 124 126 127 128
 131 137 144 145 153 155 165 167 175 178 183 189 193 195 200 219
 222 224 226 229 231 250 252 253 257 258 259 261 263 267 269 270
 271
  Unitals No      19
|Aut(G,M65)|=   8
```

```
  7   8   9  10  21  22  27  28  33  34  47  48  51  52  53  61
 62  66  71  74  79  84  85  92  93  98  99 101 102 107 108 110
111 113 120 121 128 131 136 137 142 194 195 206 207 214 215 218
219 228 229 231 232 233 234 236 237 241 244 253 256 258 268 271
272
  Unitals No      20
|Aut(G,M65)|=  8
 20  23  26  29  33  38  43  48  50  52  53  55  58  60  61  63
 65  67  70  72  73  75  78  80 113 115 118 120 121 123 126 128
163 164 167 168 169 170 173 174 179 180 183 184 185 186 189 190
209 212 214 215 218 219 221 224 242 243 245 248 249 252 254 255
259
  Unitals No      21
|Aut(G,M65)|=  8
  3  10  12  13  14  22  27  28  29  31  33  34  35  38  47  51
 54  55  56  58  67  81  84  87  90  93  98 113 116 118 123 128
132 135 156 159 164 167 170 171 173 176 178 179 181 184 188 191
199 200 202 204 210 212 219 223 225 226 232 240 248 249 253 256
259
G546CNOV, BBH1 graph in G. Royle's format, |AutG|=18432
  Unitals No       1
|Aut(G,M65)|= 16
 16  23  24  29  31  32  33  39  40  44  47  51  53  57  58  61
 65  85  92  93  94  95  97 104 105 108 111 116 131 133 139 142
144 147 148 150 151 159 164 168 169 170 176 180 182 186 189 190
193 200 201 202 208 211 225 228 230 234 235 242 243 244 245 246
258
  Unitals No       2
|Aut(G,M65)|= 16
 13  17  18  26  28  29  38  39  41  42  47  59  60  61  63  64
 73  82  85  86  91  96  97  98 107 109 111 121 129 133 134 135
137 146 150 154 156 158 165 166 168 173 174 183 185 187 190 192
197 198 199 200 201 222 225 229 230 232 240 241 242 249 250 256
258
  Unitals No       3
|Aut(G,M65)|=  8
  4   7   9  10  18  23  25  26  33  39  41  47  52  54  63  64
 70  73  77  79  91  92  94  95  99 106 108 112 122 123 124 125
132 139 140 144 150 154 157 159 168 169 172 176 180 187 191 192
195 196 205 206 229 230 232 239 247 248 251 254 257 270 271 272
273
  Unitals No       4
|Aut(G,M65)|=  8
  1   2   3   6   8  17  18  20  23  24  33  35  36  45  46  57
 71  72  76  78  80  83  85  89  92  95 100 102 107 108 109 114
118 119 120 123 141 158 162 163 172 174 176 177 181 185 191 192
193 195 202 203 207 217 219 221 222 223 228 229 231 237 240 242
262
  Unitals No       5
|Aut(G,M65)|=  8
  4   7   9  10  18  23  27  28  36  42  45  47  50  52  56  57
 82  92  93  95 106 107 109 110 114 118 119 120 134 137 139 141
149 151 158 159 161 172 173 176 179 183 191 192 196 198 203 207
217 219 221 223 231 235 239 240 242 244 254 256 257 262 263 264
265
  Unitals No       6
|Aut(G,M65)|=  8
  4   7   9  10  18  20  21  22  34  36  45  46  50  56  61  64
 66  68  72  75  87  89  95  96 102 107 109 111 114 118 119 120
134 137 143 144 149 150 151 152 163 164 170 176 177 182 188 191
199 202 203 207 217 219 221 223 231 237 238 240 244 246 247 256
257
  Unitals No       7
|Aut(G,M65)|=  8
 12  17  18  23  26  28  33  36  39  41  43  49  50  53  54  63
 72  82  83  91  92  95  99 105 108 109 112 122 130 131 132 134
140 145 148 151 154 157 161 169 172 174 176 179 180 181 191 192
195 201 205 206 207 220 226 229 230 232 240 243 248 250 251 256
263
  Unitals No       8
|Aut(G,M65)|= 32
  5  11  13  15  21  22  28  29  34  35  38  45  53  55  57  63
 71  76  78  80  86  89  91  95 100 104 107 111 122 123 124 125
133 136 138 139 146 153 156 159 162 165 167 173 178 181 187 192
196 204 205 208 218 220 222 224 225 232 234 239 242 250 251 252
257
```

Unitals No        9
|Aut(G,M65)|=   8
 18  20  28  29  33  35  36  41  52  59  62  64  66  68  72  75
 82  89  90  95  97  98 101 107 114 118 119 120 132 140 141 143
152 154 157 158 165 167 168 174 181 184 186 187 193 195 205 207
209 211 213 215 225 227 231 236 247 253 254 255 257 266 267 268
269
  Unitals No       10
|Aut(G,M65)|=  32
  4   7   9  10  19  24  27  31  34  35  38  45  53  55  57  63
 66  68  72  75  81  84  85  93  98 103 106 112 114 118 119 120
129 134 140 143 147 151 152 157 161 164 166 170 179 182 188 189
193 199 200 203 217 219 221 223 225 232 234 239 242 250 251 252
257
  Unitals No       11
|Aut(G,M65)|=  16
  6  10  14  15  33  35  43  44  50  56  62  63  67  71  75  79
 83  88  92  94  97 101 103 111 113 118 123 127 131 137 139 140
149 152 156 158 163 169 171 176 183 184 185 186 198 200 202 204
209 216 217 224 230 231 234 240 241 244 245 252 257 262 263 264
265
  Unitals No       12
|Aut(G,M65)|=  16
  8   9  12  13  34  42  46  47  51  55  59  61  70  71  74  75
 82  87  90  96  98 101 105 111 115 119 125 126 131 137 138 143
148 152 154 156 162 165 167 173 177 180 190 191 197 199 205 206
211 216 217 222 229 234 238 240 242 246 247 254 257 262 263 264
265
  Unitals No       13
|Aut(G,M65)|=   4
  3   6  12  13  16  25  38  51  53  58  59  63  65  66  68  77
 79  82  83  85  92  93 101 103 105 108 112 114 115 116 120 125
141 153 167 168 171 173 174 185 195 197 202 204 206 211 212 215
217 223 226 230 232 237 239 242 249 251 254 256 258 259 267 269
271
  Unitals No       14
|Aut(G,M65)|=   8
  6   7   8  10  16  19  24  26  27  31  36  49  53  55  58  59
 65  67  70  72  77  81  84  85  91  93 109 113 115 121 125 126
132 136 137 141 144 149 150 154 156 158 161 164 166 169 170 179
180 182 188 189 207 210 214 219 221 224 240 243 249 250 252 255
265
  Unitals No       15
|Aut(G,M65)|=   4
  2   7  15  16  17  23  25  32  34  36  39  44  53  54  58  59
 70  72  74  78  81  84  85  93  98 101 104 109 113 114 125 127
129 136 141 144 149 152 154 156 179 182 188 189 193 195 196 207
210 211 223 224 231 232 235 240 243 249 252 256 257 262 263 264
265
  Unitals No       16
|Aut(G,M65)|=   4
 11  20  22  24  27  30  33  34  38  39  47  49  53  56  57  58
 70  83  87  89  90  94  97  98 105 107 112 121 136 138 140 141
143 145 147 149 155 160 162 167 169 171 175 178 179 181 186 189
198 199 201 203 206 218 227 229 232 239 240 245 246 253 254 255
260
G546D, PG(2,16) graph in G. Royle's format, |AutG|=34217164800
  Unitals No        1
 |Aut(G,M65)|= 249600
  7   9  11  16  17  21  22  23  26  33  36  41  42  47  49  52
 68  71  74  77  80  86  89  90  91  94 105 116 117 118 121 128
131 134 135 136 139 155 167 169 171 176 177 179 181 183 188 189
195 211 213 218 219 224 226 228 229 237 241 253 259 262 269 272
273
  Unitals No        2
 |Aut(G,M65)|= 768
  2   5   6  10  12  24  27  28  30  31  34  36  38  41  47  54
 66  67  68  73  80  86  91  93  94  95 103 104 106 110 113 116
121 124 127 129 135 136 137 139 144 156 162 167 170 171 174 186
199 211 216 218 220 221 227 228 237 240 241 242 259 262 269 272
273
   G546E, BBS4 graph in G. Royle's format, |AutG|=3456
  BRRESH=        1
 |Aut(G,M65)|=   8
  2   3  12  14  17  20  26  32  37  40  44  46  54  55  58  63
 64  66  74  79  85  87  90  99 110 112 115 118 123 129 136 139

143 150 151 152 158 163 165 167 173 178 183 187 188 209 210 213
214 219 220 223 224 226 229 236 239 241 243 246 248 249 251 254
256
 BRRESH=       2
 |Aut(G,M65)|=  8
  3  12  14  15  17  20  23  32  35  37  40  44  55  58  59  60
 64  76  78  79  84  87  91  99 101 104 116 123 128 129 132 144
147 149 158 167 170 171 181 188 190 207 209 210 213 214 218 219
220 223 224 232 241 243 244 247 248 249 250 253 254 260 263 268
270
 BRRESH=       3
 |Aut(G,M65)|=  8
  3   5  12  14  17  20  29  32  37  40  41  44  49  55  58  64
 73  75  78  79  80  87  89  91  94  96  98  99 101 106 109 116
119 123 124 127 139 152 173 178 194 210 212 213 220 223 227 228
231 234 237 241 243 248 249 254 257 258 259 261 264 265 267 272
273
 BRRESH=       4
 |Aut(G,M65)|=  8
  2   3  12  14  17  20  26  32  37  40  44  46  54  55  58  63
 64  65  67  68  73  80  83  86  94  95  96  98 100 102 108 109
114 119 121 125 127 133 136 150 157 165 169 177 183 210 212 213
215 218 220 221 223 243 248 249 254 258 262 264 265 266 267 269
271
 BRRESH=       5
 |Aut(G,M65)|=  4
  8  10  12  20  27  28  34  40  43  50  53  57  60  61  63  64
 69  73  74  75  85  89  93  96  98 105 106 112 113 115 119 124
133 138 139 142 152 157 159 160 161 169 173 175 177 178 184 189
209 210 213 214 219 220 223 224 242 245 252 255 259 262 266 269
271
 Unitals No        6
 |Aut(G,M65)|=  4
  2   8  10  26  27  28  34  43  46  50  53  54  57  60  61  63
 68  73  76  80  84  94  95  96  98 102 104 109 119 121 127 128
129 132 133 136 149 150 157 158 165 167 169 171 177 183 188 190
209 210 213 214 219 220 223 224 228 231 234 237 258 259 264 265
267
 Unitals No        7
 |Aut(G,M65)|=  8
  5   8  14  17  27  29  37  41  47  49  50  52  58  66  67  72
 73  75  78  79  84  89  93  98  99 101 106 108 110 111 113 124
128 132 135 138 141 179 184 185 190 194 196 198 202 203 210 211
220 221 222 225 226 228 229 232 235 236 237 239 243 245 249 251
252
 Unitals No        8
 |Aut(G,M65)|=  24
 50  65  69  70  71  77  82  83  84  88  92 102 106 109 110 112
115 121 123 124 127 129 133 139 140 142 148 152 155 157 158 163
164 165 170 174 180 181 182 187 191 195 198 200 202 203 209 211
214 222 223 225 228 236 237 239 245 248 250 252 253 265 269 270
272
 Unitals No        9
 |Aut(G,M65)|=  12
 12  14  15  18  21  27  35  37  40  50  52  53  57  60  62  63
 65  66  68  73  74  76  81  89  97  98  99 105 106 107 114 122
131 133 134 143 146 155 156 158 168 172 173 174 179 180 181 185
196 199 200 203 204 207 209 222 227 228 232 235 239 240 245 250
259
 Unitals No       10
 |Aut(G,M65)|=  4
 14  15  17  29  44  46  50  53  55  60  63  64  72  74  80  84
 87  89  90  93  95  96  97 103 108 114 116 120 141 142 144 145
146 148 151 152 154 157 164 171 172 183 184 188 193 194 206 210
211 213 218 220 223 224 228 238 239 242 246 253 259 261 263 265
266
 Unitals No       11
 |Aut(G,M65)|=  4
  2   3  20  32  35  37  49  50  53  58  60  63  65  66  69  72
 76  77  79  88  92  96 102 106 108 115 116 127 129 131 141 146
148 151 153 154 157 160 165 171 173 181 184 186 193 204 206 210
211 212 213 220 223 224 228 232 239 242 253 256 259 260 267 269
273
 Unitals No       12
 |Aut(G,M65)|=  4
 12  14  17  26  41  44  50  53  55  59  60  63  66  69  71  73

```
    76  79  80  81  82  89 101 103 112 114 121 125 129 131 135 148
   151 153 154 156 157 160 161 165 173 181 186 190 194 206 208 210
   213 218 220 222 223 224 237 238 239 246 253 255 259 260 267 269
   273
 Unitals No        13
 |Aut(G,M65)|=  12
    5   8  13  17  20  28  38  46  47  50  51  53  58  59  60  63
   65  66  67  77  86  87  88  92 102 105 106 108 115 125 127 128
  129 131 133 138 148 150 152 159 161 170 172 174 184 187 189 191
  196 203 204 206 213 219 221 222 227 228 229 235 244 245 246 254
  259
G546F, JOWK graph in G. Royle's format, |AutG|=258048
Unitals No         1
 |Aut(G,M65)|=   8
    3   8   9  10  18  26  27  29  38  39  43  46  50  54  57  64
   67  68  70  74  82  84  88  96  98 108 110 112 116 117 119 123
  133 137 141 143 147 151 159 160 165 166 168 172 180 188 189 191
  215 217 219 220 232 237 238 239 243 245 250 254 257 258 259 262
  273
 Unitals No         2
 |Aut(G,M65)|=  16
   14  19  23  25  31  32  33  34  39  43  48  62  68  69  70  74
   75  82  84  88  92  93  97  98 103 107 112 116 117 118 122 123
  142 147 151 153 159 160 165 166 168 169 172 178 180 184 188 189
  197 198 200 201 204 222 225 227 234 237 239 241 243 250 253 255
  257
 Unitals No         3
 |Aut(G,M65)|=   4
   11  17  19  23  31  32  38  39  40  43  47  61  68  69  70  73
   74  82  84  88  89  92  98  99 101 108 112 116 117 118 121 123
  130 145 147 151 159 160 161 180 184 185 188 189 193 194 202 203
  205 218 230 231 232 237 239 243 245 250 252 256 257 261 265 267
  272
 Unitals No         4
 |Aut(G,M65)|=  48
    8  19  23  25  31  32  33  39  40  44  46  54  67  68  73  74
   80  82  83  84  89  96  97 101 102 110 112 116 119 121 123 127
  133 147 148 151 159 160 165 166 168 172 174 180 183 185 189 191
  193 197 198 200 204 220 225 229 230 238 239 241 243 248 252 254
  262
 Unitals No         5
 |Aut(G,M65)|=  24
    1   2   9  10  12  19  23  30  31  32  48  49  50  53  57  58
   67  68  71  74  78  82  84  94  95  96 103 115 116 119 123 126
  129 134 137 139 141 145 146 154 155 157 162 169 170 171 173 180
  189 190 191 192 195 196 199 207 208 209 216 217 219 221 227 255
  262
 Unitals No         6
 |Aut(G,M65)|=  32
    6   7   8   9  12  26  34  37  45  46  47  58  65  67  68  75
   80  81  83  84  91  96  97  99 100 107 112 114 117 125 126 127
  134 135 136 137 140 150 151 152 153 156 170 178 181 189 190 191
  198 199 200 201 204 218 225 227 228 235 240 242 245 253 254 255
  257
 Unitals No         7
 |Aut(G,M65)|=  12
    3   5   6   8  11  32  33  35  42  46  47  64  65  67  75  78
   79  85  87  90  92  94  97 102 104 106 112 118 119 120 123 126
  135 145 147 150 156 160 166 168 170 172 175 179 181 185 187 189
  199 213 215 216 222 223 226 228 234 236 239 241 245 251 252 256
  264
G546G, DSFP graph in G. Royle's format, |AutG|=55296
 Unitals No         1
 |Aut(G,M65)|=  12
    9  11  13  14  21  24  31  32  34  43  44  48  68  69  71  77
   82  87  88  93  98 101 102 110 116 118 120 126 129 133 135 142
  147 150 151 156 161 166 168 173 178 179 185 191 193 195 203 208
  209 217 220 223 228 235 236 239 243 244 249 256 257 258 260 266
  273
  Unitals No         2
 |Aut(G,M65)|=  24
   10  18  33  34  36  42  45  49  66  68  71  75  77  81  82  84
   90  91  99 102 103 105 108 113 118 122 123 125 134 135 137 142
  144 148 163 169 172 174 175 179 182 183 188 190 197 198 199 201
  206 211 216 217 220 222 226 228 235 236 237 241 243 250 251 253
  260
```

G546H,HALL graph in G. Royle's format,|AutG|=921600
 Unitals No     1
 |Aut(G,M65)|=  1200
  9  11  13  14  33  38  43  48  53  56  63  64  66  67  69  73
 81  83  88  93  97 101 108 110 114 120 123 124 135 136 138 142
147 150 151 156 165 166 170 173 177 183 185 191 195 202 203 207
217 218 220 224 226 230 238 239 242 247 253 256 257 258 260 266
273
 Unitals No     2
 |Aut(G,M65)|=  48
  9  11  13  14  21  24  31  32  41  42  44  47  51  54  55  60
 68  69  71  78  85  87  90  94 102 104 106 109 116 118 120 125
129 131 139 144 161 169 172 175 179 186 187 192 193 197 199 206
209 214 216 221 228 233 236 239 243 244 251 256 257 258 260 266
273
 Unitals No     3
 |Aut(G,M65)|=  32
  1   2   4  10  19  22  23  28  33  34  36  42  53  56  63  64
 67  70  71  76  89  91  93  94 105 107 109 110 115 118 119 124
129 130 132 138 153 155 157 158 161 162 164 170 177 178 180 186
195 198 199 204 211 214 215 220 233 235 237 238 249 251 253 254
273
 Unitals No     4
 |Aut(G,M65)|=  100
  1   3   9  11  12  24  34  50  51  52  59  61  67  72  73  78
 79  81  84  85  89  93  98 100 105 108 110 117 120 121 122 123
129 130 131 135 142 145 146 150 156 157 165 171 172 173 176 189
195 196 197 200 204 209 212 216 219 222 229 254 263 264 265 268
269
 Unitals No     5
 |Aut(G,M65)|=  16
  1   6   8  10  15  17  22  23  29  30  46  55  59  60  62  64
 66  68  70  75  77  88 106 113 116 121 125 128 129 132 135 136
140 148 152 154 155 160 164 166 170 172 174 183 202 205 206 207
208 215 216 219 221 223 229 230 235 236 239 241 243 252 255 256
262
 Unitals No     6
 |Aut(G,M65)|=  8
  4   6   7  10  15  29  33  35  38  41  42  49  52  54  61  62
 67  70  77  79  80  85  99 101 103 106 112 113 117 121 125 128
130 132 137 143 144 145 147 151 155 158 164 165 169 170 174 177
180 183 184 192 197 199 205 206 207 218 230 243 249 252 254 255
268
g546i,DEMP graph in G. Royle's format, |AutG|=9216000
 Unitals No     1
 |Aut(G,M65)|=  48
  1  10  11  12  13  23  33  34  37  38  40  51  53  61  62  63
 67  69  77  78  79  87  98 100 105 108 110 113 122 123 124 125
131 136 137 139 144 148 150 154 159 160 162 164 169 172 174 179
184 185 187 192 193 194 197 198 200 215 231 244 246 250 255 256
257
 Unitals No     2
 |Aut(G,M65)|=  16
  3   4   9  15  18  24  27  28  34  37  44  45  51  52  56  59
 67  68  70  74  85  88  91  93  99 100 102 106 113 119 126 128
133 137 141 143 146 153 156 159 163 164 168 171 178 185 188 191
198 200 202 203 210 214 218 220 229 230 234 237 245 249 253 255
273
 Unitals No     3
 |Aut(G,M65)|=  12
  4   6   7   8  15  19  20  22  24  29  33  34  37  43  46  49
 50  53  59  62  73  89 100 102 103 104 111 113 114 115 123 125
135 138 140 143 144 153 161 162 163 171 173 179 180 182 184 189
201 213 218 220 222 224 231 234 236 239 240 245 250 252 254 256
257
 Unitals No     4
 |Aut(G,M65)|=  24
  6   7  10  13  14  18  19  25  26  28  33  34  40  45  48  63
 65  66  72  77  80  83  84  85  87  96 101 104 107 108 110 117
120 123 124 126 130 131 137 138 140 145 148 150 153 155 166 167
170 173 174 191 207 209 212 214 217 219 239 243 244 245 247 256
257
g546j, SEMI2 graph in G. Royle's format, |AutG|=147456
 Unital No=     1
 |Aut(G,M65)|=  64
  7   9  11  12  19  22  29  30  33  34  40  47  51  54  61  62

```
 65  66  72  79  83  86  93  94 103 105 107 108 115 118 125 126
129 130 136 143 148 149 154 160 161 162 168 175 183 185 187 188
193 194 200 207 211 214 221 222 231 233 235 236 247 249 251 252
273
 Unital No=       2
 |Aut(G,M65)|=  64
  3   6  13  14  20  21  26  32  35  38  45  46  52  53  58  64
 65  66  72  79  81  82  88  95 103 105 107 108 115 118 125 126
132 133 138 144 148 149 154 160 161 162 168 175 179 182 189 190
193 194 200 207 211 214 221 222 228 229 234 240 241 242 248 255
273
 Unital No=       3
 |Aut(G,M65)|=  64
  3   6  13  14  19  22  29  30  36  37  42  48  51  54  61  62
 67  70  77  78  87  89  91  92 103 105 107 108 113 114 120 127
135 137 139 140 151 153 155 156 161 162 168 175 177 178 184 191
193 194 200 207 211 214 221 222 225 226 232 239 247 249 251 252
273
Unital No=       4
 |Aut(G,M65)|=  192
  2   7  11  15  17  24  25  28  36  37  38  46  52  53  54  62
 65  72  73  76  84  85  86  94  98 103 107 111 114 119 123 127
132 133 134 142 145 152 153 156 161 168 169 172 180 181 182 190
194 199 203 207 209 216 217 220 226 231 235 239 243 250 253 256
273
 Unital No=       5
 |Aut(G,M65)|=  48
  2   7  11  15  17  28  29  32  33  44  45  48  56  57  59  63
 65  76  77  80  84  85  86  94 107 109 111 112 114 115 119 122
131 136 137 138 146 147 151 154 161 168 169 172 184 185 187 191
193 194 199 204 210 211 215 218 232 233 235 239 243 250 253 256
273
 Unital No=       6
 |Aut(G,M65)|=  192
  2   7  11  15  17  24  25  28  33  40  41  44  49  56  57  60
 67  74  77  80  84  85  86  94  98 103 107 111 115 122 125 128
131 138 141 144 147 154 157 160 161 168 169 172 178 183 187 191
193 200 201 204 210 215 219 223 226 231 235 239 243 250 253 256
273
 Unital No=       7
 |Aut(G,M65)|=  48
  3   4   6  10  20  21  22  30  33  44  45  48  53  56  57  62
 65  68  70  76  82  83  87  90  97 104 105 108 113 114 119 124
131 132 134 138 147 152 153 154 161 168 169 172 184 185 187 191
194 197 199 206 210 211 215 218 228 229 230 238 242 245 247 254
273
 Unital No=       8
 |Aut(G,M65)|=  192
  4   5   6  14  17  24  25  28  35  42  45  48  49  56  57  60
 65  72  73  76  83  90  93  96  98 103 107 111 114 119 123 127
131 138 141 144 147 154 157 160 161 168 169 172 177 184 185 188
194 199 203 207 211 218 221 224 226 231 235 239 242 247 251 255
273
 Unital No=       9
 |Aut(G,M65)|=  48
  5  13  14  16  17  28  29  32  36  37  38  46  52  54  59  63
 69  77  78  80  91  93  95  96 107 109 111 112 115 120 121 122
133 136 137 142 145 146 151 156 161 168 169 172 180 181 182 190
193 196 198 204 209 216 217 220 232 233 235 239 244 246 251 255
273
 Unital No=      10
 |Aut(G,M65)|=  64
  2   4   8  16  18  20  24  32  35  43  44  46  49  53  58  63
 65  69  74  79  83  91  92  94  97 101 106 111 113 117 122 127
129 133 138 143 147 155 156 158 162 164 168 176 179 187 188 190
195 203 204 206 210 212 216 224 230 231 233 237 242 244 248 256
273
 Unital No=      11
 |Aut(G,M65)|=  16
  2   3  11  16  19  27  28  30  36  40  42  47  49  53  58  63
 65  69  70  73  81  85  92  94  98 103 109 112 113 117 122 127
131 138 139 143 146 147 155 160 167 172 173 174 183 188 189 190
194 199 205 208 209 210 213 224 231 234 237 239 243 251 252 254
273
 Unital No=      12
 |Aut(G,M65)|=  64
```

```
  1   5   6   9  20  24  28  30  39  42  45  47  52  56  60  62
 66  67  75  80  87  90  93  95 103 106 109 111 114 115 123 128
132 136 140 142 146 147 155 160 167 170 173 175 180 184 188 190
194 195 203 208 210 211 219 224 231 234 237 239 244 248 252 254
273
 Unital No=         13
|Aut(G,M65)|=   64
  2   8  11  12  18  24  27  28  37  38  42  45  51  52  62  64
 67  68  78  80  82  88  91  92  98 104 107 108 115 116 126 128
133 134 138 141 149 150 154 157 162 168 171 172 179 180 190 192
193 199 201 207 213 214 218 221 229 230 234 237 243 244 254 256
273
 Unital No=         14
|Aut(G,M65)|=   64
  2   4  11  14  19  24  28  32  37  38  39  47  50  52  59  62
 69  70  71  79  85  86  87  95  99 104 108 112 117 118 119 127
129 137 138 141 146 148 155 158 163 168 172 176 178 180 187 190
197 198 199 207 211 216 220 224 226 228 235 238 243 248 252 256
273
 Unital No=         15
|Aut(G,M65)|=   64
  2   8  11  12  18  24  27  28  34  40  43  44  51  52  62  64
 67  68  78  80  85  86  90  93  98 104 107 108 113 119 121 127
133 134 138 141 147 148 158 160 163 164 174 176 181 182 186 189
197 198 202 205 211 212 222 224 229 230 234 237 242 248 251 252
273
 Unital No=         16
|Aut(G,M65)|=   64
  6   8  10  11  22  24  26  27  38  40  42  43  50  53  60  61
 65  71  78  80  81  87  94  96  97 103 110 112 114 117 124 125
130 133 140 141 146 149 156 157 161 167 174 176 182 184 186 187
193 199 206 208 210 213 220 221 230 232 234 235 243 244 249 255
273
 Unital No=         17
|Aut(G,M65)|=   64
  5   7  12  16  19  22  24  31  35  38  40  47  49  50  61  62
 69  71  76  80  81  82  93  94  97  98 109 110 115 118 120 127
133 135 140 144 147 150 152 159 165 167 172 176 180 185 186 187
197 199 204 208 209 210 221 222 227 230 232 239 241 242 253 254
273
 Unital No=         18
|Aut(G,M65)|=   32
  1   7  14  16  19  20  25  31  35  36  41  47  50  53  60  61
 70  72  74  75  82  85  92  93  99 100 105 111 118 120 122 123
131 132 137 143 146 149 156 157 166 168 170 171 182 184 186 187
195 196 201 207 210 213 220 221 226 229 236 237 246 248 250 251
273
 Unital No=         19
|Aut(G,M65)|=   32
  5   7  12  16  19  22  24  31  36  41  42  43  52  57  58  59
 67  70  72  79  85  87  92  96 100 105 106 107 116 121 122 123
129 130 141 142 147 150 152 159 165 167 172 176 180 185 186 187
195 198 200 207 213 215 220 224 229 231 236 240 243 246 248 255
273
 Unital No=         20
|Aut(G,M65)|=   64
  1   4   5   8  22  25  28  30  35  39  43  45  50  58  63  64
 67  71  75  77  82  90  95  96  99 103 107 109 118 121 124 126
131 135 139 141 146 154 159 160 166 169 172 174 179 183 187 189
198 201 204 206 210 218 223 224 226 234 239 240 246 249 252 254
273
 Unital No=         21
|Aut(G,M65)|=   64
  1   5  12  14  17  21  28  30  34  39  45  48  49  53  60  62
 66  71  77  80  83  90  91  95 100 102 104 105 113 117 124 126
130 135 141 144 148 150 152 153 161 165 172 174 180 182 184 185
196 198 200 201 210 215 221 224 226 231 237 240 244 246 248 249
273
G546K, SEMI4 graph in G. Royle's format, group order=884736
 Unital No=          1
|Aut(G,M65)|=    8
  1   5   8   9  10  17  21  24  25  26  39  44  46  47  48  49
 53  56  57  58  66  67  68  70  77  87  92  94  95  96 103 108
110 111 112 114 115 116 118 125 139 145 149 152 153 154 171 183
188 190 191 192 203 219 226 227 228 230 237 242 243 244 246 253
257
```

```
 Unital No=        2
|Aut(G,M65)|=   8
  3  11  12  14  18  21  24  26  33  37  42  47  50  54  56  61
 65  69  74  79  82  85  88  90  97  99 110 111 118 123 124 125
133 135 137 138 146 150 152 157 163 171 172 174 177 179 190 191
198 203 204 205 212 219 220 224 226 227 232 238 241 246 253 255
273
 Unital No=        3
|Aut(G,M65)|= 128
  7   9  11  12  20  21  26  32  36  37  42  48  51  54  61  62
 68  69  74  80  84  85  90  96  99 102 109 110 119 121 123 124
131 134 141 142 148 149 154 160 163 166 173 174 183 185 187 188
195 198 205 206 215 217 219 220 225 226 232 239 247 249 251 252
273
 Unital No=        4
|Aut(G,M65)|=  64
  3   6  13  14  17  18  24  31  33  34  40  47  51  54  61  62
 65  66  72  79  81  82  88  95  99 102 109 110 115 118 125 126
131 134 141 142 148 149 154 160 167 169 171 172 183 185 187 188
199 201 203 204 215 217 219 220 225 226 232 239 247 249 251 252
273
 Unital No=        5
|Aut(G,M65)|=  48
  6   7   8  16  24  27  30  32  34  36  40  48  57  59  61  62
 66  68  73  77  83  86  87  92  97 102 103 106 115 118 119 124
137 139 141 142 147 155 156 158 166 167 169 173 179 181 188 191
194 195 196 204 210 212 217 221 226 228 232 240 248 251 254 256
273
 Unital No=        6
|Aut(G,M65)|= 192
  1   2   4  10  21  24  31  32  41  43  45  46  49  50  52  58
 67  70  71  76  85  88  95  96  99 102 103 108 117 120 127 128
129 130 132 138 149 152 159 160 161 162 164 170 177 178 180 186
195 198 199 204 211 214 215 220 227 230 231 236 245 248 255 256
273
 Unital No=        7
|Aut(G,M65)|=   8
  2   3  13  15  20  23  25  32  36  38  40  41  49  50  61  62
 69  70  72  76  82  85  92  93  97 106 107 110 113 119 126 128
129 130 141 142 149 151 156 160 165 167 172 176 178 180 185 189
193 198 200 206 212 215 217 224 228 230 232 233 245 246 248 252
273
 Unital No=        8
|Aut(G,M65)|=  12
  2   6  14  15  17  27  28  31  37  40  47  48  57  59  61  62
 65  69  70  73  89  90  92  96  97 101 102 105 118 120 122 123
134 136 138 139 147 149 154 158 169 170 172 176 180 184 188 190
193 199 206 208 213 216 223 224 225 235 236 239 263 266 269 271
273
 Unital No=        9
|Aut(G,M65)|= 128
  7  10  13  15  23  26  29  31  36  40  44  46  55  58  61  63
 66  67  75  80  81  85  86  89  97 101 102 105 113 117 118 121
132 136 140 142 148 152 156 158 167 170 173 175 183 186 189 191
193 197 198 201 209 213 214 217 228 232 236 238 244 248 252 254
273
 Unital No=       10
|Aut(G,M65)|= 128
  4   6   8   9  18  23  29  32  34  39  45  48  50  55  61  64
 67  74  75  79  83  90  91  95  97 101 108 110 114 119 125 128
131 138 139 143 148 150 152 153 164 166 168 169 180 182 184 185
196 198 200 201 211 218 219 223 226 231 237 240 243 250 251 255
273
 Unital No=       11
|Aut(G,M65)|=  12
  1   4   7   8  12  32  34  37  42  43  46  50  53  58  59  62
 67  70  73  77  79  83  86  89  93  95  97 100 103 104 108 128
129 132 135 136 140 147 150 153 157 159 176 192 195 198 201 205
207 210 213 218 219 222 225 228 231 232 236 242 245 250 251 254
257
 Unital No=       12
|Aut(G,M65)|=   8
 18  21  28  29  37  40  47  48  52  56  60  62  65  67  72  77
 84  86  93  96  97 107 108 111 113 114 116 122 129 135 142 144
147 149 154 158 166 168 170 171 178 182 190 191 195 198 199 204
212 213 215 219 226 227 235 240 247 250 253 255 257 263 270 272
```

273
G546L, LMRH graph in G. Royle's format, |AutG|=258048
 Unitals No     1
 |Aut(G,M65)|=  32
  3  5  6 16 17 18 25 27 33 34 41 43 49 50 57 59
 67 69 70 80 81 82 89 91 99 101 102 112 115 117 118 128
132 138 141 142 151 152 156 159 167 168 172 175 183 184 188 191
199 200 204 207 209 210 217 219 231 232 236 239 243 245 246 256
273
 Unitals No     2
 |Aut(G,M65)|=  16
  3  6 10 11 13 19 22 26 27 29 34 37 41 46 47 51
 54 58 59 61 65 82 85 89 94 95 98 101 105 110 111 113
132 135 136 140 144 147 150 154 155 157 164 167 168 172 176 178
181 185 190 191 196 199 200 204 208 212 215 216 220 224 225 241
257
 G546M, BBH2 graph in G. Royle's format, group order=3840
 Unital No     1
 |Aut(G,M65)|=  16
  2  4 15 17 21 23 24 25 30 32 33 37 38 39 42 43
 48 51 52 54 58 62 63 64 65 73 76 83 84 87 95 100
102 106 116 120 121 131 134 135 147 148 158 161 165 169 177 183
189 199 204 207 215 216 217 227 231 232 244 246 252 261 263 269
272
 Unital No     2
 |Aut(G,M65)|=  16
  4  6 13 14 17 18 20 23 25 26 30 32 35 36 41 42
 54 58 59 63 67 70 75 76 83 99 102 111 112 114 117 123
127 131 133 134 135 146 147 158 160 162 165 167 173 195 197 200
205 213 214 217 224 226 231 232 236 246 252 254 255 270 271 272
273
 Unital No     3
 |Aut(G,M65)|=  16
  2  4 15 17 18 21 23 25 31 32 34 36 44 54 56 63
 65 73 76 83 85 86 87 97 107 112 116 120 121 130 131 134
147 158 160 161 163 173 180 184 189 194 195 198 200 201 206 208
215 216 217 228 237 239 241 246 247 249 252 253 254 261 263 269
272
 Unital No     4
 |Aut(G,M65)|=  32
  1  2  3  6  7  8 10 13 14 15 28 29 30 34 44 58
 65 66 67 69 70 74 75 77 78 80 82 83 90 91 98 104
116 117 119 120 123 125 136 156 163 164 170 182 188 191 193 196
209 212 213 216 218 224 231 232 249 258 261 263 265 266 268 271
273
 Unital No     5
 |Aut(G,M65)|=  32
  1 14 19 27 28 29 30 33 34 36 40 42 43 44 46 58
 66 75 82 83 84 85 90 91 97 98 102 104 105 108 109 110
123 125 130 135 136 137 144 148 149 150 156 160 163 164 165 170
173 182 185 186 188 191 193 195 196 198 218 224 229 231 232 235
249
 Unital No     6
 |Aut(G,M65)|=  4
  2  4  8 10 33 37 38 43 52 57 62 63 71 72 75 76
 84 87 91 92 98 100 103 104 113 114 116 120 122 124 127 128
134 136 138 140 148 152 154 156 168 169 172 174 176 177 181 191
193 196 203 204 209 214 217 222 242 246 252 256 258 259 263 270
272
 Unital No     7
 |Aut(G,M65)|=  4
  2  4  8 10 36 42 45 47 49 50 56 60 71 72 75 76
 82 89 93 101 102 109 110 113 114 116 120 122 124 127 128 132
135 137 142 145 151 159 160 161 164 165 173 175 182 184 186 188
194 195 201 202 209 214 217 222 245 248 250 251 258 259 263 270
272
 Unital No     8
 |Aut(G,M65)|=  8
  2  5  9 14 17 23 25 32 36 37 38 42 54 59 62 63
 69 75 77 80 83 98 99 102 107 116 120 123 126 131 134 135
139 147 156 158 160 162 165 167 173 177 180 190 191 193 195 197
208 210 211 221 224 225 226 236 238 241 246 252 255 257 262 263
265
 Unital No     9
 |Aut(G,M65)|=  20
 51 52 55 56 59 60 61 63 69 75 77 80 82 85 88 95

```
 97  99 100 103 106 107 109 110 129 130 132 134 137 138 139 141
161 162 164 167 169 175 176 180 182 184 186 187 188 190 192 194
197 198 199 201 203 204 208 225 226 236 238 241 246 252 255 264
265
 Unital No      10
|Aut(G,M65)|=  10
  2   4   6  11  12  13  16  19  21  29  35  43  45  60  61  62
 66  67  70  71  72  77  78  87  94 101 106 108 115 116 118 129
134 144 153 154 156 162 168 173 175 181 182 183 184 188 190 194
198 200 221 222 223 228 233 235 246 249 253 263 264 265 269 270
273
 Unital No      11
|Aut(G,M65)|=   8
  3   9  16  17  18  22  23  24  25  32  39  44  46  51  59  63
 65  68  70  83  85  91  96 100 101 102 104 106 107 111 114 115
118 119 120 124 127 130 131 138 150 156 158 162 166 171 177 179
185 195 202 205 216 221 222 225 235 237 244 252 254 260 265 269
271
 Unital No      12
|Aut(G,M65)|=   8
  5   7   8  11  13  15  16  17  22  23  24  25  26  32  36  37
 39  51  56  63  68  74  79  83  85  95  96 100 101 104 107 109
111 112 116 124 126 131 138 144 148 156 158 162 163 166 179 185
190 195 196 207 214 217 222 234 235 240 252 254 255 259 263 267
273
 Unital No      13
|Aut(G,M65)|=  16
  4   5   7   9  10  15  18  19  23  24  30  35  43  45  46  56
 59  71  80  84  89  91  93  98 102 103 104 105 106 111 112 118
120 131 132 135 137 141 142 149 160 162 165 171 176 177 180 183
186 193 201 202 203 209 211 226 233 237 238 244 255 260 263 268
270
 Unital No      14
|Aut(G,M65)|=  16
  4   7   9  10  19  21  23  24  29  33  35  45  48  58  60  61
 62  82  84  87  93  98  99 103 104 106 110 111 112 131 132 133
134 139 140 141 142 146 151 156 157 162 164 165 166 179 182 183
186 193 198 203 206 229 233 236 238 241 245 246 254 260 263 268
270
 Unital No      15
|Aut(G,M65)|=  16
  4   7   9  10  18  23  24  27  35  38  43  48  56  58  59  61
 62  81  84  91  99 100 104 105 133 134 135 137 139 146 149 156
157 160 165 167 171 177 183 185 198 199 201 203 210 212 214 216
226 229 231 238 241 244 245 254 255 258 260 262 263 268 269 270
273
 Unital No      16
|Aut(G,M65)|=  16
  1   6  11  13  18  23  24  26  35  36  43  58  62  72  74  76
 78  84  91  95 104 105 109 113 114 115 116 117 119 123 128 133
139 146 156 163 165 171 177 183 190 201 203 207 210 212 213 214
215 216 221 222 226 238 240 241 254 257 258 261 262 267 269 271
273
 Unital No      17
|Aut(G,M65)|=  16
  1   6  22  23  24  28  35  39  48  54  56  59  61  72  74  84
 90  96  99 101 104 113 115 116 117 123 128 129 134 135 137 149
152 157 160 161 165 170 183 188 189 195 198 203 210 212 214 216
221 222 229 235 238 244 245 255 256 257 258 261 262 267 269 271
273
 Unital No      18
|Aut(G,M65)|=   8
  4  10  11  13  18  22  28  30  32  44  47  60  61  71  73  77
 80  86  90  95  96  98  99 101 108 110 112 114 115 116 119 129
143 146 151 153 156 157 159 163 164 166 169 177 179 180 182 205
208 209 211 214 216 218 220 222 224 234 239 248 256 257 260 266
273
 Unital No      19
|Aut(G,M65)|=   8
  6   8  26  28  29  30  32  33  35  38  39  42  48  49  54  59
 60  69  74  86  87  91  96 107 111 121 122 123 124 125 126 129
130 133 135 136 138 139 140 147 148 157 160 163 164 170 176 177
179 189 191 199 203 217 223 228 233 243 245 248 255 259 262 268
271
 Unital No      20
|Aut(G,M65)|=   8
```

```
  1  11  13  16  21  22  24  25  28  29  34  35  39  40  42  62
 63  72  75  76  78  81  99 102 104 107 110 118 120 122 123 132
136 153 154 163 169 171 173 176 177 179 180 182 186 188 189 190
191 198 200 201 203 206 209 211 218 222 228 230 235 238 240 251
254
```
 Unital No     21
 |Aut(G,M65)|=  8
```
  7   9  21  22  28  29  31  32  39  42  44  52  53  54  55  56
 62  73  77  82  84  86  87  89  90  91  95  96  99 106 110 113
117 131 133 136 142 143 144 145 148 151 153 155 156 167 177 180
186 190 191 198 205 206 214 216 225 230 235 245 250 252 253 254
255
```
 Unital No     22
 |Aut(G,M65)|=  8
```
  1   9  10  11  14  15  17  21  22  28  29  35  38  51  53  62
 63  68  70  72  73  76  80 100 104 113 116 118 119 123 125 131
132 136 139 145 146 153 154 177 180 183 184 185 186 190 191 199
203 211 212 213 214 220 222 231 238 249 250 251 254 266 267 268
269
```
 Unital No     23
 |Aut(G,M65)|=  8
```
  1   9  10  11  15  16  19  32  34  35  47  55  56  58  61  62
 63  70  72  73  75  76  80  81 103 104 107 113 116 118 119 122
123 132 136 138 140 142 144 147 148 153 154 155 159 165 184 200
203 208 211 212 213 214 218 222 228 232 238 241 251 252 254 255
256
```
 Unital No     24
 |Aut(G,M65)|=  4
```
  6  10  12  15  17  18  20  23  25  26  30  32  34  40  47  48
 56  61  62  63  66  70  74  80  82  83  87  90  96 101 102 103
107 116 119 124 128 132 136 140 144 147 148 153 154 194 200 201
208 212 213 221 224 228 232 236 240 251 254 255 256 266 267 268
269
```
 Unital No     25
 |Aut(G,M65)|=  4
```
  6  10  11  12  17  21  22  23  25  28  29  32  33  34  44  46
 51  54  56  63  66  70  74  76  83  97 106 107 108 116 118 124
128 132 139 143 144 146 148 151 154 163 169 171 176 195 200 205
207 211 212 221 224 225 226 228 229 245 249 251 255 258 259 260
261
```
 Unital No     26
 |Aut(G,M65)|=  4
```
  4   7   9  10  17  19  23  24  25  27  31  32  33  34  35  38
 39  42  45  48  49  51  52  56  58  59  62  64  71  76  78  80
 83 121 126 127 128 145 147 151 152 154 155 157 158 161 163 164
166 169 170 171 176 196 197 201 202 204 205 207 208 218 220 222
224
```
**Appendix 2: Unitals of projective planes of order 16, found with author's algorithm on Moorhouse representation of the planes [9]**
 pp16htmmoor\SEMI4.DAT
 Unital No  1
 |Aut(G,M65)|=  8
```
  1  23  27  29  32  34  38  41  45  47  48  51  56  59  61  69
 70  73  74  77  79  83  84  91  93  95 102 105 117 118 122 124
127 131 132 136 139 141 146 148 150 151 162 166 167 170 174 185
190 191 193 199 208 217 218 233 237 243 246 249 250 251 258 260
271
```
 Unital No  2
 |Aut(G,M65)|=  128
```
 17  22  26  30  31  39  42  44  46  48  54  62  63  66  70  73
 82  85  87  89  90  92  96 104 108 112 114 129 134 137 138 139
140 145 151 157 161 164 180 183 185 186 187 191 193 197 199 204
207 208 209 213 214 218 223 232 237 241 247 248 252 258 259 267
273
```
 Unital No  3
 |Aut(G,M65)|=  8
```
 19  21  28  31  32  33  34  38  39  44  49  54  58  62  68  69
 70  75  76  77  79  87  92  96 103 104 110 113 116 119 124 131
132 138 142 145 147 149 155 157 158 168 178 180 182 190 200 207
211 220 222 226 228 235 236 237 241 251 252 258 259 261 271 272
273
```
 Unital No  4
 |Aut(G,M65)|=  192
```
 17  19  21  27  32  40  42  44  46  49  59  61  69  77  78  80
 86  90  95 102 106 109 110 111 116 118 122 130 132 137 142 145
164 165 171 175 176 178 179 183 188 194 197 200 204 213 222 228
```

```
229 231 237 238 241 249 253 254 259 261 262 264 266 267 268 269
273
 Unital No  5
 |Aut(G,M65)|=  48
 19 37 40 42 44 47 49 51 61 62 69 76 78 80 87 94
101 102 106 109 111 113 116 118 122 128 130 132 133 141 144 145
147 151 161 164 166 175 176 178 183 187 188 194 210 224 228 229
232 237 238 241 243 248 253 257 259 261 262 264 266 267 268 269
273
 Unital No  6
 |Aut(G,M65)|=  128
 20 25 28 29 30 35 37 44 46 47 50 51 52 54 66 68
 73 74 75 79 94 96 97 100 104 107 109 111 120 123 124 125
132 137 139 152 157 167 168 169 171 172 173 180 187 188 194 196
205 206 207 209 211 218 222 226 228 230 231 236 240 246 270 271
273
 Unital No  7
 |Aut(G,M65)|=  8
 17 19 21 24 28 34 35 39 47 50 52 53 67 70 71 73
 75 99 102 107 119 120 125 126 128 131 137 139 143 144 148 154
157 162 163 165 169 172 174 175 181 182 187 188 192 194 195 196
207 208 218 224 226 231 232 236 237 239 249 252 253 260 265 268
273
 Unital No  8
 |Aut(G,M65)|=  64
 17 20 24 25 28 34 35 37 39 49 53 56 63 67 68 73
 74 80 95 102 103 107 120 121 123 132 137 139 143 144 152 154
157 163 164 165 166 169 170 172 173 180 188 192 196 198 199 206
207 209 211 213 216 218 221 222 225 226 228 230 231 237 252 253
273
 Unital No  9
 |Aut(G,M65)|=  12
  1 26 28 31 32 34 36 39 40 45 49 50 52 60 62 63
 64 66 68 72 82 93 96 102 103 107 110 116 124 131 132 134
136 137 140 143 144 147 151 168 169 170 182 189 194 196 197 198
199 208 214 218 219 229 230 232 235 238 241 251 252 256 257 258
270
 Unital No 10
 |Aut(G,M65)|=  128
 18 24 25 30 32 33 35 36 40 45 51 53 60 61 64 65
 75 84 93 95 97 99 108 116 122 124 130 133 146 151 156 157
159 163 168 169 170 174 183 185 188 192 194 195 197 198 203 206
210 226 230 233 234 235 237 244 251 254 255 256 261 262 268 270
273
 Unital No 11
 |Aut(G,M65)|=  12
  1  2  7 10 21 24 25 31 35 38 39 44 50 52 58 61
 65 69 71 73 75 76 79 80 89 93 95 97 99 104 105 108
109 119 122 130 133 142 152 154 166 172 178 179 183 189 211 221
222 223 227 233 235 237 239 241 243 253 257 259 260 266 268 270
273
 Unital No 12
 |Aut(G,M65)|=  8
 10 12 14 16 19 25 26 28 31 47 48 49 50 62 66 67
 68 72 81 82 88 91 92 93 103 104 107 112 129 135 136 137
139 141 143 144 152 165 166 168 172 174 177 179 194 195 199 201
211 214 219 222 224 235 237 238 243 250 252 254 255 261 265 268
273
pp16htmmoor\demp.DAT
 Unital No  1
 |Aut(G,M65)|=  48
  1 20 31 32 33 35 37 40 42 47 61 63 64 68 72 79
 82 85 88 98 99 103 104 107 111 116 118 121 126 127 133 134
136 142 147 148 150 152 153 156 160 166 172 178 191 192 193 202
204 210 215 216 220 232 234 239 240 253 255 257 258 264 265 272
273
 Unital No  2
 |Aut(G,M65)|=  16
  1  2  9 10 15 37 38 43 44 49 51 55 56 61 64 65
 67 74 81 83 84 87 88 93 96 101 102 109 112 113 118 119
121 122 130 132 136 141 143 144 147 155 161 164 171 176 178 179
201 203 210 212 214 219 220 225 227 230 234 241 245 247 250 261
271
 Unital No  3
 |Aut(G,M65)|=  24
  1 25 27 30 33 34 35 36 39 43 50 55 58 61 66 76
```

```
 77  81  82  86  87  88  92  93  94  98 104 106 110 112 117 122
124 128 135 136 138 140 146 149 152 155 160 167 169 171 180 184
191 196 198 199 210 218 229 239 242 247 248 249 252 257 263 268
273
 Unital No  4
 |Aut(G,M65)|=  12
  1  21  29  30  31  44  50  63  67  72  76  78  79  86  90  98
 99 100 103 109 110 111 113 117 118 121 124 128 131 134 137 145
146 154 161 163 164 172 173 174 179 184 190 192 194 195 201 203
207 210 211 214 220 224 227 237 242 245 246 255 256 257 266 267
272
pp16htmmoor\dsfp.DAT
 Unital No  1
 |Aut(G,M65)|=  24
  5   6   7   8  22  25  34  37  38  41  42  45  60  61  62  63
 84  94  98 103 106 107 117 118 122 123 128 135 141 144 151 158
163 164 168 171 172 175 178 181 183 185 190 198 202 217 218 219
222 227 228 231 234 241 242 248 251 252 254 255 265 266 268 271
273
 Unital No  2
 |Aut(G,M65)|=  12
 21  22  30  64  67  92  94  98 100 101 103 104 107 109 116 117
128 131 134 135 139 140 144 150 151 158 161 164 166 168 171 172
174 176 178 180 183 185 187 188 191 192 207 211 214 217 218 223
224 226 230 233 236 237 239 241 246 251 254 255 257 258 262 264
271
pp16htmmoor\hall.DAT
 Unital No  1
 |Aut(G,M65)|=  1200.00000000000
  2   8  24  38  67  82  85  87  90  96  97  98  99 105 111 113
117 118 122 125 126 128 130 137 141 144 146 148 157 159 162 163
169 170 181 183 201 205 214 215 216 217 219 221 226 229 231 233
238 240 242 243 246 247 250 251 252 253 255 258 259 263 264 265
270
 Unital No  2
 |Aut(G,M65)|=  100
 12  13  16  21  22  27  30  32  33  34  36  37  39  46  47  48
 55  56  57  60  61  63  68  71  74  87  98 112 122 124 130 132
140 145 146 150 160 166 169 174 181 207 210 211 214 215 216 217
227 229 231 232 237 238 242 243 246 247 250 251 253 255 258 263
270
 Unital No  3
 |Aut(G,M65)|=  48
  9  13  19  56  72  86  88  89  91  92  95  97  99 101 108 110
114 116 118 123 124 125 129 135 136 137 141 144 146 149 150 151
157 161 164 168 184 187 188 189 190 192 198 202 203 206 208 209
210 211 215 219 222 228 237 243 247 249 258 259 261 263 265 268
271
 Unital No  4
 |Aut(G,M65)|=  16
  4   5   8  17  19  20  37  38  39  45  46  51  54  56  60  61
 64  68  74  77  83  86  89  95  96 100 103 104 107 109 110 121
123 128 136 138 148 149 164 175 178 183 186 189 192 194 197 199
202 204 207 213 215 219 222 230 235 239 240 246 247 256 264 270
273
 Unital No  5
 |Aut(G,M65)|=  8
  1   3  12  24  28  30  31  34  39  41  50  51  59  64  65  66
 69  74  75  78  83  86  89  95  96  99 100 103 106 108 117 118
121 123 134 135 136 138 139 140 145 148 157 164 175 178 179 183
189 192 193 197 199 202 207 208 213 222 235 238 246 247 255 270
273
 Unital No  6
 |Aut(G,M65)|=  32
  3   5   6   7  18  26  46  47  50  52  54  57  61  75  78  79
 83  90  91  98 100 101 104 105 108 111 120 129 136 139 140 148
153 154 161 166 168 170 174 177 178 180 185 186 187 188 193 201
217 222 227 228 233 234 239 240 244 246 254 257 262 263 270 271
273
pp16htmmoor\john.DAT
 Unital No  1
 |Aut(G,M65)|=  48
  1   2   3   4   5   6   7   8   9  10  11  12  61  65  67  71
 76  77  78  82  85  90  98 105 109 114 115 120 123 128 136 139
140 142 148 150 154 156 157 166 172 175 179 185 189 190 200 201
202 204 211 215 219 220 225 228 233 238 242 245 253 266 271 272
```

273
 Unital No  2
 |Aut(G,M65)|=  32
  1   2   3   4   5   6   7   8  61  66  70  71  75  77  90  92
 93  96  98 104 105 107 109 121 123 132 134 136 140 141 142 145
148 154 165 169 171 172 176 181 186 189 200 202 208 210 218 221
222 225 229 230 233 239 242 244 253 266 267 268 269 270 271 272
273
 Unital No  3
 |Aut(G,M65)|=  32
 38  39  40  45  49  53  54  55  61  63  68  70  71  73  80  86
 87  89  92  93 102 104 105 111 113 119 123 124 127 132 134 136
138 144 146 155 163 168 170 171 176 182 183 193 196 203 207 208
210 212 213 214 218 222 236 237 252 254 255 256 257 259 260 261
262
 Unital No  4
 |Aut(G,M65)|=  16
  9  10  11  12  38  39  40  45  61  62  65  66  67  70  75  78
 82  85  92  96 107 110 114 115 118 121 128 129 134 139 141 145
149 156 157 162 165 169 179 181 185 186 190 192 204 211 218 219
220 221 228 229 230 238 239 241 244 254 255 256 257 259 260 261
262
 Unital No  5
 |Aut(G,M65)|=  4
  3  25  28  57  63  66  68  70  73  74  78  87  89  91  92  98
101 104 105 113 115 116 117 121 126 129 134 137 141 142 150 153
154 160 165 169 170 173 174 178 184 187 194 203 204 208 210 212
215 217 218 220 223 225 226 228 232 234 235 239 244 246 247 250
268
 Unital No  6
 |Aut(G,M65)|=  4
  3  23  29  60  66  67  71  72  73  78  84  88  89  91  93 101
103 105 107 114 115 118 119 121 129 136 137 138 141 145 146 150
160 162 167 169 170 173 179 183 187 188 194 195 196 197 203 209
210 215 219 221 222 223 225 232 235 236 237 240 244 247 249 251
270
 Unital No  7
 |Aut(G,M65)|=  4
  1  26  28  59  68  70  71  72  75  77  78  89  90  91  93  96
101 105 108 116 117 118 119 122 124 125 127 129 136 137 141 147
150 152 157 159 160 162 169 170 173 179 180 181 187 189 191 194
203 207 210 214 215 218 219 221 222 226 234 239 240 244 247 248
270
 Unital No  8
 |Aut(G,M65)|=  8
 23  24  25  26  37  46  47  48  61  66  67  71  78  84  86  88
 90  95  98 101 105 114 118 121 122 129 130 138 143 145 146 147
148 150 166 167 173 175 183 185 188 189 194 203 204 210 212 215
219 222 223 235 236 237 244 246 247 254 255 256 257 266 271 272
273
 Unital No  9
 |Aut(G,M65)|=  16
 14  23  24  25  26  34  35  36  38  39  40  45  50  51  52  56
 61  77  78  79  89  90  92  95  97  98 102 106 109 122 123 124
128 131 135 136 137 146 149 154 159 161 162 163 164 171 173 174
175 178 186 188 191 195 202 203 205 206 211 212 214 243 244 248
251
 Unital No 10
 |Aut(G,M65)|=  24
 14  16  18  21  22  23  24  25  26  34  35  36  61  74  80  81
 85  90  92  93  95  99 100 103 109 113 114 116 118 123 126 127
128 132 137 138 140 148 149 150 152 157 163 167 169 196 203 204
206 207 211 215 217 223 225 229 231 233 234 238 240 267 268 269
270
 Unital No 11
 |Aut(G,M65)|=  16
  9  10  11  12  38  39  40  45  49  53  54  55  61  64  65  66
 70  73  74  80  81  85  90  92  99 103 109 113 114 116 123 124
128 129 148 150 151 161 170 172 174 177 180 181 182 185 187 188
210 211 212 215 224 232 233 238 242 246 248 249 252 254 255 256
257
 Unital No 12
 |Aut(G,M65)|=  16
  1   2   3   4   9  10  11  12  61  65  68  70  78  84  86  88
 90  95  98 101 105 114 118 121 125 129 131 137 142 144 148 149
157 163 166 167 177 180 186 187 192 193 194 210 212 214 216 217

218 219 227 228 229 239 243 250 252 266 267 268 269 270 271 272 273
 Unital No 13
|Aut(G,M65)|=  8
  1  2  3  4 23 24 25 26 37 46 47 48 49 53 54 55
 61 63 68 69 84 85 87 93 94 97 105 106 109 111 116 121
122 128 130 134 139 145 148 149 152 159 160 161 173 177 178 180
182 187 194 197 202 204 207 208 219 221 225 230 237 238 245 251 252
 Unital No 14
|Aut(G,M65)|=  8
  5  6  7  8 57 58 59 60 61 72 75 77 82 85 103 105
107 108 109 110 117 118 122 126 127 135 136 138 154 160 161 163
164 167 171 172 173 174 175 185 187 193 199 212 213 217 218 220
223 227 230 232 237 238 243 244 252 258 263 264 265 267 268 269 270
 Unital No 15
|Aut(G,M65)|=  8
 16 18 21 22 23 24 25 26 27 28 29 30 38 39 40 45
 61 63 82 83 87 90 91 93 95 97 106 107 109 112 123 131
135 136 142 143 144 145 150 151 152 153 158 160 165 166 173 174
176 177 183 185 188 193 196 198 200 203 204 211 223 236 244 246 248
 Unital No 16
|Aut(G,M65)|=  8
  9 10 11 12 16 18 21 22 61 62 65 68 69 75 77 91
 94 98 101 104 105 110 113 114 120 122 125 130 131 134 136 141
142 151 153 155 159 162 164 165 166 178 183 185 189 195 206 207
213 220 221 232 236 247 248 250 252 259 260 261 262 266 271 272 273
 Unital No 17
|Aut(G,M65)|=  8
  1  2  3  4 57 58 59 60 61 64 65 66 68 75 81 91
 95 98 99 101 109 110 111 117 120 126 127 128 130 131 134 136
140 146 151 155 156 159 160 165 168 174 178 179 189 206 207 208
210 220 231 236 238 244 246 247 250 258 259 260 261 262 263 264 265
 Unital No 18
|Aut(G,M65)|=  8
 14 34 35 36 50 51 52 56 61 62 63 67 69 76 77 94
 97 102 104 105 108 113 114 115 119 122 124 125 133 141 142 143
145 147 152 153 154 162 164 166 169 173 183 185 192 195 201 209
213 218 221 222 232 242 248 251 252 254 255 256 257 259 260 261 262
 Unital No 19
|Aut(G,M65)|=  8
  4  6 13 15 17 19 20 31 32 33 37 45 46 47 48 49
 50 51 52 56 65 70 73 74 75 80 82 88 89 96 102 111
115 118 122 125 128 138 139 140 145 146 147 148 157 165 168 169
170 178 184 185 186 187 190 207 216 221 234 238 239 240 241 243 263
 Unital No 20
|Aut(G,M65)|=  4
  2  6 39 53 65 66 67 68 70 74 75 76 80 81 82 88
 92 93 107 114 115 118 119 121 122 125 129 130 132 135 139 140
145 146 148 149 150 161 164 165 168 169 170 171 174 175 179 184
185 187 190 204 211 220 221 227 228 229 230 232 238 239 241 244 265
 Unital No 21
|Aut(G,M65)|=  8
  2  8 40 54 63 65 67 71 73 74 76 78 80 82 87 91
 92 95 101 111 112 115 116 121 122 123 125 129 130 134 143 144
145 146 149 151 152 158 162 164 165 168 170 171 177 179 183 185
187 192 193 212 213 220 221 223 224 229 231 234 236 239 240 251 263
 Unital No 22
|Aut(G,M65)|=  8
  2  8 37 40 46 47 48 50 51 52 54 56 66 67 68 73
 88 89 90 94 95 96 102 111 114 119 121 123 124 128 138 140
142 147 148 155 157 158 160 177 178 184 186 188 191 207 216 220
227 228 229 230 232 234 240 243 254 255 256 257 259 260 261 262 263
 Unital No 23
|Aut(G,M65)|=  8
  3  5 13 15 17 19 20 31 32 33 45 54 62 64 67 76
 80 83 84 87 90 93 95 110 113 121 123 125 130 135 136 137

```
142 143 150 151 153 154 158 160 164 168 171 173 174 177 183 187
188 189 218 219 220 229 235 237 265 266 267 268 269 270 271 272
273
 Unital No 24
 |Aut(G,M65)|=  8
 17  29  37  49  64  65  66  72  75  76  80  86  90  93  95  98
 99 107 109 110 111 115 120 127 129 130 131 133 135 140 143 147
151 152 154 155 163 167 169 173 176 178 179 187 188 191 194 202
205 206 208 209 217 218 223 225 226 230 231 236 244 245 249 251
261
 Unital No 25
 |Aut(G,M65)|=  8
 19  28  46  54  62  63  64  70  72  74  76  86  88  90  98 100
104 108 109 111 113 115 118 120 127 131 132 133 134 135 137 140
143 144 149 152 154 155 160 161 162 164 167 168 172 179 181 187
191 194 202 208 212 216 218 223 225 229 240 244 245 247 248 251
261
pp16htmmoor\jowk.DAT
 Unital No  1
 |Aut(G,M65)|=  32
  1   4  10  17  24  33  37  41  49  51  52  54  56  64  70  71
 73  74  78  79  81  83 103 104 105 109 110 112 118 119 132 136
137 144 145 152 155 157 159 161 163 168 173 180 185 190 194 196
206 207 208 213 216 223 224 225 230 233 243 247 254 258 259 262
271
 Unital No  2
 |Aut(G,M65)|=  16
  1   3   7  16  17  23  24  33  49  51  56  58  61  64  71  73
 75  78  79  88  90 102 104 105 107 108 109 119 120 121 125 132
136 137 142 145 148 152 161 172 175 179 180 181 185 194 205 206
208 212 223 224 228 230 236 237 242 244 246 252 254 256 259 268
271
 Unital No  3
 |Aut(G,M65)|=  48
  1   2   4   7  15  16  18  24  25  35  39  41  49  52  53  61
 81  82  86  87  93  97 100 108 111 115 117 131 132 133 135 138
144 155 156 158 159 160 168 171 174 183 186 195 196 206 212 215
217 220 225 229 235 237 247 249 250 255 260 262 264 268 270 272
273
 Unital No  4
 |Aut(G,M65)|=  24
  1  14  15  17  19  21  24  33  35  36  38  40  49  51  56  59
 61  63  66  77  79  80  90  92  94 103 104 105 108 117 119 121
125 126 132 140 141 144 158 171 175 181 185 187 189 191 194 202
204 207 210 215 220 225 228 229 230 234 237 239 243 245 248 249
254
 Unital No  5
 |Aut(G,M65)|=  12
  1  18  20  22  25  27  28  32  44  45  46  47  48  51  54  56
 59  63  64  66  68  70  71  77  80  81  86  96  98 104 110 111
117 119 120 123 125 132 142 145 149 155 161 172 175 179 188 194
195 200 202 203 204 207 212 213 217 219 220 225 229 230 246 250
254
 Unital No  6
 |Aut(G,M65)|=  8
 40  43  46  53  55  57  58  63  64  68  69  70  73  88  94 101
103 104 106 109 116 117 123 126 127 128 131 134 136 139 140 141
144 146 148 149 154 159 160 165 181 187 195 203 209 211 213 221
224 229 232 233 236 238 240 242 244 248 249 258 260 263 265 267
270
 Unital No  7
 |Aut(G,M65)|=  4
  8  11  20  29  30  40  43  44  45  46  47  59  64  66  67  72
 75  79  80  83  87  88  93  96 103 104 115 116 117 118 121 124
130 138 139 140 143 147 149 153 156 160 169 170 172 179 203 208
209 218 222 223 235 237 242 246 248 249 250 258 263 265 269 272
273
pp16htmmoor\lmrh.DAT
 Unital No  1
 |Aut(G,M65)|=  32
  1   3  11  27  30  31  36  39  46  56  57  61  62  72  80  83
 89  90  97 101 102 107 109 111 113 114 116 118 120 126 133 139
147 155 159 164 165 166 172 175 176 179 182 185 199 205 213 214
217 228 230 231 238 242 246 247 253 254 255 257 263 264 266 269
270
pp16htmmoor\math.DAT
```

Unital No  1
|Aut(G,M65)|=  16
 31  33  35  36  38  39  48  51  55  61  69  71  72  78  81  82
 89  90  94  99 101 109 111 113 114 115 116 117 121 131 134 137
140 144 149 150 154 156 161 162 174 178 181 185 189 195 214 218
224 229 230 231 232 238 241 243 244 255 260 263 265 268 270 272
273

Unital No  2
|Aut(G,M65)|=  128
 19  33  35  37  43  45  49  53  58  59  68  70  71  73  79  81
 82  86  94  95  99 103 107 108 115 119 122 124 131 134 137 141
142 149 150 152 153 154 160 161 173 177 179 184 189 194 201 203
207 210 235 238 240 244 250 251 257 258 260 261 266 269 270 271
273

Unital No  3
|Aut(G,M65)|=  128
 21  25  26  28  30  35  37  39  45  46  48  54  61  62  64  67
 69  72  84  86  87  88  90  91  96 101 110 117 120 131 134 137
151 152 153 155 160 164 166 172 175 180 187 193 201 203 204 207
214 215 217 222 223 226 227 228 230 232 250 253 255 257 258 262
273

Unital No  4
|Aut(G,M65)|=  64
 17  21  26  28  30  39  42  44  46  52  57  60  64  65  67  69
 74  77  78  83  84  88  90  91  96 101 104 116 123 126 130 136
147 151 155 156 158 164 170 172 175 180 187 191 193 199 202 204
206 209 214 215 222 223 225 226 227 229 231 252 253 255 256 263
273

Unital No  5
|Aut(G,M65)|=  128
 19  35  41  42  44  45  48  54  61  62  70  72  74  75  76  78
 81  94  95  99 114 117 123 126 131 134 137 139 140 142 147 153
157 158 159 162 163 166 169 170 171 178 181 182 186 189 196 198
199 200 202 203 205 206 207 209 211 225 228 229 241 246 247 252
273

Unital No  6
|Aut(G,M65)|=  16
 19  21  26  28  31  33  35  36  39  43  51  54  55  62  67  69
 70  72  78  86  87  88  91  95  96 105 109 110 111 112 116 117
120 122 132 135 137 142 151 157 160 161 167 169 172 173 174 188
190 196 212 217 222 229 237 238 244 246 249 251 259 263 264 268
273

Unital No  7
|Aut(G,M65)|=  128
 20  21  22  26  27  38  48  54  61  62  63  66  67  68  85  86
 88  89  92  97 102 104 107 111 125 132 135 138 139 140 145 148
154 156 158 160 163 169 170 175 180 190 191 197 198 199 201 202
214 216 221 224 231 234 237 241 245 246 247 255 257 260 268 272
273

Unital No  8
|Aut(G,M65)|=  64
 20  36  48  54  61  62  71  72  73  75  78  79  82  86  90  91
 96 101 106 112 114 115 117 118 121 126 132 135 138 139 144 149
150 157 160 161 162 166 167 169 171 174 176 183 185 186 196 198
200 201 219 220 225 228 229 238 240 246 249 257 264 265 269 270
273

Unital No  9
|Aut(G,M65)|=  128
 19  20  21  25  28  46  48  54  61  62  64  67  70  72  81  84
 87  94  95  99 108 110 111 112 117 119 120 122 132 135 138 142
144 146 148 155 166 168 172 173 175 180 189 193 194 195 204 215
217 226 227 228 230 232 234 235 243 253 262 264 265 266 268 271
273

Unital No 10
|Aut(G,M65)|=  16
 19  20  27  29  31  33  36  38  40  46  48  52  60  61  63  68
 71  78  79  87  95 102 107 109 114 119 127 128 138 139 144 148
156 158 162 164 166 170 171 174 175 180 193 194 197 198 200 201
215 217 223 224 227 228 229 231 234 235 244 257 265 266 269 270
273

Unital No 11
|Aut(G,M65)|=  128
 21  43  49  53  58  59  65  67  70  71  75  77  79  81  83  86
 90  91  96 101 110 111 114 115 116 120 121 130 136 139 141 142
148 149 157 160 161 162 167 169 171 175 180 185 186 189 196 198
200 201 230 232 234 238 246 249 251 256 257 261 263 268 269 270

273
 Unital No 12
 |Aut(G,M65)|=  64
 17  19  24  25  32  34  37  38  39  48  57  61  69  71  74  75
 78  81  89  95  98 113 116 129 132 133 135 140 141 147 153 169
176 178 181 183 186 189 193 194 199 201 202 206 207 209 215 218
224 229 230 232 233 235 239 241 246 250 257 261 262 263 270 272
273
 Unital No 13
 |Aut(G,M65)|=  16
 18  20  24  26  32  33  40  41  43  44  49  50  51  59  66  69
 73  76  82  84  86  89  95  98  99 100 105 111 112 118 123 125
128 129 135 140 143 145 147 148 163 164 165 174 176 178 187 190
197 201 205 206 210 212 214 217 219 222 226 236 261 262 265 267
273
 Unital No 14
 |Aut(G,M65)|=  16
 17  19  21  30  32  33  40  43  44  60  61  62  64  65  69  73
 79  80  83  87  93  98  99 103 106 109 114 115 116 118 122 124
134 136 137 141 145 148 150 151 164 165 168 172 175 177 186 188
192 193 194 195 203 204 207 209 226 227 238 242 262 268 270 272
273
 Unital No 15
 |Aut(G,M65)|=  32
 18  22  26  28  32  37  48  53  59  61  63  74  75  83  87  88
 89  94  99 102 103 104 105 113 119 127 128 131 134 144 146 154
155 158 161 164 168 170 172 178 181 184 186 188 191 205 206 211
212 217 218 221 223 226 238 239 242 250 256 259 260 265 266 272
273
 Unital No 16
 |Aut(G,M65)|=  16
 22  24  28  30  31  34  37  38  39  44  48  54  61  62  65  69
 73  75  77  82  83  84  93  98 100 103 104 105 109 114 116 130
136 147 150 152 155 156 162 172 177 179 184 186 187 191 192 193
204 209 213 215 221 224 226 231 233 240 250 252 256 258 259 263
273
pp16htmmoor\semi2.DAT
 Unital No  1
 |Aut(G,M65)|=  64
 27  28  29  31  32  37  38  42  45  46  50  52  55  59  66  67
 70  73  75  80  81  82  86  89  91  93  94  95 115 116 128 130
133 137 140 142 145 146 148 161 173 181 184 185 191 192 197 202
205 209 210 213 215 217 220 225 232 233 238 257 258 263 266 269
273
 Unital No  2
 |Aut(G,M65)|=  64
 18  24  30  31  32  33  36  37  38  46  50  52  55  59  67  68
 75  78  79  80  81  86  91  93  94  95 104 109 111 112 115 119
128 130 133 137 143 144 146 148 151 152 157 164 168 181 184 189
192 201 202 232 233 234 235 236 239 240 242 250 254 261 265 269
273
 Unital No  3
 |Aut(G,M65)|=  64
 17  19  20  28  31  34  37  44  45  48  50  52  64  66  71  72
 81  88  89  92  96 104 109 110 113 128 130 131 139 140 141 147
154 157 165 167 171 191 192 194 196 202 203 208 210 214 215 220
221 232 233 234 238 241 243 245 248 249 255 257 258 267 269 272
273
 Unital No  4
 |Aut(G,M65)|=  64
 17  18  20  26  32  33  36  37  38  48  54  57  74  75  78  79
 80  81  91  93  94  95 103 110 111 115 124 131 134 136 143 148
151 154 156 164 168 171 172 175 181 184 188 189 193 194 196 198
204 224 227 232 233 236 239 241 246 248 249 250 252 254 261 262
273
 Unital No  5
 |Aut(G,M65)|=  64
 24  41  50  52  54  57  66  84  92  96  97 102 105 107 111 112
114 117 119 120 123 124 128 130 134 136 140 151 152 153 155 158
160 169 175 180 182 186 187 190 194 196 201 212 218 222 227 229
230 235 236 238 240 241 242 244 245 248 257 260 261 262 270 272
273
 Unital No  6
 |Aut(G,M65)|=  192
 19  22  29  31  32  33  39  40  42  46  51  52  55  58  64  66
 74  75  76  80  84  89  90 100 102 104 106 109 113 114 124 126

```
 128 132 135 137 161 162 163 175 176 178 182 193 195 201 205 207
 209 220 229 231 232 236 242 247 249 251 252 253 257 258 260 272
 273
 Unital No  7
 |Aut(G,M65)|=  192
  18  21  25  31  32  46  51  52  55  58  72  73  74  75  76  81
  84  85  89  91  92  94  99 114 116 119 120 126 128 132 135 137
 140 150 156 157 166 169 172 175 181 191 194 195 205 206 209 211
 215 217 218 230 249 251 252 254 255 258 260 262 263 266 268 272
 273
 Unital No  8
 |Aut(G,M65)|=  48
  18  21  27  31  32  46  51  52  55  58  66  72  73  75  80  83
  84  88  89  90  91  94 100 114 118 119 123 126 128 132 135 137
 150 153 157 161 166 169 170 174 180 183 196 205 209 210 211 214
 215 217 218 232 237 243 246 249 251 252 254 258 259 263 264 268
 273
 Unital No  9
 |Aut(G,M65)|=  192
  20  23  25  29  32  34  35  36  39  42  51  52  55  58  63  65
  67  74  75  76  81  85  86  89  97  98 104 105 106 117 119 124
 136 147 150 157 158 164 167 169 175 181 191 193 195 200 202 205
 207 209 216 218 220 242 244 253 257 258 260 266 268 269 271 272
 273
 Unital No 10
 |Aut(G,M65)|=  48
  20  21  22  29  32  33  36  40  41  42  51  52  55  58  63  64
  72  75  77  83  86  88  89  91  94  98 104 109 113 117 118 122
 123 127 136 153 158 164 166 170 171 174 176 178 182 200 201 204
 205 209 211 214 216 220 236 237 242 243 247 257 258 259 264 269
 273
 Unital No 11
 |Aut(G,M65)|=  48
  18  23  26  27  30  34  35  39  43  47  51  52  55  58  65  66
  67  73  80  90  97 100 105 106 110 111 115 119 121 124 131 134
 139 147 148 149 150 152 157 161 167 169 180 183 193 196 202 207
 210 215 217 218 221 232 239 244 245 246 253 254 263 268 270 271
 273
 Unital No 12
 |Aut(G,M65)|=  16
  21  23  25  28  31  37  40  41  42  47  48  49  53  57  63  68
  72  73  76  84  86  89  94  95 100 108 111 112 121 122 125 131
 132 136 138 144 153 154 157 162 163 166 169 179 181 185 189 201
 204 205 210 213 214 218 229 234 238 240 242 245 251 255 266 272
 273
 Unital No 13
 |Aut(G,M65)|=  64
  23  26  28  30  31  44  50  55  60  62  63  66  72  73  75  79
  81  85  86  89  92 100 101 106 108 111 115 116 131 132 136 138
 149 153 164 174 176 182 183 184 185 187 189 190 193 205 208 212
 213 216 222 229 234 235 242 246 247 248 249 255 256 262 267 268
 273
 Unital No 14
 |Aut(G,M65)|=  64
  17  18  21  23  24  35  41  44  47  48  54  57  65  70  72  76
  77  84  88  92  96  99 101 104 107 111 113 114 121 122 126 127
 143 145 147 150 162 163 166 169 173 174 176 179 189 211 214 221
 225 226 228 234 239 242 244 245 250 251 254 261 264 267 271 272
 273
 Unital No 15
 |Aut(G,M65)|=  64
  23  24  28  30  32  33  35  36  41  44  50  52  55  59  64  65
  72  75  77  78  79  84  88  92  96  99 101 104 106 109 114 115
 119 122 126 127 141 147 148 151 152 157 164 168 174 177 178 197
 199 201 205 211 226 234 235 236 242 244 251 253 256 258 265 271
 273
 Unital No 16
 |Aut(G,M65)|=  64
  26  38  41  42  44  47  51  54  55  62  64  67  68  69  71  81
  90  94  99 102 103 104 105 109 113 116 120 125 127 130 131 138
 142 148 150 154 155 159 164 166 182 185 186 195 206 207 210 211
 215 216 225 235 236 241 244 247 248 253 254 258 261 262 264 269
 273
 Unital No 17
 |Aut(G,M65)|=  64
  18  21  26  27  30  34  35  37  39  41  50  51  55  56  80  90
```

```
  98 102 108 118 123 124 127 131 136 140 145 149 153 156 159 161
 163 169 170 171 179 181 186 190 191 195 197 198 203 212 220 221
 225 229 233 234 235 242 246 248 249 250 251 255 258 263 264 268
 273
 Unital No 18
 |Aut(G,M65)|=  64
  18  20  24  29  32  37  48  52  53  58  65  71  72  78  83  84
  86  90  94  99 102 103 106 108 109 128 135 141 146 160 161 163
 168 170 172 174 179 180 185 190 193 195 198 203 217 218 219 225
 228 230 233 235 240 247 248 252 255 257 258 261 264 266 271 272
 273
 Unital No 19
 |Aut(G,M65)|=  32
  18  20  21  23  27  34  36  37  41  46  51  57  59  62  66  68
  71  73  78  80  85  88  89  91  96 100 104 106 111 115 118 127
 128 138 141 145 146 150 155 165 176 182 188 190 201 203 204 206
 216 218 220 221 223 227 238 241 243 252 258 259 262 265 271 272
 273
 Unital No 20
 |Aut(G,M65)|=  32
  30  34  37  40  46  47  48  51  59  60  68  71  72  75  77  78
  80  88  93  99 100 103 104 106 110 111 117 118 119 120 125 130
 135 136 140 146 149 150 154 155 160 162 164 172 173 176 178 183
 191 193 197 203 209 213 214 216 223 238 248 250 251 260 263 271
 273
 Unital No 21
 |Aut(G,M65)|=  64
  31  36  37  38  45  47  48  57  60  62  67  68  73  75  82  85
  86  88  91  93  96  99 103 105 107 108 115 126 128 132 133 138
 150 157 158 165 168 172 175 181 192 197 199 201 209 210 212 213
 219 226 229 232 240 241 244 247 253 257 258 259 260 266 268 272
 273
pp16htmmoor\bbh2.DAT
 Unital No  1
 |Aut(G,M65)|=  16
  21  23  24  25  26  31  32  33  34  35  36  37  38  56  57  59
  94 102 103 105 107 111 113 115 130 133 135 137 148 150 157 159
 170 174 176 180 182 184 186 188 189 190 192 195 201 206 210 211
 214 222 226 236 238 239 242 244 248 253 257 260 261 262 267 268
 273
 Unital No  2
 |Aut(G,M65)|=  16
  27  28  29  30  48  49  51  53  62  67  68  70  73  75  76  79
  93  94  95  96  99 105 107 111 130 133 135 137 144 145 150 159
 165 167 170 176 180 182 184 186 190 192 194 196 201 205 207 210
 221 222 226 229 238 239 242 244 248 249 252 257 261 262 267 268
 273
 Unital No  3
 |Aut(G,M65)|=  16
   1   2   3   4   5   7   9  11  40  42  43  45  47  50  52  54
  98 100 102 103 105 111 113 115 121 125 133 137 141 148 151 157
 172 174 175 176 180 182 186 188 189 190 192 195 206 208 209 211
 214 222 226 236 237 240 241 244 248 253 258 260 261 262 269 272
 273
 Unital No  4
 |Aut(G,M65)|=  16
   6   8  10  12  15  18  19  20  61  63  64  65  66  69  71  72
  93  95  96  98  99 100 105 111 121 125 133 137 141 144 145 151
 165 167 172 175 176 180 182 186 190 192 194 196 205 207 208 209
 221 222 226 229 237 240 241 244 248 249 252 258 261 262 269 272
 273
 Unital No  5
 |Aut(G,M65)|=   8
  27  28  29  30  31  33  34  37  40  42  43  45  61  69  71  72
  89  90  91  92  95  99 103 106 108 112 113 114 118 120 121 122
 124 125 126 130 132 135 136 139 143 144 145 146 148 153 154 157
 178 180 182 187 243 245 247 256 261 262 263 267 268 269 271 272
 273
 Unital No  6
 |Aut(G,M65)|=  32
  21  27  28  29  30  32  35  36  38  39  41  44  46  56  57  59
  82  84  90  91 104 109 110 116 126 127 129 131 132 136 139 140
 141 148 151 157 172 174 175 176 180 182 186 188 194 196 199 200
 202 203 205 207 224 227 230 232 238 239 242 245 250 251 256 257
 273
 Unital No  7
```

```
 |Aut(G,M65)|=  32
 61  63  64  65  66  69  71  72  73  74  75  76  77  78  79  80
 82  84  90  91 104 109 110 116 118 119 122 123 131 132 139 140
141 144 145 151 165 167 172 175 176 180 182 186 189 195 199 200
202 203 206 211 224 228 231 232 238 239 242 245 254 256 257 259
273
 Unital No  8
 |Aut(G,M65)|=   4
  4   6   7  11  12  15  16  17  19  23  25  27  28  34  38  41
 45  61  62  68  71  73  74  77  79  85  90  92  93  97 105 110
112 115 120 122 131 132 136 138 143 149 151 152 153 156 165 175
177 178 180 188 196 201 205 209 220 224 232 234 237 240 246 255
264
 Unital No  9
 |Aut(G,M65)|=   4
 11  12  17  19  22  33  36  39  43  47  51  59  61  63  64  65
 66  68  73  77  83  84  88  97 100 106 107 110 111 117 118 119
120 125 126 127 132 135 140 141 146 149 154 160 161 167 173 174
181 184 186 189 190 191 211 222 226 228 231 244 248 249 252 265
273
 Unital No 10
 |Aut(G,M65)|=  20
  1  11  12  17  19  21  24  26  29  30  33  34  36  38  39  41
 43  45  48  50  52  53  55  56  58  83  84  88  96  97 102 110
114 116 120 123 128 129 132 137 150 164 173 175 190 191 194 197
199 207 215 217 228 231 242 244 245 248 249 252 256 257 265 268
272
 Unital No 11
 |Aut(G,M65)|=   8
  4   9  10  13  18  84  85  88  89  92  95  96 102 103 105 108
114 116 122 123 128 129 131 134 136 137 149 151 175 177 189 191
194 197 198 202 207 208 209 211 215 216 222 224 226 231 232 234
235 236 239 240 241 242 244 245 248 250 252 256 263 264 265 269
272
 Unital No 12
 |Aut(G,M65)|=  16
  1  35  37  40  44  87  88  89  90  93  99 105 109 112 113 114
115 122 123 128 129 131 133 136 138 144 150 154 155 161 164 167
173 175 182 183 186 189 195 197 200 201 206 208 211 213 215 216
219 220 222 223 226 237 238 241 242 243 244 247 248 265 267 270
272
 Unital No 13
 |Aut(G,M65)|=  16
  4   9  10  13  18  84  86  87  90  91  93  99 106 109 111 112
113 115 117 118 119 126 127 133 138 140 149 151 175 177 191 195
196 199 200 201 202 205 206 210 213 217 218 219 220 223 225 231
233 236 237 238 243 246 247 250 252 255 257 258 263 266 267 268
270
 Unital No 14
 |Aut(G,M65)|=  16
  1  35  37  40  44  85  86  91  92  95  96 102 103 106 108 111
116 117 118 119 126 127 134 137 140 144 150 154 155 161 164 167
173 175 182 183 186 194 196 198 199 205 207 209 210 217 218 224
225 232 233 234 235 239 240 245 246 255 256 257 258 264 266 268
269
 Unital No 15
 |Aut(G,M65)|=   8
 15  18  19  20  47  50  52  54  63  64  65  66  74  77  78  80
 83  85  86  90  92 102 108 111 115 124 125 126 128 130 132 133
136 143 146 149 152 158 160 170 174 175 178 181 188 191 194 196
212 218 221 223 226 227 229 230 235 238 241 244 247 263 268 272
273
 Unital No 16
 |Aut(G,M65)|=  10
  1  11  12  17  19  28  30  39  42  44  45  48  50  51  54  62
 63  66  68  69  72  73  74  76  80  86  88  90  93  96  98 106
107 109 119 123 125 130 138 140 141 143 153 159 168 170 172 176
178 179 181 182 190 192 194 195 200 201 229 231 249 259 261 262
267
 Unital No 17
 |Aut(G,M65)|=  16
  5   7   9  11  15  18  19  20  22  23  24  25  26  55  58  60
 81  84  86  88  96  98 102 107 108 109 111 114 122 124 132 136
142 145 148 149 151 155 158 159 168 182 183 187 189 196 197 199
217 225 226 234 239 245 247 250 252 254 258 260 262 267 269 271
273
```

Unital No 18
 |Aut(G,M65)|=  8
 20  24  33  61  74  87  88  91  92  95  99 100 102 104 107 116
119 124 127 131 132 133 134 138 143 150 152 154 155 157 163 165
167 171 172 173 180 182 184 192 193 194 195 202 203 208 210 217
219 225 230 232 233 234 239 241 243 245 252 254 259 260 261 264
268
 Unital No 19
 |Aut(G,M65)|=  8
 20  24  33  61  74  85  87  89  92  94  96  99 101 107 109 113
117 119 124 129 133 137 138 139 142 144 146 151 152 154 161 165
170 171 172 174 177 182 187 191 193 195 196 201 203 208 212 215
218 224 225 226 228 234 238 239 243 247 249 250 254 260 262 269
270
 Unital No 20
 |Aut(G,M65)|=  4
 20  24  33  61  74  83  84  89  92  94  97 107 111 112 114 115
117 118 119 121 122 129 135 138 141 146 148 149 154 158 160 169
171 177 179 182 183 186 187 195 196 197 199 203 206 211 212 213
214 216 225 226 227 236 238 239 240 244 246 249 257 260 266 271
272
 Unital No 21
 |Aut(G,M65)|=  4
 18  23  31  61  78  82  83  89  90  94  98 102 106 110 111 112
117 118 119 121 125 127 128 136 143 147 154 157 158 159 160 166
169 176 177 178 183 185 187 194 195 199 200 204 207 211 212 213
221 222 226 230 235 236 237 238 240 244 252 255 258 260 264 267
271
 Unital No 22
 |Aut(G,M65)|=  8
 13  14  15  16  17  18  19  20  27  28  29  30  62  67  68  70
 86  87  89  90  94  97  98  99 103 110 111 116 117 118 134 136
141 159 160 163 165 166 177 180 183 185 187 188 202 207 211 212
213 216 221 226 232 233 235 236 238 251 254 258 261 267 271 272
273
 Unital No 23
 |Aut(G,M65)|=  8
 22  27  28  29  30  31  33  34  37  40  42  43  45  55  58  60
 81  83  87  88  95  99 104 108 109 110 112 116 121 125 127 129
143 146 147 150 158 159 161 163 165 167 169 181 194 196 199 200
219 224 232 235 237 240 244 248 249 252 254 259 261 262 267 268
273
 Unital No 24
 |Aut(G,M65)|=  8
 21  27  28  29  30  47  50  52  54  56  57  59  74  77  78  80
 85  86  90  91  98 100 103 104 109 110 113 116 117 126 128 136
150 159 161 163 166 168 173 174 178 185 187 188 204 205 207 212
219 221 222 224 226 229 232 235 241 250 251 258 261 262 267 268
273
 Unital No 25
 |Aut(G,M65)|=  8
 32  35  36  38  39  41  44  46  63  64  65  66  73  75  76  79
 83  84  87  88  89  93  96  98 104 117 121 127 129 131 133 134
139 148 150 153 154 157 158 160 163 169 170 185 186 198 205 207
208 214 218 226 233 234 236 238 242 243 245 250 251 261 269 271
273
pp16htmmoor\BBH1.DAT
 Unital No  1
 |Aut(G,M65)|=  8
 14  20  33  39  40  41  43  46  48  50  52  53  55  60  63  64
 68  71  74  79  81  83  86  87  89  94  95  99 107 109 112 113
114 117 127 130 134 139 142 148 150 162 170 171 173 175 180 188
194 196 197 199 207 214 222 226 243 248 249 252 258 261 262 263
264
 Unital No  2
 |Aut(G,M65)|=  8
 14  30  45  67  75  83  86  88  92  93  94  95  99 100 102 111
113 114 119 120 125 127 128 129 131 137 138 139 148 150 152 155
161 170 171 173 175 178 182 188 190 194 198 199 200 204 205 206
214 215 226 230 233 237 243 244 248 249 252 258 262 263 264 269
272
 Unital No  3
 |Aut(G,M65)|=  8
 17  22  27  31  35  39  43  46  52  55  57  60  63  64  65  67
 68  83  86  89  92  93  95  99 104 107 108 109 110 112 117 118
119 134 137 139 142 153 160 162 168 173 175 184 190 195 196 197

```
 199 200 207 212 220 222 225 226 233 237 251 252 254 258 262 266
 268
 Unital No  4
 |Aut(G,M65)|=  8
 13  14  15  16  17  22  27  31  35  39  43  46  52  55  60  63
 68  84  85  89  90  93  94 102 104 106 109 110 119 122 124 125
 127 129 138 139 149 150 152 155 168 172 173 182 188 191 195 196
 200 201 204 211 220 221 222 225 233 239 242 245 246 252 262 266
 268
 Unital No  5
 |Aut(G,M65)|=  8
  2   3   4   6  17  22  27  30  31  33  34  35  36  38  42  44
 47  84  87  90  92  93 100 111 112 115 119 121 122 124 128 130
 132 134 137 139 142 144 147 161 169 172 173 174 178 179 180 190
 197 200 205 211 217 233 237 242 244 245 252 255 259 261 262 267
 269
 Unital No  6
 |Aut(G,M65)|=  32
 17  22  27  31  35  39  43  46  52  55  57  60  63  64  65  67
 68  86  96  99 102 107 109 110 112 118 119 123 125 131 134 137
 139 141 145 154 166 173 174 175 196 200 203 205 207 208 209 212
 215 216 217 219 220 225 235 236 237 240 241 244 251 254 256 258
 266
pp16htmmoor\bbs4.DAT
 Unital No  1
 |Aut(G,M65)|=  8
  3   9  10  12  16  18  19  20  57  61  63  65  67  70  71  72
 88  91  92  93  95 100 105 109 111 113 116 125 129 131 133 140
 141 142 145 146 147 150 152 153 154 159 163 168 178 191 202 209
 210 211 213 222 224 227 234 239 246 252 255 259 260 266 267 270
 273
 Unital No  2
 |Aut(G,M65)|=  8
 16  18  19  20  28  30  34  36  58  60  62  69  73  74  75  76
 83  84  86  90  91  97 101 104 106 109 116 118 135 138 141 145
 146 148 152 156 160 168 170 171 175 178 179 181 188 198 199 200
 205 208 214 223 227 231 232 239 244 250 251 257 258 259 263 268
 273
 Unital No  3
 |Aut(G,M65)|=  8
 14  15  17  24  38  44  46  47  53  54  55  56  58  60  62  69
 82  84  89  91  92  97  98 102 108 109 115 124 139 141 146 149
 152 158 161 162 164 169 170 172 176 179 181 183 186 188 189 192
 212 219 221 224 226 227 232 239 244 245 247 248 252 259 264 267
 273
 Unital No  4
 |Aut(G,M65)|=  4
  4  21  28  54  61  82  86  89  91  92  94  98 100 103 106 108
 109 111 116 124 129 132 149 152 158 159 160 161 162 164 168 169
 172 174 175 179 181 183 185 186 187 192 197 199 205 206 208 209
 212 221 224 232 234 235 239 241 244 245 248 250 252 257 264 267
 268
 Unital No  5
 |Aut(G,M65)|=  4
  4  23  34  55  65  82  86  87  88  89  93  94  97 103 105 107
 108 109 110 125 129 131 132 140 146 149 150 154 160 164 165 166
 174 180 181 182 183 186 188 190 191 192 197 199 203 217 218 220
 221 227 232 233 234 235 239 240 242 245 246 248 254 256 262 264
 268
 Unital No  6
 |Aut(G,M65)|=  8
 25  26  28  30  31  32  34  36  40  41  42  48  58  60  62  69
 82  86  88  90  91  94  95  97 105 108 109 124 129 136 140 146
 152 156 160 161 164 165 179 181 182 184 186 188 192 197 199 200
 203 221 222 224 227 232 234 239 244 246 248 256 263 264 267 268
 273
 Unital No  7
 |Aut(G,M65)|=  12
 37  42  44  51  53  82  85  86  87  90  98 101 103 105 108 115
 118 122 125 127 129 133 135 136 138 152 158 162 164 165 170 172
 175 179 180 181 187 188 193 197 200 202 203 204 209 216 220 221
 222 224 225 228 229 230 232 233 239 240 242 243 254 258 263 266
 272
 Unital No  8
 |Aut(G,M65)|=  12
 40  43  47  51  55  83  97  99 100 102 104 109 111 112 113 116
```

118 120 123 125 130 131 135 138 139 147 148 149 150 151 156 157
158 162 163 168 172 173 181 182 184 185 186 190 192 193 200 211
212 213 241 245 246 247 249 254 255 256 260 261 264 265 269 270
272
 Unital No  9
 |Aut(G,M65)|=   24
  2  5  6 11 16 18 19 20 28 30 34 36 53 54 55 56
 81 85 86 90 107 112 113 116 117 122 123 126 129 130 131 135
138 139 143 151 152 158 162 163 165 170 174 175 182 183 185 187
188 192 196 198 216 217 219 220 221 222 223 224 229 230 232 233
273
 Unital No 10
 |Aut(G,M65)|=   4
 22 31 37 57 69 81 83 84 91 93 97 100 104 105 109 111
114 117 119 121 127 132 134 135 139 141 143 149 150 151 154 155
168 169 173 180 185 186 192 193 194 196 199 212 214 216 217 218
223 225 226 231 232 236 240 245 246 247 248 254 258 260 261 268
272
 Unital No 11
 |Aut(G,M65)|=   4
 13 25 37 58 70 84 100 109 121 132 134 135 136 139 142 146
147 149 152 155 157 158 159 167 169 170 172 173 176 179 180 182
184 186 188 192 199 209 210 212 213 214 216 217 218 224 226 230
233 234 235 239 241 245 247 248 255 256 261 264 267 268 270 271
272
 Unital No 12
 |Aut(G,M65)|=   4
 23 31 39 60 70 84 100 109 121 129 131 132 138 139 143 144
147 149 151 153 154 155 157 165 167 168 169 173 176 178 181 192
193 196 198 199 200 202 204 205 206 207 214 216 217 218 223 225
231 234 235 238 240 243 247 248 254 256 257 258 259 261 262 264
272
 Unital No 13
 |Aut(G,M65)|=   8
  3  9 10 12 38 40 41 42 44 46 47 48 58 60 62 69
131 132 134 136 138 139 142 143 149 150 155 157 169 172 173 176
178 180 182 192 194 196 198 199 204 205 212 214 216 217 218 224
231 233 235 236 243 245 246 247 248 255 258 261 262 264 270 272
273

**Appendix 3: Unitals of projective planes of order 16, found with author's algorithm on dreanaut format of the planes [10]**

pp16\pp-16-1.DAT,  |Aut(G546)|=  34217164800
 Unital No  1
 |Aut(G,M65)|=  249600
 17 29 39 41 43 48 49 50 52 54 55 63 75 84 89 90
 95 97 98 103 107 109 110 114 115 116 124 128 130 138 142 144
145 154 163 170 173 174 175 178 179 187 191 193 195 196 201 203
205 212 215 218 221 224 241 243 248 249 254 256 261 263 265 270
271
 Unital No  2
 |Aut(G,M65)|=  768
  1 21 22 27 31 34 36 39 46 50 53 61 62 72 77 79
 81 86 87 93 95 100 101 105 110 116 119 120 122 131 133 137
138 146 149 151 152 162 167 168 170 180 182 186 193 194 200 201
207 214 217 220 221 230 238 239 241 250 253 254 257 260 265 272
273
pp16\pp-16-2.DAT,  |Aut(G546)|=  147456
 Unital No  1
 |Aut(G,M65)|=   64
  7  9 11 12 19 22 29 30 33 34 40 47 55 57 59 60
 68 69 74 80 83 86 93 94 97 98 104 111 115 118 125 126
135 137 139 140 145 146 152 159 167 169 171 172 177 178 184 191
193 194 200 207 211 214 221 222 231 233 235 236 243 246 253 254
273
 Unital No  2
 |Aut(G,M65)|=   64
  3  6 13 14 20 21 26 32 35 38 45 46 49 50 56 63
 68 69 74 80 81 82 88 95 100 101 106 112 115 118 125 126
135 137 139 140 145 146 152 159 163 166 173 174 177 178 184 191
193 194 200 207 211 214 221 222 228 229 234 240 244 245 250 256
273
 Unital No  3
 |Aut(G,M65)|=   64
  3  6 13 14 19 22 29 30 36 37 42 48 52 53 58 64
 68 69 74 80 87 89 91 92 100 101 106 112 113 114 120 127
132 133 138 144 145 146 152 159 163 166 173 174 179 182 189 190
193 194 200 207 211 214 221 222 225 226 232 239 241 242 248 255

273
 Unital No  4
|Aut(G,M65)|=  192
  1   8   9  12  20  21  22  30  33  40  41  44  50  55  59  63
 65  72  73  76  84  85  86  94 100 101 102 110 113 120 121 124
130 135 139 143 148 149 150 158 162 167 171 175 178 183 187 191
196 197 198 206 211 218 221 224 226 231 235 239 241 248 249 252
273
 Unital No  5
|Aut(G,M65)|=  4
  4   6  11  15  18  21  23  30  36  38  43  47  53  61  62  64
 75  77  79  80  84  86  91  95 101 109 110 112 113 120 121 124
130 131 135 138 147 148 150 154 162 163 167 170 181 189 190 192
196 197 198 206 211 218 221 224 226 231 235 239 242 243 247 250
273
 Unital No  6
|Aut(G,M65)|=  192
  1   8   9  12  18  23  27  31  34  39  43  47  50  55  59  63
 67  74  77  80  81  88  89  92  99 106 109 112 115 122 125 128
129 136 137 140 146 151 155 159 163 170 173 176 179 186 189 192
193 200 201 204 212 213 214 222 226 231 235 239 241 248 249 252
273
 Unital No  7
|Aut(G,M65)|=  48
  4   6  11  15  21  24  25  30  40  41  43  47  53  61  62  64
 65  76  77  80  81  88  89  92  98 101 103 110 117 125 126 128
129 136 137 140 146 151 155 159 161 162 167 172 187 189 191 192
193 204 205 208 213 216 217 222 226 231 235 239 242 243 247 250
273
 Unital No  8
|Aut(G,M65)|=  192
  1   8   9  12  18  23  27  31  33  40  41  44  50  55  59  63
 65  72  73  76  82  87  91  95  99 106 109 112 116 117 118 126
131 138 141 144 147 154 157 160 163 170 173 176 178 183 187 191
195 202 205 208 209 216 217 220 226 231 235 239 241 248 249 252
273
 Unital No  9
|Aut(G,M65)|=  48
  1   8   9  12  17  18  23  28  33  36  38  44  50  55  59  63
 72  73  75  79  84  86  91  95 107 109 111 112 115 116 118 122
130 131 135 138 147 148 150 154 163 168 169 170 178 181 183 190
195 200 201 202 209 212 214 220 226 231 235 239 241 248 249 252
273
 Unital No  10
|Aut(G,M65)|=  64
  3   5  10  14  18  23  24  25  33  43  44  47  52  54  61  64
 67  69  74  78  82  87  88  89 100 102 109 112 114 119 120 121
131 133 138 142 146 151 152 153 164 166 173 176 180 182 189 192
195 197 202 206 211 213 218 222 228 230 237 240 242 247 248 249
273
 Unital No  11
|Aut(G,M65)|=  64
  3   5  10  14  17  27  28  31  36  38  45  48  51  53  58  62
 67  69  74  78  84  86  93  96  99 101 106 110 113 123 124 127
129 139 140 143 146 151 152 153 161 171 172 175 180 182 189 192
196 198 205 208 211 213 218 222 228 230 237 240 241 251 252 255
273
 Unital No  12
|Aut(G,M65)|=  16
  1  11  12  15  17  27  28  31  34  39  40  41  51  59  60  62
 67  71  73  78  84  86  93  96  99 107 108 110 114 120 123 124
131 135 137 142 146 148 152 160 161 164 175 176 178 180 184 192
195 196 206 208 209 215 217 223 225 231 233 239 242 245 248 250
273
 Unital No  13
|Aut(G,M65)|=  64
  2   3   5   9  17  22  27  32  39  40  42  46  55  56  58  62
 66  67  69  73  82  83  85  89  97 102 107 112 114 115 117 121
135 136 138 142 145 150 155 160 164 172 173 175 177 182 187 192
199 200 202 206 209 214 219 224 226 227 229 233 247 248 250 254
273
 Unital No  14
|Aut(G,M65)|=  64
  3   4   9  15  18  21  28  29  34  37  44  45  54  56  58  59
 67  68  73  79  83  84  89  95  99 100 105 111 114 117 124 125
134 136 138 139 146 149 156 157 163 164 169 175 177 183 190 192

```
198 200 202 203 214 216 218 219 226 229 236 237 246 248 250 251
273
 Unital No  15
|Aut(G,M65)|=  64
  3  11  12  14  22  23  25  29  38  39  41  45  50  52  56  64
 66  68  72  80  83  91  92  94  98 100 104 112 115 123 124 126
129 133 138 143 150 151 153 157 166 167 169 173 182 183 185 189
194 196 200 208 210 212 216 224 227 235 236 238 243 251 252 254
273
 Unital No  16
|Aut(G,M65)|=  64
  6   7  11  14  18  20  21  31  33  40  42  48  51  57  60  61
 65  72  74  80  82  84  85  95  97 104 106 112 115 121 124 125
129 136 138 144 146 148 149 159 163 169 172 173 177 184 186 192
194 196 197 207 211 217 220 221 227 233 236 237 242 244 245 255
273
 Unital No  17
|Aut(G,M65)|=  64
  7   9  11  12  23  25  27  28  35  38  45  46  55  57  59  60
 68  69  74  80  84  85  90  96 103 105 107 108 115 118 125 126
132 133 138 144 148 149 154 160 161 162 168 175 179 182 189 190
195 198 205 206 212 213 218 224 227 230 237 238 247 249 251 252
273
 Unital No  18
|Aut(G,M65)|=  64
  3  10  11  15  19  26  27  31  36  38  40  41  51  58  59  63
 66  71  77  80  83  90  91  95 100 102 104 105 113 117 124 126
129 133 140 142 145 149 156 158 164 166 168 169 177 181 188 190
195 202 203 207 212 214 216 217 225 229 236 238 244 246 248 249
273
 Unital No  19
|Aut(G,M65)|=  32
  1   7  14  16  22  24  26  27  33  39  46  48  54  56  58  59
 66  69  76  77  86  88  90  91  97 103 110 112 114 117 124 125
129 135 142 144 147 148 153 159 161 167 174 176 178 181 188 189
198 200 202 203 210 213 220 221 230 232 234 235 242 245 252 253
273
 Unital No  20
|Aut(G,M65)|=  64
  3   8  12  16  21  22  23  31  37  38  39  47  49  57  58  61
 67  72  76  80  81  89  90  93 101 102 103 111 113 121 122 125
131 136 140 144 147 152 156 160 161 169 170 173 181 182 183 191
193 201 202 205 210 212 219 222 229 230 231 239 243 248 252 256
273
 Unital No  21
|Aut(G,M65)|=  32
  5   8  11  13  18  22  26  28  41  46  47  48  53  56  59  61
 66  70  74  76  89  94  95  96 105 110 111 112 117 120 123 125
137 142 143 144 153 158 159 160 162 166 170 172 178 182 186 188
197 200 203 205 213 216 219 221 226 230 234 236 241 243 244 247
273
pp16\pp-16-3a
 Unital No   1
|Aut(G,M65)|=  64
  4   5  10  16  17  18  24  31  39  41  43  44  52  53  58  64
 71  73  75  76  81  82  88  95 103 105 107 108 116 117 122 128
135 137 139 140 148 149 154 160 163 166 173 174 180 181 186 192
193 194 200 207 209 210 216 223 231 233 235 236 241 242 248 255
273
 Unital No   2
|Aut(G,M65)|=  128
  7   9  11  12  19  22  29  30  36  37  42  48  52  53  58  64
 71  73  75  76  83  86  93  94 103 105 107 108 119 121 123 124
132 133 138 144 148 149 154 160 163 166 173 174 180 181 186 192
195 198 205 206 209 210 216 223 231 233 235 236 243 246 253 254
273
 Unital No   3
|Aut(G,M65)|=   8
  4   7   9  16  18  22  24  29  36  39  41  48  54  59  60  61
 65  69  74  79  82  86  88  93  97 100 111 112 117 119 121 122
133 135 137 138 146 148 152 160 162 164 168 176 177 181 186 191
193 198 205 207 210 213 216 218 227 231 233 238 241 246 253 255
273
 Unital No   4
|Aut(G,M65)|=   8
 21  24  31  32  41  42  44  48  57  59  61  62  66  71  77  80
```

```
 84  92  93  95  97 103 105 111 116 117 119 123 131 132 142 144
153 155 157 158 164 172 173 175 177 181 188 190 196 197 199 203
209 213 220 222 225 231 233 239 241 246 251 256 259 262 263 268
273
 Unital No   5
 |Aut(G,M65)|=  48
  1   2   4  10  21  25  29  31  37  40  47  48  49  58  59  62
 70  71  72  80  88  91  94  96  99 101 108 111 113 114 116 122
131 139 140 142 145 149 154 159 161 170 171 174 179 187 188 190
197 200 207 208 210 211 212 220 226 227 228 236 242 244 248 256
273
 Unital No   6
 |Aut(G,M65)|=  8
  9  14  15  16  17  21  23  29  33  35  36  39  49  51  52  55
 69  77  78  80  83  84  88  91 105 110 111 112 117 125 126 128
131 132 142 144 146 153 156 159 162 163 164 172 178 185 188 191
193 198 199 202 210 213 220 221 225 231 233 239 242 245 252 253
273
 Unital No   7
 |Aut(G,M65)|=  192
  1   3   8  13  17  19  24  29  36  37  39  43  52  53  55  59
 65  67  72  77  89  90  92  96 100 101 103 107 114 118 126 127
137 138 140 144 148 149 151 155 161 163 168 173 185 186 188 192
201 202 204 208 209 211 216 221 233 234 236 240 244 245 247 251
273
 Unital No   8
 |Aut(G,M65)|=  128
  1   6  11  16  18  19  21  25  33  38  43  48  50  51  53  57
 65  70  75  80  81  86  91  96  98  99 101 105 116 124 125 127
130 131 133 137 148 156 157 159 164 172 173 175 183 184 186 190
193 198 203 208 210 211 213 217 228 236 237 239 244 252 253 255
273
 Unital No   9
 |Aut(G,M65)|=  128
  1   7  14  16  18  21  28  29  33  39  46  48  54  56  58  59
 65  71  78  80  86  88  90  91  97 103 110 112 115 116 121 127
134 136 138 139 147 148 153 159 166 168 170 171 182 184 186 187
195 196 201 207 211 212 217 223 225 231 238 240 243 244 249 255
273
 Unital No  10
 |Aut(G,M65)|=  12
  9  10  12  16  17  21  28  30  36  44  45  47  49  50  52  58
 65  67  72  77  99 100 110 112 116 124 125 127 129 131 136 141
147 154 155 159 165 166 170 173 177 181 188 190 197 200 207 208
211 212 222 224 231 232 234 238 245 248 255 256 263 264 266 270
273
 Unital No  11
 |Aut(G,M65)|=  12
  3  18  19  22  25  26  37  38  43  44  46  52  53  57  58  61
 66  72  74  77  80  82  83  88  92  95  99 101 104 105 107 116
120 123 125 127 132 134 138 140 143 146 150 151 155 157 164 185
195 197 198 205 207 209 213 216 218 220 239 242 244 249 251 252
264
pp16\pp-16-4.DA
 Unital No   1
 |Aut(G,M65)|=  1200.
  1  17  36  37  41  45  48  49  66  69  70  71  72  84  88  91
 94  95 100 101 105 109 112 114 117 118 119 120 135 138 139 140
141 145 163 166 169 170 174 180 184 187 190 191 195 198 201 202
206 215 218 219 220 221 226 227 236 239 240 242 243 252 255 256
257
 Unital No   2
 |Aut(G,M65)|=  48
  2  17  36  39  42  47  48  58  65  66  67  76  79  84  85  88
 90  92  99 100 101 104 106 113 114 115 124 128 131 135 136 137
144 148 165 166 171 172 175 180 182 186 191 192 195 197 199 205
207 214 216 220 222 224 225 226 229 230 231 241 242 246 247 248
258
 Unital No   3
 |Aut(G,M65)|=  8
  1   2   4  10  17  23  25  31  41  42  44  48  57  58  60  64
 66  67  69  73  99 100 110 112 113 114 116 122 133 136 143 144
148 156 157 159 161 165 172 174 177 181 188 190 195 202 203 207
211 212 222 224 226 227 229 233 242 246 254 255 258 263 269 272
273
 Unital No   4
```

```
 |Aut(G,M65)|=  16
  5  8 15 16 17 23 25 31 33 35 40 45 57 58 60 64
 66 67 69 73 97 103 105 111 115 118 119 124 133 136 143 144
148 156 157 159 161 165 172 174 178 183 189 192 194 199 205 208
210 216 219 220 226 227 229 233 241 243 248 253 258 263 269 272
273
 Unital No   5
 |Aut(G,M65)|=  32
 21 24 31 32 36 37 39 43 51 54 55 60 67 74 75 79
 81 86 91 96 101 102 106 109 116 117 119 123 131 132 142 144
147 150 151 156 164 172 173 175 178 183 189 192 201 202 204 208
211 218 219 223 229 230 234 237 244 252 253 255 265 267 269 270
273
 Unital No   6
 |Aut(G,M65)|=  100
  1  6  9 12 13 22 24 26 28 31 33 36 39 45 46 52
 54 59 60 62 80 85 87 93 94 95 98 101 103 106 111 114
120 121 122 125 139 146 149 150 156 160 169 183 193 196 200 203
207 209 212 213 217 224 232 242 250 251 254 256 257 262 265 268
269
pp16\pp-16-5
 Unital No   1
 |Aut(G,M65)|=  48
  9 11 13 14 21 24 31 32 38 41 42 48 51 54 55 60
 66 69 70 78 83 88 90 94 101 106 108 109 114 119 120 125
132 140 142 143 163 164 173 176 183 186 187 191 196 198 200 203
212 213 215 217 226 227 233 239 242 251 252 256 257 258 260 266
273
 Unital No   2
 |Aut(G,M65)|=  1200.
  9 11 13 14 19 22 23 28 37 38 42 45 66 70 78 79
 83 90 91 95 105 106 108 112 114 119 125 128 132 134 136 137
149 152 159 160 164 165 167 171 183 184 186 190 196 204 205 207
211 212 222 224 226 227 229 233 242 248 251 252 257 258 260 266
273
 Unital No   3
 |Aut(G,M65)|=  32
  5  8 15 16 17 18 20 26 41 43 45 46 51 54 55 60
 69 72 79 80 85 88 95 96 101 104 111 112 117 120 127 128
131 134 135 140 153 155 157 158 163 166 167 172 185 187 189 190
201 203 205 206 217 219 221 222 227 230 231 236 243 246 247 252
273
 Unital No   4
 |Aut(G,M65)|=  8
  7 12 14 15 16 19 27 28 29 32 33 35 37 39 48 49
 56 59 62 63 69 73 75 76 79 82 85 87 90 95 97 99
102 106 110 113 114 117 123 125 141 151 163 164 165 170 172 192
194 202 204 205 206 210 211 215 219 222 225 234 237 239 240 242
259
 Unital No   5
 |Aut(G,M65)|=  16
  4  6  7 13 14 21 27 28 29 32 33 34 40 41 48 54
 56 60 62 64 65 68 69 73 77 82 85 87 90 92 97 101
102 103 104 114 116 117 118 128 132 152 162 164 169 172 174 182
199 200 201 204 205 209 210 211 215 222 225 237 238 239 240 249
259
 Unital No   6
 |Aut(G,M65)|=  100
  3  7  8  9 16 17 22 24 26 31 33 37 39 42 48 49
 52 59 62 63 70 84 85 89 90 95 100 113 117 120 121 123
135 150 151 153 154 155 165 166 171 172 175 192 207 210 213 214
215 216 225 228 233 237 240 244 248 250 251 256 263 267 268 270
272
pp16\pp-16-6
 Unital No   1
 |Aut(G,M65)|=  32
  7  8 12 15 39 42 43 48 50 51 52 60 67 72 73 74
 85 88 91 93 98 99 103 106 114 118 120 125 135 136 140 143
150 155 156 157 161 162 169 171 182 185 186 191 195 202 203 207
210 214 222 223 227 233 236 237 246 247 249 253 260 266 269 270
273
 Unital No   2
 |Aut(G,M65)|=  16
  8 13 14 15 20 23 26 28 33 34 35 38 52 55 58 60
 72 77 78 79 85 89 91 96 97 98 99 102 120 125 126 127
129 130 131 134 145 146 147 150 168 173 174 175 180 183 186 188
```

196 199 202 204 212 215 218 220 232 237 238 239 241 242 243 246
273
pp16\pp-16-7
 Unital No  1
 |Aut(G,M65)|=  32
  1  3  6 12 15 18 19 22 23 24 38 40 41 43 44 49
 51 53 57 61 65 66 71 79 80 83 87 89 91 95 111 124
133 135 136 137 144 146 150 155 158 160 165 171 172 175 176 184
194 195 196 197 203 209 214 217 218 224 231 241 242 245 248 252
258
 Unital No  2
 |Aut(G,M65)|=  16
  8 13 14 15 20 23 26 28 33 34 35 38 52 55 58 60
 72 77 78 79 85 89 91 96 97 98 99 102 120 125 126 127
129 130 131 134 145 146 147 150 168 173 174 175 180 183 186 188
196 199 202 204 212 215 218 220 232 237 238 239 241 242 243 246
273
pp16\pp-16-8
 Unital No  1
 |Aut(G,M65)|=  8
  1 29 39 49 60 62 63 64 66 68 69 76 78 85 97 98
 99 100 111 113 114 116 119 123 132 133 134 137 144 146 150 156
157 160 167 170 171 172 175 182 184 187 189 190 193 197 198 207
208 214 215 220 221 223 228 229 235 237 238 242 247 251 254 256
269
 Unital No  2
 |Aut(G,M65)|=  16
 18 26 27 29 37 38 40 44 49 56 58 64 68 69 73 76
 83 87 95 96 99 100 104 107 113 119 121 127 132 133 137 140
145 146 147 150 163 170 171 175 180 182 189 192 193 199 201 207
216 218 220 221 226 229 230 235 242 247 253 256 257 260 265 270
273
 Unital No  3
 |Aut(G,M65)|=  4.
 14 19 23 26 31 32 35 37 38 40 44 49 58 60 61 62
 76 81 84 85 89 94 98 99 104 108 111 114 116 119 121 127
132 145 146 153 157 160 170 179 182 184 190 192 193 197 199 201
208 214 216 218 221 223 226 228 229 230 237 247 257 260 265 267
270
 Unital No  4
 |Aut(G,M65)|=  48
  2 10 11 13 17 20 25 30 35 39 47 48 49 59 60 63
 68 69 73 76 97 98 99 102 116 119 126 127 132 133 137 140
146 156 158 160 163 170 171 175 182 183 187 190 193 195 201 208
213 214 218 221 229 230 234 237 242 247 253 256 257 260 265 270
273
 Unital No  5
 |Aut(G,M65)|=  24
  5  6  8 12 19 23 31 32 35 39 47 48 56 59 62 64
 67 68 78 80 100 102 107 111 115 117 119 120 133 134 138 141
147 150 157 158 164 167 174 175 180 183 184 189 193 196 203 205
217 219 221 222 229 232 239 240 242 245 246 251 259 263 271 272
273
 Unital No  6
 |Aut(G,M65)|=  32
 17 25 26 29 33 39 46 48 52 53 58 64 66 68 75 78
 82 90 95 96 97 99 101 107 116 124 125 127 131 133 138 142
149 152 155 157 161 162 163 166 177 180 191 192 194 195 205 207
213 221 222 224 227 234 235 239 242 244 251 254 258 260 267 270
273
 Unital No  7
 |Aut(G,M65)|=  12
  2  4  9 13 18 20 25 29 42 44 46 47 49 55 62 64
 68 76 77 79 81 82 85 96 103 106 109 111 130 134 142 143
153 154 156 160 161 164 175 176 177 178 183 188 195 199 201 206
216 217 221 224 225 234 235 238 244 247 250 252 258 260 265 269
273
pp16\pp-16-9
 Unital No  1
 |Aut(G,M65)|=  8
  2  4  8 16 21 22 24 28 34 35 41 44 51 52 54 58
 69 71 77 78 91 92 94 95 102 105 107 111 114 115 118 126
132 133 135 139 147 151 159 160 184 186 190 192 194 202 203 205
212 220 221 223 229 233 234 240 247 248 249 253 257 258 259 262
273
 Unital No  2

|Aut(G,M65)|=  4.
  4  5  7  8 16 30 39 51 52 54 60 63 67 70 73 78
 79 87 88 89 94 96 112 117 119 121 126 128 132 133 135 136
144 149 150 152 153 156 164 165 166 168 172 179 185 188 190 191
193 194 202 203 205 211 212 214 220 223 239 243 258 261 265 267
272
 Unital No   3
 |Aut(G,M65)|=  24
  1  2  5  9 10 17 18 26 27 29 39 49 50 57 58 60
 67 68 71 75 78 84 93 94 95 96 112 115 116 119 122 126
129 136 137 139 141 147 151 158 159 160 163 164 167 175 176 178
180 190 191 192 194 201 202 203 205 209 214 217 219 221 239 243
262
 Unital No   4
 |Aut(G,M65)|=  16
  6 12 15 16 17 20 25 30 34 38 42 44 51 53 55 56
 67 71 75 77 91 93 95 96 102 107 108 109 114 115 119 122
131 134 135 140 146 154 155 157 165 166 168 172 178 186 191 192
195 199 207 208 213 216 223 224 226 229 232 234 245 248 251 253
273
 Unital No   5
 |Aut(G,M65)|=  48
 15 20 21 22 24 28 34 36 38 41 44 55 65 67 71 76
 78 81 85 94 95 96 100 102 105 107 108 113 115 119 120 126
131 147 151 158 159 160 165 166 168 169 172 177 182 190 191 192
193 195 199 207 208 224 228 229 232 233 234 244 245 248 249 253
259
 Unital No   6
 |Aut(G,M65)|=  32
  6  7  8  9 12 26 34 37 45 46 47 58 65 67 68 75
 80 81 83 84 91 96 97 99 100 107 112 114 117 125 126 127
134 135 136 137 140 150 151 152 153 156 170 178 181 189 190 191
198 199 200 201 204 218 225 227 228 235 240 242 245 253 254 255
257
 Unital No   7
 |Aut(G,M65)|=  12
  3  5  6  8 11 32 33 35 42 46 47 64 65 67 75 78
 79 85 87 90 92 94 97 102 104 106 112 118 119 120 123 126
135 145 147 150 156 160 166 168 170 172 175 179 181 185 187 189
199 213 215 216 222 223 226 228 234 236 239 241 245 251 252 256
264
pp16\pp-16-10
 Unital No   1
 |Aut(G,M65)|=  12
 10 18 34 37 38 39 42 49 69 71 76 77 79 82 83 90
 92 95 98 99 106 108 112 115 118 120 126 128 129 132 133 135
144 148 161 164 166 168 175 178 182 183 184 186 193 195 196 197
208 209 212 216 220 223 229 231 235 236 239 243 246 248 249 256
260
 Unital No   2
 |Aut(G,M65)|=  24
  5  8 15 16 25 27 29 30 34 44 46 48 67 71 75 77
 82 87 88 89 97 99 105 112 113 116 117 120 131 133 140 143
147 150 151 156 162 164 171 174 177 182 184 190 194 196 201 205
213 214 215 223 230 235 236 237 241 244 255 256 257 258 260 266
273
pp16\pp-16-11
 Unital No   1
 |Aut(G,M65)|=  24
  1  9 11 13 14 20 25 27 29 30 34 38 39 42 46 51
 54 55 58 60 67 82 83 89 90 92 98 99 106 107 108 124
129 132 137 141 144 146 147 150 151 156 161 164 171 174 175 178
182 183 186 189 193 196 197 201 205 209 212 216 219 222 230 247
260
 Unital No   2
 |Aut(G,M65)|=  12
  2  5 12 14 21 23 26 27 34 39 43 48 52 53 54 61
 67 68 77 79 81 86 93 96 106 108 110 112 122 123 124 127
129 130 132 138 145 147 149 153 163 166 167 172 179 180 185 192
217 219 221 222 225 230 233 239 242 247 254 255 257 258 260 266
273
pp16\pp-16-12
 Unital No   1
 |Aut(G,M65)|=  48
 20 21 22 30 33 39 46 48 50 54 57 64 69 72 79 80
 82 83 88 94 101 104 111 112 115 116 121 127 130 131 136 142

149 151 153 154 163 166 167 172 180 181 182 190 195 196 201 207
212 215 216 221 226 231 235 239 242 246 249 256 258 262 270 271
273
 Unital No   2
 |Aut(G,M65)|=  16
  1  7 14 16 18 28 30 32 35 36 41 47 50 56 59 60
 70 74 78 80 85 88 91 93 101 102 106 109 114 117 124 125
133 141 142 144 145 151 152 155 161 165 167 173 178 182 186 188
198 200 202 203 209 210 215 220 225 230 231 234 248 251 254 256
273
 Unital No   3
 |Aut(G,M65)|=  12
 18 28 30 32 37 41 43 48 51 58 61 64 67 69 76 79
 84 89 90 91 98 102 105 112 113 120 121 124 129 131 133 139
149 154 155 156 161 162 173 174 183 188 189 190 194 201 202 206
211 218 221 224 225 227 229 235 241 242 253 254 261 262 266 269
273
pp16\pp-16-13
 Unital No   1
 |Aut(G,M65)|=  16
  2  3 13 15 20 21 26 32 38 42 46 48 50 55 59 63
 70 74 78 80 83 87 91 93 98 100 101 111 114 118 126 127
132 133 138 144 150 151 155 158 163 167 171 173 180 181 182 190
193 200 201 204 212 213 215 219 227 234 237 240 242 243 253 255
273
 Unital No   2
 |Aut(G,M65)|=  48
  3  4  5 10 12 19 23 25 27 32 35 39 41 43 48 49
 52 59 62 63 66 73 74 77 79 83 84 85 90 92 104 120
130 137 138 141 143 152 166 172 173 174 176 177 180 187 190 191
193 194 197 198 199 216 225 226 229 230 231 246 252 253 254 256
257
 Unital No   3
 |Aut(G,M65)|=  24
  9 19 24 27 29 31 35 40 43 45 47 49 53 58 61 62
 66 71 72 74 80 89 98 102 107 108 110 116 117 119 124 127
130 135 136 138 144 146 150 155 156 158 169 177 181 186 189 190
193 195 196 198 208 212 213 215 220 223 225 227 228 230 240 249
257
pp16\pp-16-14
 Unital No   1
 |Aut(G,M65)|=  16
  2  4 10 15 18 19 28 29 33 42 43 47 51 55 60 62
 71 74 76 77 82 88 92 95 97 99 102 110 113 115 123 125
129 135 136 143 145 150 152 158 168 171 172 173 182 185 186 190
194 198 199 203 211 213 216 218 231 234 235 239 242 246 247 251
273
 Unital No   2
 |Aut(G,M65)|=  32
  2  4 10 15 20 22 26 32 35 39 41 46 51 55 60 62
 67 68 75 78 82 88 92 95 98 105 106 112 113 115 123 125
129 135 136 143 145 150 152 158 168 171 172 173 182 185 186 190
205 206 207 208 210 214 215 219 225 226 227 228 241 244 249 252
273
 Unital No   3
 |Aut(G,M65)|=  32
  7  8  9 14 18 19 28 29 33 42 43 47 51 55 60 62
 71 74 76 77 81 87 90 96 97 99 102 110 113 115 123 125
129 135 136 143 154 155 156 160 168 171 172 173 180 184 187 192
195 197 200 202 221 222 223 224 229 233 236 238 253 254 255 256
273
 Unital No   4
 |Aut(G,M65)|=  48
  3  8  9 12 25 28 29 31 33 36 47 48 53 58 62 64
 73 76 77 79 81 84 85 90 97 100 110 111 117 122 125 126
131 136 141 144 149 153 154 156 163 168 175 176 177 179 180 184
205 206 207 208 210 214 215 219 241 244 249 252 259 261 264 266
273
 Unital No   5
 |Aut(G,M65)|=  16
  3  8  9 12 18 19 24 27 41 44 45 46 50 57 59 60
 65 68 78 80 83 88 89 92 99 102 103 104 115 120 127 128
134 135 137 140 150 151 157 158 165 170 173 174 178 187 189 192
194 198 199 203 210 214 215 219 253 254 255 256 269 270 271 272
273
 Unital No   6

```
 |Aut(G,M65)|=  16
  2  9 10 14 17 21 24 25 36 37 40 44 49 50 52 59
 66 70 79 80 87 90 92 95 103 107 109 111 117 122 127 128
133 138 141 144 150 152 156 160 161 164 166 167 184 185 187 190
205 206 207 208 210 214 215 219 225 226 227 228 241 244 249 252
273
 Unital No   7
 |Aut(G,M65)|=  16
  3  5  6 11 17 21 27 31 34 40 44 48 51 56 62 64
 67 68 71 77 82 87 88 90 102 105 106 110 115 120 125 126
133 138 141 144 146 147 149 151 165 170 175 176 182 184 186 187
194 198 199 203 210 214 215 219 230 232 237 240 253 254 255 256
273
 Unital No   8
 |Aut(G,M65)|=  24
  2  6 13 15 18 21 22 24 37 39 40 43 51 58 61 62
 66 67 70 74 87 91 93 96 99 103 106 107 117 120 126 127
133 136 141 143 146 150 158 159 163 170 174 176 183 187 191 192
205 206 207 208 209 212 217 220 243 245 248 250 258 262 263 267
273
 Unital No   9
 |Aut(G,M65)|=   4
  2  4 12 13 16 19 26 27 29 31 39 40 41 43 44 56
 57 58 59 60 65 67 72 77 78 81 83 85 87 91 97 102
103 105 111 114 121 123 126 128 130 135 138 142 143 146 147 151
154 160 161 163 168 172 176 177 186 188 189 190 206 217 238 242
264
 Unital No  10
 |Aut(G,M65)|=   4
  1  2  7 14 16 17 19 25 29 31 39 41 44 45 48 50
 55 57 60 63 66 67 70 72 78 81 85 86 88 92 97 101
102 103 104 113 114 117 123 124 131 133 137 142 143 146 149 150
158 160 163 167 168 174 175 180 184 186 188 191 207 217 229 246
259
 Unital No  11
 |Aut(G,M65)|=   4
  6 11 12 14 16 17 21 25 29 31 34 41 44 45 48 54
 57 59 60 62 66 67 70 72 79 84 85 88 89 91 101 102
103 104 108 114 115 116 121 123 132 136 138 142 143 146 147 150
157 159 162 163 168 173 176 179 180 181 188 189 206 212 229 251
261
 Unital No  12
 |Aut(G,M65)|=   8
  1  3  4  8 19 22 28 32 33 39 42 46 50 53 54 56
 74 75 76 78 86 87 93 95 97 98 99 112 116 124 126 127
133 135 136 139 145 148 149 154 164 172 173 175 178 187 190 191
195 197 200 202 225 226 227 228 242 246 247 251 269 270 271 272
273
 Unital No  13
 |Aut(G,M65)|=   8
  1  2  3 14 17 28 31 32 35 37 38 43 50 54 61 62
 65 76 78 79 90 91 92 94 102 104 106 107 113 117 120 121
130 134 141 143 146 152 156 157 163 164 170 172 177 181 187 190
195 197 200 202 230 232 237 240 242 246 247 251 259 261 264 266
273
 Unital No  14
 |Aut(G,M65)|=   8
  1  4  5 10 18 20 23 25 33 34 39 44 51 54 57 62
 72 74 78 79 83 88 89 92 104 106 109 110 114 115 116 125
130 133 140 141 150 151 157 158 161 165 166 175 182 183 190 191
193 196 201 204 221 222 223 224 225 226 227 228 259 261 264 266
273
 Unital No  15
 |Aut(G,M65)|=   8
  3  4  5 12 16 17 22 25 26 27 33 37 40 45 47 54
 71 72 74 79 80 84 86 91 93 96 100 103 105 107 112 120
143 149 151 155 157 159 172 177 179 183 186 188 201 205 206 207
208 211 213 214 216 218 225 226 227 228 239 255 257 259 260 265
268
 Unital No  16
 |Aut(G,M65)|=   8
  2  3  4  8 12 17 37 50 51 54 56 63 79 82 90 91
 93 96 107 114 117 120 123 124 129 132 139 141 143 151 155 156
157 158 164 165 172 173 174 177 179 180 186 190 193 195 197 200
202 210 211 213 216 218 230 232 237 239 240 253 257 260 265 266
268
```

Unital No 17
|Aut(G,M65)|= 4
  1  3  4 10 16 21 22 25 26 27 33 36 40 45 47 49
 53 54 57 59 70 71 72 74 80 84 87 91 94 96 100 103
105 109 112 118 120 122 126 127 133 136 137 143 144 146 149 155
157 160 161 167 171 172 175 177 181 183 185 186 196 214 235 253
264
 Unital No 18
|Aut(G,M65)|= 8
  2  3  7  8 12 17 18 19 23 25 37 40 41 46 48 53
 56 57 60 63 67 69 75 79 80 82 88 90 91 93 97 103
107 108 110 114 118 124 125 126 129 131 136 137 139 145 151 156
157 158 162 165 171 173 175 177 179 180 190 192 201 215 234 253
261
 Unital No 19
|Aut(G,M65)|= 8
  3  4 10 12 19 23 26 27 33 41 45 48 49 50 54 57
 69 71 72 75 87 91 93 96 100 108 110 112 117 120 126 127
133 136 141 143 151 155 157 160 164 167 171 172 177 179 185 186
193 196 201 204 211 213 216 218 243 245 248 250 257 260 265 268
273
 Unital No 20
|Aut(G,M65)|= 8
  4  5  6 10 11 17 19 20 21 22 36 37 41 42 47 51
 70 74 75 77 79 83 85 86 87 96 98 105 107 109 110 123
132 148 153 155 159 160 174 185 186 188 189 191 195 197 200 201
202 214 221 222 223 224 230 232 237 239 240 255 259 269 270 271
272
 Unital No 21
|Aut(G,M65)|= 8
  3  4  9 13 15 20 21 26 30 31 33 38 39 41 45 50
 51 53 57 60 67 70 74 77 79 82 85 87 88 94 97 99
100 103 107 116 118 119 126 128 130 135 138 142 143 145 146 147
149 150 164 165 169 170 174 177 182 185 189 191 193 221 229 250
258
 Unital No 22
|Aut(G,M65)|= 8
  2  9 11 14 15 20 22 28 29 30 36 40 41 45 48 50
 51 54 57 60 66 70 72 73 74 87 88 89 90 92 98 99
103 106 110 116 118 119 124 128 131 135 138 142 143 146 147 151
159 160 164 165 170 174 176 177 179 184 188 191 196 221 233 248
262
 Unital No 23
|Aut(G,M65)|= 8
  5  8 13 15 18 22 30 31 33 37 40 41 50 51 54 58
 68 69 72 76 85 88 93 96 98 102 109 111 113 121 126 127
130 133 134 136 145 151 153 155 164 172 174 176 178 180 182 188
205 206 207 208 221 222 223 224 241 244 249 252 257 260 265 268
273
 Unital No 24
|Aut(G,M65)|= 8
  5  7  8 11 17 23 25 27 35 42 46 47 49 51 57 58
 65 66 70 73 81 89 93 96 99 106 110 112 119 123 126 127
130 134 141 143 145 153 158 159 163 164 170 172 178 179 182 186
194 198 199 203 221 222 223 224 242 246 247 251 269 270 271 272
273
 Unital No 25
|Aut(G,M65)|= 8
  2  7 14 15 17 28 29 30 36 41 46 48 49 54 59 60
 66 71 72 74 81 88 90 92 98 103 104 106 113 114 119 124
131 133 143 144 146 151 157 159 163 165 173 176 177 179 181 188
195 197 200 202 221 222 223 224 253 254 255 256 259 261 264 266
273
 Unital No 26
|Aut(G,M65)|= 8
  3  4  5  9 21 24 30 31 33 39 41 43 51 59 60 61
 67 74 77 78 82 87 93 94 100 103 107 108 113 115 119 128
133 135 137 143 145 147 149 156 164 165 171 173 182 187 189 191
193 196 201 204 231 234 235 239 241 244 249 252 269 270 271 272
273
 Unital No 27
|Aut(G,M65)|= 8
  3  5  7  9 12 17 19 28 30 31 33 42 44 45 48 53
 57 58 63 64 68 69 73 77 78 81 84 86 88 90 100 104
105 110 112 113 117 122 125 128 131 136 140 142 144 147 149 153
155 156 163 164 168 173 175 177 178 180 184 186 208 223 228 254

269
pp16\pp-16-15
 Unital No  1
 |Aut(G,M65)|=  16
  2  3  8 10 15 17 18 20 22 29 34 35 37 38 46 49
 50 52 55 62 64 72 79 80 90 91 94 98 104 109 113 125
127 130 132 135 138 141 143 151 154 163 174 179 182 184 185 187
190 195 196 199 200 202 205 219 224 225 230 235 236 239 240 244
247
 Unital No  2
 |Aut(G,M65)|=  32
  5  9 10 12 14 17 20 22 28 29 35 43 44 46 47 52
 54 55 57 59 77 85 97 121 130 132 133 134 135 138 141 147
148 151 169 173 174 179 184 185 187 188 190 192 195 196 199 202
205 206 208 209 211 224 225 229 230 234 235 236 240 245 247 250
259
 Unital No  3
 |Aut(G,M65)|=  48
 18 21 28 31 34 37 44 47 49 54 59 60 64 66 69 76
 79 84 87 90 93 98 101 108 111 116 119 122 125 132 135 138
141 145 150 155 160 162 165 172 175 179 184 185 190 196 199 202
205 212 215 218 221 225 230 235 240 241 246 251 256 257 261 272
273
 Unital No  4
 |Aut(G,M65)|=  32
 10 28 43 57 66 69 72 73 77 84 85 86 87 91 97 99
101 104 108 113 118 119 121 122 131 132 134 135 136 138 141 147
155 160 172 173 175 179 180 183 184 185 188 190 196 197 199 202
204 205 208 211 212 221 225 230 232 233 234 235 240 241 245 256
259
 Unital No  5
 |Aut(G,M65)|=  16
 15 29 46 63 64 65 67 69 73 79 84 86 90 94 96 101
103 104 110 111 113 116 122 123 124 131 132 135 138 141 142 143
145 151 160 161 163 165 166 171 172 176 183 186 187 196 197 198
199 202 204 205 210 211 212 213 220 221 223 227 237 238 246 251
252
 Unital No  6
 |Aut(G,M65)|=  16
  3  8  9 14 17 20 22 23 26 27 29 32 35 40 41 46
 63 84 87 90 93 130 132 133 135 138 140 141 143 161 166 171
176 177 178 181 182 187 188 191 192 195 196 199 200 201 202 205
206 210 213 220 223 226 227 229 232 233 236 238 239 258 266 270
273
 Unital No  7
 |Aut(G,M65)|=  24
 50 81 86 91 96 97 100 102 103 106 107 109 112 114 116 117
119 122 124 125 127 129 130 133 134 139 140 143 144 146 149 156
159 161 164 166 167 170 171 173 176 194 197 204 207 211 212 215
216 217 218 221 222 226 228 229 231 234 236 237 239 265 269 270
272
 Unital No  8
 |Aut(G,M65)|=  8
  5 26 35 53 64 66 68 85 90 107 110 118 126 131 134 142
144 145 147 148 151 153 154 157 160 165 167 172 173 177 178 182
183 186 187 188 192 193 196 197 203 204 210 212 216 221 222 225
226 227 229 231 236 237 238 239 242 243 244 246 247 250 251 252
253
 Unital No  9
 |Aut(G,M65)|=  4
  2  9 13 29 31 32 34 40 48 50 53 55 59 60 62 63
 67 75 76 83 84 91 103 104 111 119 127 128 129 131 132 146
148 151 152 154 157 160 161 165 170 179 186 191 196 203 207 210
216 218 225 232 237 241 242 243 244 247 250 253 257 259 266 267
268
 Unital No  10
 |Aut(G,M65)|=  4
  2  9 13 29 31 32 34 40 48 50 53 55 59 60 62 63
 66 67 75 83 90 91 103 110 111 118 119 127 132 134 142 145
146 147 161 172 173 177 179 182 183 187 188 192 193 204 205 221
222 223 226 227 229 231 236 239 240 246 252 254 259 263 264 271
273
 Unital No  11
 |Aut(G,M65)|=  4
  2  7  9 21 29 32 34 38 40 50 53 55 56 59 60 63
 65 67 69 89 91 93 105 109 111 113 117 119 131 132 144 146

```
 153 160 161 165 167 177 178 179 182 186 187 192 195 199 200 201
 203 204 206 213 216 221 227 230 237 242 246 248 257 259 266 267
 268
 Unital No  12
 |Aut(G,M65)|=   8
  5   9  12  13  14  17  18  22  23  26  34  35  37  43  46  49
 56  58  61  63  64  65  66  70  75  79  80  83  87  88  89  90
 94  99 110 118 123 132 135 136 138 141 150 165 179 183 184 185
 190 194 197 201 205 213 218 221 224 225 226 227 232 247 249 251
 256
 Unital No  13
 |Aut(G,M65)|=  16
  8  12  14  17  26  27  37  41  47  52  58  59  63  68  70  74
 85  89  95 100 102 112 115 121 127 129 134 136 137 139 140 144
 148 150 155 165 168 172 178 181 183 186 188 191 192 193 196 205
 213 216 220 225 234 240 248 249 255 257 260 262 263 265 267 271
 272
 Unital No  14
 |Aut(G,M65)|=   8
  4   9  15  20  22  28  41  47  48  50  52  53  54  56  60  63
 65  72  74  83  91  95  97  98 100 114 120 125 135 136 143 146
 154 160 168 172 176 177 178 180 181 185 188 191 198 201 202 219
 220 222 234 235 236 241 242 246 249 250 251 256 259 262 263 267
 272
 Unital No  15
 |Aut(G,M65)|=   8
  3  10  15  18  20  32  37  38  41  50  53  54  58  60  62  63
 65  67  68  83  85  96 100 101 107 114 122 126 131 132 143 149
 154 155 168 171 175 177 178 181 183 188 190 191 195 198 205 209
 222 223 226 234 240 241 242 244 246 251 254 256 259 262 263 267
 272
 Unital No  16
 |Aut(G,M65)|=   8
  3   5   8  15  20  26  27  32  34  35  37  41  54  58  60  61
 64  68  74  79  84  85  95  97 105 107 118 120 126 129 132 134
 139 140 143 144 147 149 152 153 154 157 158 163 168 171 177 182
 190 193 206 212 215 217 218 221 224 227 232 233 237 238 239 245
 253
 Unital No  17
 |Aut(G,M65)|=   8
  7  31  43  50  51  74  76  77  78  80  81  92  93  94  95  97
 99 102 108 109 115 118 119 120 127 129 134 135 138 139 144 146
 147 152 153 158 159 161 176 184 185 200 201 202 205 212 214 215
 218 219 220 221 223 227 229 232 233 235 236 238 240 247 249 250
 254
 Unital No  18
 |Aut(G,M65)|=   8
  3   5   7   8  12  23  26  27  31  32  34  35  37  43  46  49
 50  51  58  61  64  76  78  80  81  93  94  99 108 109 118 119
 127 132 138 140 143 149 154 157 159 161 163 168 171 177 182 184
 190 199 200 202 206 213 214 220 224 229 230 235 239 247 248 249
 253
 Unital No  19
 |Aut(G,M65)|=   8
  5   8  14  15  17  20  26  27  34  35  41  44  54  55  61  63
 64  66  69  76  79 130 132 138 140 147 149 151 152 153 157 158
 159 161 163 169 171 182 184 190 192 196 200 202 206 210 212 214
 215 218 220 221 224 225 229 235 239 243 247 249 253 258 266 270
 273
 Unital No  20
 |Aut(G,M65)|=   4
  4   8  12  23  27  28  34  46  48  49  50  53  56  60  61  63
 68  73  75  85  86  89 106 107 111 116 124 126 130 137 138 145
 151 159 161 164 165 167 169 170 173 184 189 192 194 197 198 199
 204 206 207 213 222 224 230 234 239 242 248 253 258 259 266 270
 273
 Unital No  21
 |Aut(G,M65)|=   4
  4   8  15  20  27  28  34  41  48  50  53  54  56  60  61  63
 65  68  72  83  85  91  98 100 107 114 125 126 129 134 138 139
 142 143 144 150 154 159 161 162 168 177 184 186 198 205 206 212
 215 218 221 222 223 224 234 239 240 242 253 254 257 259 260 265
 271
 Unital No  22
 |Aut(G,M65)|=   8
 10  12  15  20  23  28  35  40  43  57  59  64  68  72  73  81
```

```
 84  87  90  92  93  96  98 111 112 113 118 119 120 122 123 128
129 132 133 135 140 151 152 154 156 159 161 166 168 169 170 182
185 187 190 191 198 217 228 255 257 259 260 261 264 265 268 269
271
 Unital No  23
 |Aut(G,M65)|=   8
  3   8   9  14  33  35  38  40  41  43  46  48  50  53  60  63
 81  82  85  86  91  92  95  96  99 100 103 104 105 106 109 110
132 135 138 141 147 152 153 158 163 168 169 174 180 183 186 189
225 228 230 231 234 235 237 240 241 244 246 247 250 251 253 256
259
 Unital No  24
 |Aut(G,M65)|=   8
  2   3   4   6   8  11  13  17  19  20  21  22  25  28  33  34
 35  37  48  49  52  55  56  57  65  70  71  73  74  75  80  85
 86  92 102 107 111 116 120 121 140 143 148 151 163 168 187 192
194 195 210 213 215 219 220 223 233 239 243 247 248 249 254 256
259
 Unital No  25
 |Aut(G,M65)|=   8
  6  11  15  19  25  29  46  64  66  69  78  81  84  87 105 108
110 114 117 119 123 124 127 128 134 137 139 142 145 150 152 153
164 172 173 175 178 180 183 191 194 196 198 199 201 202 203 205
209 211 215 222 225 228 229 230 233 235 237 240 245 252 253 256
259
 Unital No  26
 |Aut(G,M65)|=   8
  4   8  15  22  28  29  37  46  48  50  53  55  56  60  63  64
 68  69  78  81  82  84  97 108 110 115 119 123 137 139 143 146
149 150 151 152 156 159 163 173 175 178 180 187 194 198 201 209
211 215 229 233 237 243 248 249 252 253 254 256 257 258 259 263
269
 Unital No  27
 |Aut(G,M65)|=   8
  4   8   9  16  22  27  28  30  36  37  42  44  48  55  56  58
 68  73  78  81  82  86  97 108 111 114 115 116 117 119 124 127
130 143 144 147 151 154 163 170 174 181 182 187 193 194 197 212
215 216 233 234 238 251 255 256 257 258 259 260 263 264 266 269
272
 Unital No  28
 |Aut(G,M65)|=   4
  3   4  15  17  28  29  34  46  48  50  52  53  56  60  63  64
 68  73  79  82  86  90  97 104 111 113 115 116 130 137 144 147
150 154 170 174 175 179 180 181 182 184 185 190 193 194 206 214
215 216 226 233 234 243 248 249 250 254 255 256 259 260 264 266
272
 Unital No  29
 |Aut(G,M65)|=   8
  5  15  23  29  40  46  54  60  64  68  75  82  83  97 103 115
127 137 139 144 147 150 152 170 173 175 178 180 181 193 194 196
199 201 202 203 205 206 209 214 215 216 222 226 228 229 233 234
243 245 248 249 250 253 254 255 256 257 258 261 263 267 269 270
271
pp16\pp-16-16
 Unital No   1
 |Aut(G,M65)|=   8
  7  10  12  14  18  24  28  32  36  40  44  46  52  53  54  64
 67  70  73  77  83  89  90  93 100 102 106 109 113 115 119 126
131 135 140 141 151 152 156 160 163 166 169 170 180 181 183 192
195 197 200 202 217 221 222 224 237 238 239 240 259 261 262 263
273
 Unital No   2
 |Aut(G,M65)|=   8
  5  10  14  15  17  19  20  24  35  36  41  44  54  55  63  64
 70  71  77  78  81  83  86  93  97  98 100 107 113 116 125 126
133 134 135 138 147 152 158 160 168 169 172 175 177 180 190 192
195 197 200 202 217 221 222 224 250 252 253 255 259 261 262 263
273
 Unital No   3
 |Aut(G,M65)|=   8
  5   6  12  14  17  18  23  28  36  40  41  43  49  53  55  64
 67  70  71  72  86  87  90  93 102 104 106 107 113 116 118 119
129 131 133 140 147 152 153 156 166 168 170 171 177 183 186 192
195 197 200 202 218 219 220 223 227 229 231 235 250 252 253 255
273
 Unital No   4
```

|Aut(G,M65)|=  8
  6  7 14 15 18 22 30 32 35 36 43 46 49 52 63 64
 67 71 76 77 83 86 90 92 99 106 109 111 116 119 120 126
132 140 141 143 147 150 156 160 163 164 170 174 181 186 190 192
195 197 200 202 218 219 220 223 227 229 231 235 265 267 268 269
273
 Unital No   5
|Aut(G,M65)|=  4
  3  5  7 13 14 18 21 23 31 32 34 36 41 46 48 50
 52 55 62 64 65 69 71 73 74 85 86 89 91 94 104 106
108 109 112 114 115 116 122 123 129 139 140 141 144 145 147 148
151 155 161 163 165 170 171 177 181 188 189 192 204 211 237 255
265
 Unital No   6
|Aut(G,M65)|=  4
  2  3  7 14 15 17 18 20 30 32 35 36 39 46 48 50
 52 57 63 64 66 67 71 74 79 84 86 88 90 94 98 106
107 109 111 116 119 121 123 128 133 138 140 141 143 147 148 149
156 157 163 164 168 170 175 179 181 188 190 192 201 211 240 253
269
 Unital No   7
|Aut(G,M65)|=  8
  1  3 11 16 18 21 25 29 36 40 41 43 49 53 55 64
 65 66 69 79 81 87 89 93 103 106 108 110 114 121 122 128
130 131 137 143 151 152 153 158 164 166 167 176 178 184 188 191
209 212 215 216 227 229 231 235 250 252 253 255 259 261 262 263
273
 Unital No   8
|Aut(G,M65)|=  24
  4  5  6 13 18 20 21 32 33 38 44 46 53 57 59 64
 67 70 76 78 84 88 94 96 97 103 104 110 114 122 124 127
131 132 135 143 147 151 153 160 161 167 168 176 182 184 188 190
210 211 213 214 237 238 239 240 241 243 244 247 257 258 260 264
273
 Unital No   9
|Aut(G,M65)|=  12
  2  6 10 11 16 17 19 21 23 29 34 39 41 42 46 52
 53 54 56 61 65 69 70 74 80 81 85 89 95 96 98 102
103 104 110 117 120 121 122 126 131 132 133 136 143 145 146 147
148 157 164 169 174 175 176 177 178 182 186 189 208 212 237 246
263
 Unital No  10
|Aut(G,M65)|=  4
  1  3  9 14 15 17 18 23 25 26 34 39 41 43 44 49
 50 51 52 54 68 70 71 72 73 83 86 88 89 92 100 105
106 109 111 114 119 120 125 126 135 138 139 143 144 152 153 156
158 159 161 167 171 173 175 179 182 183 189 192 197 213 227 241
261
 Unital No  11
|Aut(G,M65)|=  4
  2  6 12 14 15 19 21 23 25 27 38 41 44 46 47 49
 52 55 62 64 70 71 72 74 79 83 84 89 92 94 99 100
101 108 111 120 122 124 125 126 131 133 135 141 143 152 153 154
155 156 162 163 168 171 173 180 182 183 188 190 197 211 235 247
263
 Unital No  12
|Aut(G,M65)|=  4
  1  2  6 14 15 20 23 24 26 27 38 39 41 44 47 49
 50 52 62 64 67 69 70 71 77 81 82 87 89 92 100 104
106 108 112 113 115 117 125 126 130 136 138 141 143 147 149 150
153 156 161 164 168 170 173 180 182 183 188 189 195 213 235 244
262
 Unital No  13
|Aut(G,M65)|=  12
  2  9 12 13 14 19 21 27 29 31 33 35 40 42 47 49
 51 56 62 64 65 67 70 71 76 88 89 94 95 96 97 99
108 109 111 114 120 121 122 124 129 130 135 142 144 146 147 151
152 153 163 164 169 170 176 178 181 187 190 191 208 215 239 248
261
pp16\pp-16-17
 Unital No   1
|Aut(G,M65)|=  8
  2  8 10 14 20 21 22 32 34 35 44 48 49 51 52 58
 63 66 69 72 76 79 82 85 92 95 96 97 102 107 108 112
116 133 138 141 143 148 154 156 159 161 163 166 169 179 182 184
192 196 199 210 213 225 230 243 248 257 261 262 266 269 271 272

273
 Unital No   2
 |Aut(G,M65)|=  8
  5  7  8 14 22 23 28 32 33 34 40 44 52 54 58 62
 63 66 69 76 78 79 86 98 122 129 131 133 136 143 144 147
 148 154 155 158 160 163 164 169 172 173 175 178 180 182 183 191
 192 193 196 199 208 210 211 213 222 225 230 231 234 243 245 248
 252
 Unital No   3
 |Aut(G,M65)|=  8
  5 10 21 23 40 48 51 54 63 67 82 85 91 92 95 97
 102 107 111 112 113 118 119 123 128 131 134 136 138 139 141 152
 153 155 156 159 160 161 166 167 170 172 175 179 180 181 183 184
 188 193 208 211 222 231 234 245 252 257 261 262 266 269 271 272
 273
 Unital No   4
 |Aut(G,M65)|=  8
  7 10 21 28 33 48 51 60 62 66 69 76 78 79 82 85
 91 92 95 99 101 104 105 110 115 116 120 121 126 129 130 138
 139 143 147 148 153 157 159 161 164 168 169 170 178 182 184 187
 188 196 197 202 210 212 220 225 235 238 243 249 251 257 261 272
 273
 Unital No   5
 |Aut(G,M65)|=  4
  2  8 10 20 21 32 35 44 48 49 50 51 52 53 60 63
 66 69 76 79 82 85 92 95 100 103 106 109 129 133 135 136
 146 147 154 160 163 164 172 176 178 180 185 192 193 195 199 207
 213 218 219 222 230 231 232 239 241 247 248 252 257 259 261 272
 273
 Unital No   6
 |Aut(G,M65)|=  4
  2 10 14 20 21 22 34 35 48 49 50 51 53 58 60 63
 66 69 76 79 83 88 89 94 115 120 121 126 129 138 142 143
 147 148 150 159 161 162 164 169 178 182 184 186 197 198 199 201
 209 212 213 217 228 229 230 238 248 251 253 255 257 259 261 272
 273
 Unital No   7
 |Aut(G,M65)|=  24
  2  4  5  8 13 14 17 20 21 23 27 28 37 38 41 43
 46 47 50 51 52 58 59 62 64 65 70 79 83 88 90 104
 106 109 113 124 127 129 131 139 150 152 158 164 170 172 181 189
 191 198 203 205 211 220 222 228 230 237 243 245 252 265 269 270
 272
 Unital No   8
 |Aut(G,M65)|=  12
  2  5  9 10 14 31 33 34 37 41 46 50 53 56 60 63
 70 76 77 79 81 83 86 92 99 102 104 109 115 122 124 125
 130 134 135 140 145 155 158 159 163 169 173 176 180 181 184 186
 193 194 200 204 205 206 213 217 227 229 230 233 234 239 242 254
 259
 Unital No   9
 |Aut(G,M65)|=  8
 14 22 34 50 58 78 81 82 85 92 95 108 119 131 134 136
 137 138 141 142 146 149 152 161 162 165 166 167 172 175 181 185
 190 193 195 196 197 198 199 201 202 205 207 211 214 215 216 221
 224 226 227 233 234 236 237 241 243 247 248 249 251 252 253 254
 255
 Unital No  10
 |Aut(G,M65)|=  4
  8 10 21 22 37 48 50 51 53 55 60 63 66 68 69 72
 75 76 79 82 89 96 97 106 108 115 116 124 130 134 141 147
 149 151 161 169 173 180 183 185 186 188 189 192 203 205 207 209
 210 211 212 214 219 224 231 233 235 246 248 255 257 258 259 263
 269
 Unital No  11
 |Aut(G,M65)|=  4
  8 10 21 22 37 48 50 51 53 55 60 63 71 73 75 81
 85 89 99 101 102 104 105 106 110 120 124 125 134 140 141 147
 149 157 161 163 173 180 182 183 185 186 188 189 195 198 199 200
 201 204 206 215 220 222 225 234 238 241 242 254 259 260 264 266
 272
 Unital No  12
 |Aut(G,M65)|=  4
  3 10 17 21 34 48 50 51 52 53 60 63 73 75 77 81
 89 95 99 101 104 105 106 110 112 124 125 126 132 139 140 156
 157 158 162 163 164 165 172 175 176 181 182 184 198 205 207 209

210 212 214 219 222 224 233 234 235 242 246 248 257 258 259 263
269
 Unital No  13
 |Aut(G,M65)|=  12
  5  8 10 17 20 31 33 46 47 50 53 56 58 59 60 63
 70 77 78 79 83 90 91 92 101 102 104 109 113 115 116 124
131 139 141 143 150 154 156 158 164 166 168 172 177 179 181 189
194 203 205 206 215 219 220 222 228 229 230 233 243 244 245 256
259
pp16\pp-16-18
 Unital No   1
 |Aut(G,M65)|=  16
  9 19 37 49 66 89 91 92 93 96 97 101 106 108 112 113
116 118 119 124 131 132 133 134 142 146 148 151 153 158 164 165
166 167 171 177 178 184 188 189 194 195 200 201 202 217 218 219
221 224 225 228 235 237 240 242 244 248 249 254 259 262 263 265
270
 Unital No   2
 |Aut(G,M65)|=  16
  9 24 39 56 71 85 88 91 94 96 102 106 109 110 112 117
118 120 121 124 129 135 140 141 143 147 151 156 158 160 163 167
170 171 172 179 184 186 187 188 197 199 200 202 205 211 212 213
220 223 226 227 233 236 239 243 245 246 252 255 260 265 268 270
272
 Unital No   3
 |Aut(G,M65)|=  16
  9 31 47 63 79 81 82 84 91 93 99 101 103 104 106 113
116 118 124 128 129 132 133 138 139 145 148 149 150 151 165 166
167 171 176 178 184 188 189 192 195 200 201 202 208 211 214 216
217 219 225 226 229 232 237 242 243 249 250 253 258 259 265 266
269
 Unital No   4
 |Aut(G,M65)|=   8
  9 23 34 54 65 83 87 90 94 96 100 105 107 110 112 115
117 118 120 124 129 130 133 136 143 152 155 157 158 160 163 167
170 171 173 179 180 184 186 188 196 199 200 202 205 209 210 211
217 223 225 226 230 234 239 241 242 243 247 255 261 264 267 270
272
 Unital No   5
 |Aut(G,M65)|=   8
 17 18 21 22 33 36 40 41 58 59 60 61 74 75 76 77
 88 90 94 96 105 109 110 112 114 117 118 120 129 136 141 144
149 151 153 154 162 164 173 174 178 179 181 185 193 194 197 206
216 217 218 219 225 228 232 237 248 249 250 252 260 264 270 271
273
 Unital No   6
 |Aut(G,M65)|=  16
 19 30 31 32 42 43 44 45 58 59 60 61 66 78 79 80
 85 87 88 90 97 99 103 108 116 117 118 125 129 142 143 144
150 158 159 160 163 165 169 171 177 190 191 192 193 196 197 203
214 216 218 221 227 230 234 236 241 245 246 247 260 264 267 268
273
 Unital No   7
 |Aut(G,M65)|=  16
 26 27 28 29 34 35 38 39 52 54 56 57 74 75 76 77
 81 94 95 96 98 100 101 106 113 114 120 124 129 142 143 144
150 158 159 160 161 166 167 170 180 182 183 187 194 206 207 208
213 215 217 219 226 238 239 240 241 245 246 247 260 264 267 268
273
 Unital No   8
 |Aut(G,M65)|=  16
 20 23 24 25 33 36 40 41 50 51 53 55 65 71 72 73
 85 87 88 90 97 99 103 108 116 117 118 125 132 133 137 140
148 151 154 157 163 165 169 171 177 190 191 192 193 196 197 203
214 216 218 221 227 230 234 236 244 248 250 253 262 263 266 269
273
 Unital No   9
 |Aut(G,M65)|=  16
 17 18 21 22 37 46 47 48 49 62 63 64 67 68 69 70
 81 94 95 96 98 100 101 106 113 114 120 124 132 133 137 140
148 151 154 157 161 166 167 170 180 182 183 187 194 206 207 208
213 215 217 219 226 238 239 240 244 248 250 253 262 263 266 269
273
 Unital No  10
 |Aut(G,M65)|=   8
 26 27 28 29 33 36 40 41 49 62 63 64 66 78 79 80

```
  82  87  88  89 100 101 111 112 117 119 125 127 135 136 137 140
 146 149 150 159 161 164 170 175 179 180 183 188 197 198 201 203
 209 216 220 221 225 226 231 238 245 247 251 252 258 262 266 271
 273
 Unital No  11
 |Aut(G,M65)|=  32
  26  27  28  29  42  43  44  45  58  59  60  61  74  75  76  77
  85  87  88  90  98 100 101 106 115 121 122 123 130 135 136 141
 147 152 155 156 163 165 169 171 178 184 185 186 194 206 207 208
 213 215 217 219 227 230 234 236 242 249 251 252 257 259 261 265
 273
 Unital No  12
 |Aut(G,M65)|=  20
   1   3   5   7   9  17  18  25  30  32  36  39  40  46  48  50
  53  54  62  64  68  70  73  78  80  82  83  90  91  92 102 113
 116 118 124 128 132 134 135 137 139 146 147 148 149 154 162 163
 168 169 176 180 185 186 187 192 208 212 225 228 232 233 234 241
 267
 Unital No  13
 |Aut(G,M65)|=  32
  17  18  21  22  33  36  40  41  50  51  53  55  67  68  69  70
  81  94  95  96  97  99 103 108 115 121 122 123 129 142 143 144
 150 158 159 160 163 165 169 171 178 184 185 186 194 206 207 208
 214 216 218 221 226 238 239 240 244 248 250 253 262 263 266 269
 273
 Unital No  14
 |Aut(G,M65)|=  4
  10  11  12  13  18  20  26  30  34  40  42  48  50  56  59  64
  68  72  77  78  84  87  93  95 101 103 109 112 133 138 141 143
 145 151 156 159 162 168 172 173 180 182 183 187 194 206 207 208
 214 219 220 224 229 230 237 239 242 247 253 254 264 265 269 270
 273
 Unital No  15
 |Aut(G,M65)|=  4
   2   4   6   8  15  20  25  26  29  32  34  39  42  43  46  54
  56  58  59  62  72  73  74  77  80  81  97  98 100 108 110 119
 129 150 161 162 168 170 174 179 180 181 187 191 195 200 201 202
 208 213 215 218 221 222 226 244 248 251 252 256 257 261 262 263
 272
 Unital No  16
 |Aut(G,M65)|=  10
   4   7   9  11  13  19  20  26  29  32  33  39  40  42  44  50
  52  55  59  61  65  74  77  78  79  86  88  89  90  91  97 100
 101 104 112 115 122 124 127 128 131 134 137 140 144 160 169 177
 179 180 183 190 201 209 211 212 214 224 232 255 259 262 263 265
 271
 Unital No  17
 |Aut(G,M65)|=  8
  17  18  21  22  34  35  38  39  58  59  60  61  65  71  72  73
  84  85  87  89  97  99 105 109 114 120 121 122 130 133 135 137
 149 150 153 160 161 164 166 174 179 180 181 182 194 200 205 206
 211 217 219 223 231 233 238 240 243 244 253 256 257 260 264 265
 273
 Unital No  18
 |Aut(G,M65)|=  8
  19  30  31  32  37  46  47  48  58  59  60  61  67  68  69  70
  85  86  88  93  97 108 110 111 113 120 126 127 130 136 138 139
 148 150 154 158 169 171 174 175 180 184 185 187 193 200 203 204
 212 215 219 220 228 235 238 239 241 244 247 248 257 259 271 272
 273
 Unital No  19
 |Aut(G,M65)|=  8
  17  18  21  22  33  36  40  41  50  51  53  55  67  68  69  70
  82  84  85  88  97 102 108 109 113 120 122 123 130 132 136 137
 146 150 153 158 166 170 174 175 179 180 187 189 200 204 206 207
 211 215 219 224 231 237 238 239 243 244 248 254 257 259 264 267
 273
 Unital No  20
 |Aut(G,M65)|=  8
  19  30  31  32  42  43  44  45  52  54  56  57  74  75  76  77
  81  83  86  95 101 106 110 112 117 118 127 128 130 132 135 140
 148 155 156 157 161 164 166 174 177 184 186 190 195 204 207 208
 211 213 215 223 225 237 238 240 241 245 249 252 259 260 261 264
 273
 Unital No  21
 |Aut(G,M65)|=  8
```

```
 16  18  45  58  69  81  83  85  86  90 101 106 107 109 110 113
114 115 126 127 129 132 136 137 144 148 152 154 158 159 163 164
171 172 176 179 182 188 190 191 193 194 198 201 208 210 212 213
215 223 225 227 230 237 240 242 246 252 254 255 258 261 265 268
272
 Unital No  22
 |Aut(G,M65)|=   8
 16  17  42  61  67  83  86  87  88  94 101 102 104 105 106 120
122 124 126 127 129 132 135 137 144 148 154 156 158 159 162 164
165 169 176 181 187 189 190 191 194 196 199 202 208 209 213 215
220 224 225 234 236 237 239 241 242 252 254 255 258 261 264 265
272
 Unital No  23
 |Aut(G,M65)|=   4
  8  21  40  64  65  81  84  91  93  94  97 103 105 109 110 120
122 123 124 125 130 133 135 138 140 147 151 154 156 158 164 166
167 171 176 178 186 187 190 192 194 197 198 202 203 209 213 214
220 221 225 226 231 234 239 242 244 248 254 255 260 263 266 271
272
 Unital No  24
 |Aut(G,M65)|=  16
  2   3   5   6  20  23  24  25  34  35  38  39  67  68  69  70
 83  84  85  94  98 103 107 110 118 120 122 128 130 134 140 143
146 148 156 160 163 168 170 175 177 186 187 189 197 200 201 206
212 214 217 223 229 230 237 239 245 248 249 256 261 263 264 271
273
 Unital No  25
 |Aut(G,M65)|=   4
 16  26  41  51  77  82  87  88  92  93 100 103 107 108 109 114
115 122 124 128 130 131 133 135 140 148 153 154 155 156 163 164
165 168 173 179 182 187 189 191 197 198 202 203 206 210 212 213
216 218 229 234 235 236 237 242 246 247 250 252 257 259 266 267
268
 Unital No  26
 |Aut(G,M65)|=   4
 16  29  36  53  74  86  88  89  90  92  97 102 106 107 108 113
114 121 122 128 130 136 137 139 140 147 148 149 151 156 162 164
165 169 173 179 180 182 188 191 193 197 199 202 206 210 214 216
217 220 229 230 231 232 236 241 242 247 248 249 257 260 261 263
268
pp16\pp-16-19
 BRRESH1
 |Aut(G,M65)|=   8
  5   8   9  16  17  23  25  32  34  42  43  44  49  52  57  60
 65  68  77  80  83  97  99 100 103 116 120 124 128 130 138 142
143 148 153 154 157 162 163 165 167 169 171 173 176 195 202 203
204 210 211 222 223 233 234 236 240 243 249 250 256 258 259 260
261
 BRRESH2
 |Aut(G,M65)|=  16
  6  13  16  17  23  25  26  27  29  32  37  38  41  42  44  46
 47  49  53  54  55  57  62  63  69  72  78  83  85  89  90 107
109 110 120 121 125 133 137 140 150 157 159 165 166 171 186 189
191 193 198 201 209 223 224 227 235 238 248 249 255 260 263 268
270
 BRRESH3
 |Aut(G,M65)|=  16
  8  14  15  17  23  25  29  30  31  32  38  46  48  56  58  61
 66  74  80  81  83  90  95 101 107 109 119 122 126 130 139 143
152 153 160 162 166 169 177 183 189 193 194 201 213 217 222 225
227 232 234 236 238 239 242 246 247 250 251 252 253 258 262 269
273
 BRRESH4
 |Aut(G,M65)|=  16
  1   2   9  11  17  23  25  32  34  37  42  43  52  54  57  64
 66  69  77  78  81  83  84  85  93  97  98  99 100 116 119 125
126 134 138 142 144 147 150 153 154 163 169 171 176 195 202 203
208 210 215 218 221 225 233 236 240 246 247 249 250 257 262 263
265
 BRRESH5
 |Aut(G,M65)|=   8
  1   7   9  12  19  22  24  28  36  37  40  41  43  46  49  51
 52  56  59  64  66  70  72  79  82  83  84  87  93  98 112 115
116 123 124 136 138 146 156 162 164 165 166 183 186 188 189 205
208 214 216 217 219 221 222 223 224 225 232 241 254 258 263 266
271
```

```
 BRRESH6
|Aut(G,M65)|=  4
  1  5  7 17 18 22 23 25 31 32 36 38 39 40 41 43
 46 49 50 52 55 56 59 64 72 76 79 82 83 85 95 100
110 112 116 119 123 132 134 136 146 151 158 164 173 176 184 189
191 194 199 205 209 214 222 229 231 232 241 250 256 258 263 266
271
 BRRESH7
|Aut(G,M65)|=  4
  1  3  4  6  9 14 15 18 21 28 35 36 38 53 58 62
 69 76 78 85 86 95 104 108 110 113 117 118 119 122 126 127
131 132 144 147 156 159 161 170 177 181 182 190 197 198 207 216
218 220 225 227 228 229 230 231 237 251 252 254 257 265 266 268
270
 BRRESH8
|Aut(G,M65)|= 32
  1  3  5  7  8  9 11 12 14 16 22 27 28 32 47 50
 59 66 67 68 69 71 72 73 75 76 77 82 87 93 103 115
116 119 120 121 128 129 130 150 157 161 170 173 179 182 183 208
209 210 212 213 217 221 232 244 246 258 261 263 265 267 269 270
272
 BRRESH9
|Aut(G,M65)|= 32
  7  9 20 22 27 28 30 32 47 50 51 52 53 55 58 59
 62 73 77 82 86 87 89 93 103 115 116 129 130 132 141 150
154 157 158 161 163 170 171 173 177 179 182 183 190 194 195 198
206 208 210 212 229 230 232 235 236 241 244 246 249 250 252 253
254
 BRRESH   10
|Aut(G,M65)|= 20
 13 18 21 27 28 31 34 41 43 44 45 49 53 58 59 62
 74 86 93 95 100 102 108 109 110 123 129 131 132 135 143 148
150 151 156 159 161 162 167 170 172 174 177 178 182 187 190 192
193 198 199 204 208 210 226 229 230 232 240 246 251 252 254 256
271
 BRRESH   11
|Aut(G,M65)|=  4
  1  4  6  9 11 12 15 21 24 31 34 38 39 40 41 43
 44 60 61 63 65 66 69 91 94 100 105 111 113 117 119 135
137 139 148 151 152 162 165 168 175 181 184 196 200 205 209 211
216 217 220 221 222 227 234 240 243 255 256 260 263 264 265 266
267
 BRRESH   12
|Aut(G,M65)|=  8
  1  5  9 12 15 21 24 27 32 34 43 49 55 56 57 59
 64 68 75 76 77 78 81 84 96 97 100 116 118 120 122 125
130 135 137 141 142 144 155 158 162 164 173 184 186 188 202 203
210 213 215 222 224 233 240 242 252 257 258 259 263 266 269 270
271
 BRRESH   13
|Aut(G,M65)|= 16
  4 10 18 19 24 30 32 35 45 51 52 55 58 60 61 62
 63 70 79 81 85 86 95 104 111 113 117 122 123 124 125 132
134 140 142 145 147 152 158 162 165 171 176 184 185 190 191 194
195 198 199 200 206 214 216 228 231 242 248 250 256 257 260 266
273
 BRRESH   14
|Aut(G,M65)|=  8
  1 12 13 14 15 16 19 21 27 28 32 37 38 45 47 51
 52 57 58 62 63 78 80 84 85 87 90 98 104 106 107 119
120 132 141 145 155 165 166 170 173 179 185 186 188 196 200 203
208 209 213 225 229 230 231 233 235 236 239 250 252 258 265 267
270
 BRRESH   15
|Aut(G,M65)|= 16
  2  5  8 15 19 21 24 29 32 33 39 42 43 46 48 60
 61 66 67 68 72 74 75 76 78 81 85 90 96 97 109 114
119 121 127 134 140 146 147 152 153 156 159 162 164 165 166 184
185 188 189 201 202 209 211 221 223 227 240 248 256 258 265 267
270
 BRRESH   16
|Aut(G,M65)|= 16
 18 21 29 30 32 33 39 42 43 44 46 47 48 55 56 59
 63 66 68 72 75 86 90 95 96 97 103 106 109 130 132 142
144 145 146 148 150 153 156 158 159 164 166 171 176 188 189 190
191 201 202 205 208 225 227 232 240 242 244 250 255 258 265 267
```

270
 BRRESH    17
 |Aut(G,M65)|=  16
  2  4  5  8 10 15 32 35 44 45 47 51 52 56 58 59
 62 67 70 74 76 78 79 103 104 106 111 113 114 117 119 121
122 123 124 125 127 130 144 148 150 194 195 198 199 200 205 206
208 209 211 214 216 221 223 225 228 231 232 244 255 257 260 266
273
 BRRESH    18
 |Aut(G,M65)|=  8
  2  5  8 15 19 25 34 35 37 40 41 44 47 52 53 59
 64 67 70 73 74 76 77 78 79 85 100 102 111 114 119 121
127 130 131 133 137 145 148 149 156 167 183 193 199 201 209 210
211 212 214 216 221 223 225 228 232 234 238 239 240 247 250 253
255
 BRRESH    19
 |Aut(G,M65)|=  8
  7  9 19 22 25 29 40 52 53 59 60 61 64 70 79 85
 87 90 99 101 102 108 110 115 116 129 131 133 135 143 144 145
147 150 152 156 160 161 164 167 179 183 188 201 214 216 229 230
235 236 240 245 246 247 250 253 255 257 258 260 265 266 267 270
273
 BRRESH    20
 |Aut(G,M65)|=  8
  2  8 18 20 22 25 27 29 34 43 45 51 53 59 60 61
 64 67 74 84 87 89 90 95 100 104 109 126 128 129 131 135
139 143 144 145 146 147 150 152 160 161 162 164 169 171 179 180
184 188 190 193 197 200 217 219 231 237 238 245 246 247 249 250
255
 BRRESH    21
 |Aut(G,M65)|=  10
  6  8  9 10 18 20 28 29 36 37 43 47 57 59 61 62
 68 75 76 78 82 94 96 102 103 105 106 113 115 119 120 121
122 123 128 131 132 137 144 150 154 158 160 162 163 164 168 175
181 190 191 210 213 217 222 233 234 235 237 259 261 262 264 272
273
 BRRESH    22
 |Aut(G,M65)|=  8
  1  4  6 11 12 14 16 17 18 23 25 29 31 32 44 45
 46 49 54 57 69 77 78 83 87 89 93 97 98 99 100 101
108 110 113 119 121 130 140 141 148 157 158 164 165 169 179 185
190 196 200 207 212 215 223 226 232 238 242 243 248 261 263 269
272
 BRRESH    23
 |Aut(G,M65)|=  8
  4  8 11 17 23 25 27 29 30 32 33 34 39 40 42 44
 48 49 57 61 74 77 78 83 84 87 91 98 102 104 113 119
126 130 141 143 148 152 158 162 163 164 179 184 191 193 196 201
212 215 217 218 220 222 224 231 232 240 242 243 246 261 263 269
272
 BRRESH    24
 |Aut(G,M65)|=  16
  1  6  8 14 23 26 36 37 39 41 44 47 49 51 52 53
 54 55 57 63 64 68 69 72 74 95 101 109 123 125 127 128
129 130 136 138 144 145 148 150 152 153 154 155 158 159 163 190
195 196 204 205 207 208 219 220 221 222 235 240 243 247 249 253
256
 BRRESH    25
 |Aut(G,M65)|=  4
  2  9 13 16 17 23 25 32 35 36 39 44 49 54 62 63
 72 74 79 80 83 98 101 107 109 117 118 125 127 129 137 138
142 145 146 150 152 179 182 183 184 185 186 188 189 195 196 199
207 212 215 223 224 232 233 235 240 242 253 255 256 258 259 260
261
 BRRESH    26
 |Aut(G,M65)|=  4
  1  3  9 11 17 21 22 23 25 28 29 32 33 35 44 46
 59 60 62 63 69 71 75 79 83 98 102 107 110 117 120 124
128 138 140 142 144 145 146 149 151 162 165 167 173 196 199 202
206 212 213 217 220 232 233 236 237 242 243 246 253 258 259 260
261
pp16\pp-16-20
 Unital No   1
 |Aut(G,M65)|=  16
  4  7  9 14 18 21 27 32 34 37 43 48 50 56 58 64
 66 67 77 80 83 85 91 93 102 106 111 112 113 115 116 120

131 134 139 143 145 150 156 159 165 168 170 171 193 199 202 205
210 214 221 223 225 234 237 238 241 243 248 249 261 262 264 271
273
 Unital No  2
 |Aut(G,M65)|=  128
  7  8  9 10 18 21 28 31 35 37 43 45 49 55 57 63
 65 66 67 68 83 85 91 93 97 98 101 110 116 120 125 127
139 142 143 144 147 150 154 159 161 165 170 174 199 202 203 206
212 215 218 219 225 227 228 234 242 247 253 254 258 260 263 269
273
 Unital No  3
 |Aut(G,M65)|=  64
  7  8  9 10 17 22 27 32 34 40 42 48 50 56 58 64
 65 66 67 68 83 85 91 93 102 105 106 109 116 120 125 127
139 142 143 144 147 150 154 159 164 168 171 175 194 195 198 207
209 214 221 224 225 227 228 234 241 248 251 252 258 260 263 269
273
 Unital No  4
 |Aut(G,M65)|=  128
  1  2 15 16 17 22 27 32 35 37 43 45 49 55 57 63
 69 70 71 72 82 88 90 96 102 105 106 109 113 115 122 126
129 130 132 133 148 149 153 160 164 168 171 175 194 195 198 207
212 215 218 219 232 237 238 239 242 247 253 254 259 265 270 272
273
 Unital No  5
 |Aut(G,M65)|=  16
  7  8 13 14 19 23 25 29 34 35 42 43 57 58 63 64
 70 72 73 75 83 87 89 93 100 103 104 111 116 118 121 125
138 139 141 142 147 152 156 159 164 168 170 174 195 198 199 206
209 212 214 219 229 230 238 239 248 251 253 254 260 263 266 271
273
 Unital No  6
 |Aut(G,M65)|=  16
  9 10 13 14 19 20 25 26 33 38 41 46 51 53 60 62
 65 67 77 79 81 83 93 95 98 101 105 106 114 116 123 125
138 141 143 144 146 148 158 160 161 164 165 168 193 196 201 208
210 211 213 222 231 238 239 240 245 249 250 256 258 266 269 271
273
 Unital No  7
 |Aut(G,M65)|=  64
  7  8  9 10 20 23 26 29 36 38 44 46 50 56 58 64
 69 70 71 72 84 86 92 94 100 107 111 112 117 119 124 128
139 142 143 144 147 150 154 159 164 168 171 175 199 202 203 206
213 216 217 222 229 230 231 240 241 248 251 252 261 266 268 271
273
 Unital No  8
 |Aut(G,M65)|=  128
  1  2 15 16 20 23 26 29 36 38 44 46 50 56 58 64
 77 78 79 80 81 87 89 95 102 105 106 109 114 118 121 123
129 130 132 133 146 151 155 158 164 168 171 175 193 204 205 208
213 216 217 222 229 230 231 240 241 248 251 252 261 266 268 271
273
 Unital No  9
 |Aut(G,M65)|=  128
  3  4 13 14 20 23 26 29 34 40 42 48 52 54 60 62
 65 66 67 68 84 86 92 94 97 98 101 110 117 119 124 128
137 138 140 141 146 151 155 158 161 165 170 174 193 204 205 208
209 214 221 224 229 230 231 240 245 246 250 255 258 260 263 269
273
 Unital No  10
 |Aut(G,M65)|=  128
  3  4 13 14 17 22 27 32 36 38 44 46 50 56 58 64
 77 78 79 80 82 88 90 96 102 105 106 109 113 115 122 126
137 138 140 141 146 151 155 158 162 166 169 173 193 204 205 208
213 216 217 222 226 233 235 236 241 248 251 252 257 262 264 267
273
 Unital No  11
 |Aut(G,M65)|=  128
  1  5 11 15 19 20 23 24 33 36 46 47 49 56 57 64
 66 67 77 80 83 85 89 95 101 105 109 110 114 116 120 121
139 140 141 143 149 151 158 160 165 168 170 171 195 199 202 207
209 215 219 224 225 226 234 236 242 251 252 254 260 266 268 269
273
 Unital No  12
 |Aut(G,M65)|=  16
  1  7  9 15 18 19 21 24 35 38 44 45 51 54 59 62

```
  73  74  75  76  81  87  91  93 102 104 106 108 114 116 120 121
 130 132 142 144 151 152 157 158 161 165 170 174 193 204 205 208
 213 215 217 219 225 227 228 234 244 245 246 256 258 263 265 272
 273
```
 Unital No 13
 |Aut(G,M65)|= 32
```
  3  5 11 13 18 21 26 29 34 39 41 48 54 56 62 64
 65 68 73 76 81 87 92 94 101 104 108 110 114 117 121 124
 135 136 137 138 147 151 154 158 165 166 169 170 203 204 206 208
 209 210 223 224 225 230 234 240 241 245 246 248 258 262 263 264
 273
```
 Unital No 14
 |Aut(G,M65)|= 64
```
  1  5 11 15 25 26 29 30 34 35 45 48 51 54 59 62
 70 71 73 76 81 87 91 93 101 105 109 110 114 116 120 121
 130 132 135 136 145 147 154 156 165 168 170 171 196 197 204 208
 209 215 219 224 225 226 234 236 244 245 246 256 259 262 264 270
 273
```
 Unital No 15
 |Aut(G,M65)|= 16
```
  5  8 10 11 18 21 25 30 34 36 46 48 52 56 60 64
 69 72 77 80 82 83 85 88 99 103 105 109 116 117 120 124
 130 132 139 143 146 152 155 157 165 166 169 170 203 204 206 208
 209 211 220 224 233 235 237 238 243 249 251 252 259 260 269 270
 273
```
 Unital No 16
 |Aut(G,M65)|= 16
```
  4  8 10 14 19 20 23 24 38 39 41 44 50 55 58 63
 70 71 73 76 83 85 89 95 97 98 102 106 113 115 117 124
 129 131 133 134 146 148 153 155 161 164 174 175 196 197 204 208
 212 214 218 221 225 226 234 236 243 249 250 255 258 261 263 271
 273
```
 pp16\pp-16-21
 Unital No 1
 |Aut(G,M65)|= 16
```
  1 21 23 28 30 31 34 35 39 41 43 49 50 54 55 64
 68 73 76 79 80 92 97 98 102 110 112 117 119 121 122 128
 129 133 134 137 139 146 148 155 157 158 167 171 172 174 175 191
 196 197 200 206 208 209 212 213 214 219 230 242 244 249 252 255
 268
```
 Unital No 2
 |Aut(G,M65)|= 16
```
  7  9 11 14 15 17 20 21 26 30 35 41 42 44 46 51
 52 56 57 64 69 72 74 77 79 83 85 89 93 95 107 127
 130 135 138 141 142 148 150 151 152 157 164 168 171 174 175 177
 180 181 183 185 193 213 225 227 231 235 237 241 243 248 250 251
 268
```
 Unital No 3
 |Aut(G,M65)|= 128
```
  1  4  8 11 12 20 23 25 29 31 45 51 60 61 62 63
 66 67 71 75 76 83 88 90 94 96 110 116 123 125 126 128
 137 151 153 154 155 160 164 165 168 175 176 179 184 185 187 189
 202 216 217 218 220 223 227 230 231 239 240 244 247 250 252 254
 268
```
 Unital No 4
 |Aut(G,M65)|= 128
```
  1  4  5  7 10 17 18 23 25 30 34 51 55 61 62 63
 66 67 70 72 73 81 82 88 90 93 97 116 120 125 126 128
 134 147 151 153 154 155 161 163 165 168 174 179 181 182 186 189
 197 212 216 217 218 220 226 228 230 231 237 244 245 246 249 254
 268
```
 Unital No 5
 |Aut(G,M65)|= 16
```
  6  7 11 12 14 18 21 29 30 32 33 42 43 44 46 50
 65 71 74 79 80 82 84 85 87 91 101 114 115 119 124 127
 133 138 142 143 144 150 162 163 170 175 176 177 182 185 186 188
 193 211 214 215 219 224 227 229 235 236 238 241 243 246 248 255
 268
```
 Unital No 6
 |Aut(G,M65)|= 16
```
  3  4  5  8 16 20 21 23 27 29 37 51 55 57 58 60
 66 72 73 78 80 85 89 90 93 96 97 99 100 106 110 128
 134 135 136 137 141 155 163 165 169 171 174 178 186 187 189 190
 194 212 216 221 222 223 226 227 231 232 235 242 244 247 250 256
 268
```
 Unital No 7

|Aut(G,M65)|=  64
 4  7 11 12 14 19 25 28 30 32 38 39 40 42 46 59
67 72 75 76 77 84 90 91 93 95 101 103 104 105 109 124
130 131 132 138 142 159 163 168 170 175 176 183 186 188 189 192
193 195 196 201 205 224 228 231 233 239 240 248 249 251 254 255
268
 Unital No   8
|Aut(G,M65)|=  32
 3  4  8  9 11 18 20 21 24 32 34 36 43 44 47 55
57 59 62 63 65 66 72 78 79 96 101 104 105 106 111 119
130 131 140 141 142 148 163 165 166 169 172 187 197 199 204 207
208 212 217 220 222 224 227 231 232 238 240 242 243 245 247 251
268
 Unital No   9
|Aut(G,M65)|=  128
 4  7  9 13 16 18 20 24 31 32 36 41 42 44 47 62
67 72 74 78 79 81 83 87 95 96 99 105 106 107 112 125
136 139 141 142 144 154 163 168 169 172 173 180 182 184 187 188
199 204 205 206 207 217 228 231 234 235 238 243 245 247 251 252
268
 Unital No  10
|Aut(G,M65)|=  128
 5  7  9 13 14 19 23 24 29 31 36 41 44 45 47 50
51 53 55 58 66 67 73 79 80 83 85 86 89 92 100 122
135 157 165 168 171 172 174 177 178 184 190 191 199 202 204 206
207 210 212 213 216 221 226 228 233 234 238 243 244 248 250 252
268
 Unital No  11
|Aut(G,M65)|=  64
 2  4 11 15 16 18 19 23 25 31 33 52 53 54 55 63
65 67 76 79 80 81 84 88 90 96 98 115 117 118 120 128
133 145 146 147 152 155 166 168 171 172 175 179 182 183 187 189
198 209 210 212 215 220 229 231 235 236 240 244 245 248 252 254
268
 Unital No  12
|Aut(G,M65)|=  128
 1  3 10 11 14 19 21 22 26 29 38 51 53 55 56 57
66 68 73 76 77 84 85 86 89 94 101 116 118 119 120 122
130 145 147 148 151 157 165 167 170 174 175 177 178 183 185 190
193 210 211 212 216 222 230 232 233 237 240 241 242 248 250 253
268
 Unital No  13
|Aut(G,M65)|=  64
 4 17 18 19 21 31 33 39 41 44 48 52 53 55 63 64
70 74 75 76 80 89 101 102 105 110 112 113 116 117 124 126
129 130 137 139 142 146 150 151 153 159 161 171 173 175 176 190
196 198 203 206 207 210 212 215 219 220 231 242 245 246 248 252
268
 Unital No  14
|Aut(G,M65)|=  128
 2 20 21 25 27 31 38 39 43 47 48 55 57 58 60 63
66 68 70 73 80 91 97 99 100 102 110 114 117 118 126 127
129 134 135 136 137 145 146 149 153 156 161 165 167 171 174 192
193 196 203 204 208 212 220 221 222 223 229 242 247 252 254 256
268
 Unital No  15
|Aut(G,M65)|=  16
 1  4  7  9 13 17 25 27 30 31 35 39 40 45 47 54
60 61 62 63 68 81 83 88 89 95 98 99 102 106 109 116
129 133 136 138 141 151 167 179 182 184 188 190 195 196 200 202
204 209 217 218 220 223 228 230 231 234 238 246 249 252 254 256
268
 Unital No  16
|Aut(G,M65)|=  16
 3  7 11 13 16 27 36 49 51 52 60 63 67 70 71 75
77 83 84 86 91 95 98 104 109 110 111 113 118 122 125 127
131 133 137 138 140 145 150 154 156 157 161 164 168 170 176 177
183 184 188 192 199 214 215 216 220 223 228 232 234 235 240 256
268
pp16\pp-16-22
 Unital No   1
|Aut(G,M65)|=  16
 1  5  8 14 36 39 41 42 50 51 54 56 69 75 77 79
81 83 90 93 97 102 105 111 113 114 119 123 129 132 133 136
146 149 153 158 168 170 171 174 179 180 190 191 198 199 205 206
211 213 215 220 228 230 235 236 248 249 252 253 258 266 268 271

273
 Unital No  2
 |Aut(G,M65)|=  16
  1  5  8 14 34 35 38 40 52 55 57 58 65 76 78 80
 82 89 91 92 99 103 107 110 113 118 122 123 132 136 139 144
149 152 153 160 161 165 166 167 179 181 186 190 194 196 197 204
214 218 222 223 228 232 236 239 241 243 247 255 258 265 271 272
273
 Unital No  3
 |Aut(G,M65)|=  8
 14 17 22 24 27 32 33 38 40 45 48 49 50 51 63 64
 65 66 67 69 80 82 83 86 88 94 98 99 102 104 108 117
123 125 127 128 129 133 139 141 143 149 154 156 157 158 164 165
172 173 174 185 187 188 190 191 199 203 204 206 207 214 232 243
258
 Unital No  4
 |Aut(G,M65)|=  8
  4  9 10 16 21 30 31 32 37 46 47 48 59 61 62 64
 75 77 78 80 97 108 110 112 116 119 121 122 132 135 137 138
146 148 149 156 166 170 172 173 183 184 187 188 195 201 204 207
213 217 218 221 228 231 235 239 249 250 251 255 260 261 263 269
273
 Unital No  5
 |Aut(G,M65)|=  8
  5 17 19 21 27 30 33 34 45 46 47 49 56 61 62 63
 65 69 70 75 78 92 97 128 142 146 149 152 156 160 163 166
172 173 176 178 184 187 188 192 195 198 204 207 208 210 214 218
219 221 227 231 232 235 237 243 245 248 249 255 258 260 261 262
271
 Unital No  6
 |Aut(G,M65)|=  8
 17 27 30 31 33 37 45 46 49 59 62 63 65 69 77 78
 81 92 94 96 101 107 109 111 116 119 121 122 130 131 134 136
146 147 156 160 166 168 172 176 182 184 188 192 194 195 204 208
210 211 213 221 226 227 235 239 246 248 251 255 261 262 264 269
273
 Unital No  7
 |Aut(G,M65)|=  8
  1  5  8 14 20 24 27 32 38 41 45 48 51 58 63 64
 66 69 71 80 82 83 86 88 116 119 121 122 133 139 141 143
145 146 149 154 161 164 166 173 177 184 185 187 193 195 199 207
214 217 219 222 228 232 237 238 243 245 250 254 258 263 270 271
273
 Unital No  8
 |Aut(G,M65)|=  32
 18 19 22 24 34 35 38 40 49 60 62 64 65 76 78 80
 84 87 89 90 98 99 102 104 114 115 118 120 132 135 137 138
148 151 153 154 165 171 173 175 178 179 182 184 193 204 206 208
209 220 222 224 228 231 233 234 244 247 249 250 257 268 270 272
273
 Unital No  9
 |Aut(G,M65)|=  8
  1 21 27 29 31 32 37 43 45 47 48 52 55 57 58 64
 68 71 73 74 80 82 83 84 89 94 102 103 104 106 110 114
115 116 121 126 134 135 136 138 142 156 172 188 204 209 210 213
214 221 225 227 232 235 239 241 242 245 246 253 257 259 264 267
271
 Unital No  10
 |Aut(G,M65)|=  32
 17 28 30 32 33 44 46 48 53 59 61 63 69 75 77 79
 84 87 89 90 101 107 109 111 114 115 118 120 129 140 142 144
145 156 158 160 162 163 166 168 178 179 182 184 193 204 206 208
210 211 214 216 226 227 230 232 245 251 253 255 261 267 269 271
273
 Unital No  11
 |Aut(G,M65)|=  16
  2 17 21 23 25 29 33 34 35 36 44 50 56 59 60 64
 67 68 71 75 77 82 97 123 135 145 149 151 153 160 161 162
164 169 172 178 181 184 187 192 195 199 200 203 205 212 213 220
221 222 227 236 237 239 240 244 246 248 249 256 259 261 264 265
266
 Unital No  12
 |Aut(G,M65)|=  16
  2 17 19 22 27 30 35 37 41 42 47 49 52 58 59 63
 68 70 72 78 80 87 107 113 130 147 150 157 159 160 162 165
166 167 175 180 181 186 188 190 194 199 202 206 208 209 211 213

219 222 225 228 235 239 240 242 244 245 246 247 258 259 263 266 272

Unital No 13
|Aut(G,M65)|= 4
 2  3  9 10 12 26 37 49 53 57 58 59 66 68 73 75
 79 82 83 85 90 91 100 103 104 108 110 113 114 120 126 127
129 132 138 139 143 153 154 157 158 160 167 177 183 184 185 190
193 213 214 215 216 219 238 242 244 248 249 255 258 260 261 263
271

Unital No 14
|Aut(G,M65)|= 8
 3  7 13 15 17 20 23 30 33 41 42 46 50 51 53 61
 66 67 75 79 85 91 93 95 101 107 109 111 113 124 126 128
147 153 156 159 161 162 169 173 177 179 180 187 194 196 197 204
210 212 214 217 227 228 232 233 245 251 252 254 268 269 270 271
273

Unital No 15
|Aut(G,M65)|= 4
11 17 20 22 26 30 33 41 42 43 44 51 53 54 57 61
 66 67 75 77 78 91 106 120 142 147 148 149 150 153 161 168
171 173 176 179 180 184 186 188 193 197 202 203 204 212 213 214
217 221 225 227 232 233 238 245 246 248 252 254 260 263 268 269
271

Unital No 16
|Aut(G,M65)|= 4
 1  5  8 14 17 28 30 32 34 35 38 40 68 71 73 74
 81 88 90 93 100 102 109 112 115 121 124 125 130 135 141 142
150 153 156 159 167 168 174 175 177 178 186 191 195 196 207 208
209 213 215 216 229 230 233 240 242 245 250 254 259 260 261 268
273